%% file: main.tex
\documentclass[preprint,3p]{elsarticle}
\journal{Journal of Mechanics}
\input{preamble/preamble.tex}
\input{preamble/mathshortcuts.tex}

\biboptions{sort&compress}
\begin{document}
%
\begin{frontmatter}
	\title{Residual-based error estimation and adaptivity for stabilized immersed isogeometric analysis using truncated hierarchical B-splines}
	%
	\author[eindhoven,pavia]{Sai~C.~Divi \corref{mycorrespondingauthor}}
	\cortext[mycorrespondingauthor]{Corresponding author}
	\ead{s.c.divi@tue.nl}
	\author[eindhoven]{Pieter~H.~van~Zuijlen}
	\ead{pietervanzuijlen@outlook.com}
	\author[ut]{Tuong~Hoang}
	\ead{t.hoang-1@utwente.nl}
	\author[reden]{Frits de Prenter}
	\ead{f.deprenter@reden.nl}
	\author[pavia]{Ferdinando~Auricchio}
	\ead{auricchio@unipv.it}
	\author[pavia]{Alessandro~Reali}
	\ead{alessandro.reali@unipv.it}
	\author[eindhoven]{E.~Harald~van~Brummelen}
	\ead{e.h.v.brummelen@tue.nl}
	\author[eindhoven]{Clemens~V.~Verhoosel}
	\ead{c.v.verhoosel@tue.nl}
	\address[eindhoven]{Department of Mechanical Engineering, Eindhoven University of Technology, 5600MB Eindhoven, The Netherlands}
	\address[pavia]{Department of Civil Engineering and Architecture, University of Pavia, 27100 Pavia, Italy}
	\address[ut]{University of Twente, P.O.~Box 217, 7500 AE Enschede, The Netherlands}
	\address[reden]{REDEN -- Research Development Netherlands, F.~Hazemeijerstraat 800, 7555 RJ Hengelo, The Netherlands}
	%

\subfile{sections/abstract}

	\begin{keyword}
		{Isogeometric analysis, Immersed methods, Error estimation, Adaptivity, Stabilization, Scan-based analysis}
	\end{keyword}
\end{frontmatter}

%
%
\newpage

\subfile{sections/introduction}

\subfile{sections/fundamentals}
\subfile{sections/eea}

\subfile{sections/results}

\subfile{sections/conclusions}
\subfile{sections/acknowledgement}
\bibliography{bibliography/references}
%
\end{document}

%% file: preamble/preamble.tex
\usepackage{graphicx}
\graphicspath{{images/}}
\usepackage{subcaption}
\usepackage{float}
\usepackage{pdflscape}
\usepackage{rotating}
\usepackage[inline]{enumitem}

\usepackage{amssymb}
\usepackage{mathrsfs}
\usepackage{mathtools}
\usepackage{stmaryrd}
\usepackage{amsthm}
\newtheorem*{example*}{Example}

\newtheorem*{proof*}{Proof}
\newtheorem{remark*}{Remark}

\usepackage[lined,boxed,commentsnumbered,linesnumbered]{algorithm2e}

\SetCommentSty{mycommfont}
\usepackage{color}
\usepackage{subfiles}
\usepackage[hidelinks]{hyperref}
\usepackage{booktabs}
\usepackage{multirow}
\usepackage{tikz}
\usetikzlibrary{math}
\usetikzlibrary{spy}
\usetikzlibrary{calc}
\usetikzlibrary{decorations.markings}
\usetikzlibrary{shapes.misc, positioning}
\usepackage{epstopdf}

%% file: preamble/mathshortcuts.tex
\def\bm#1{\mathbf{#1}}
\newcommand{\R}{\mathbb{R}}			
\newcommand{\gradn}{\partial_n}
\renewcommand{\sup}{\operatornamewithlimits{sup\vphantom{p}}}

\newcommand{\domain}{\Omega}
\newcommand{\boundary}{\partial \domain}
\newcommand{\dirichletboundary}{\boundary_{D}}
\newcommand{\neumannboundary}{\boundary_{N}}
\newcommand{\ambientdomain}{\mathcal{A}}
\newcommand{\meshsize}{h}
\newcommand{\mesh}{\mathcal{T}}
\newcommand{\cutmesh}{\mathcal{T}_{\domain}}      
\newcommand{\ambientmesh}{\mesh_{\ambientdomain}} 
\newcommand{\boundarymesh}{\mesh_{\boundary}}

\newcommand{\dirichletmesh}{\mesh_{\dirichletboundary}}
\newcommand{\element}{K}
\newcommand{\face}{F}
\newcommand{\edge}{E}
\newcommand{\skeleton}{\mathcal{F}_{\rm skeleton}}
\newcommand{\ghost}{\mathcal{F}_{\rm ghost}}

\newcommand{\V}{V}		
\newcommand{\Vh}{{\V^{\meshsize}}}		
\newcommand{\Uh}{{U^{\meshsize}}}		
\renewcommand{\v}{v}		
\newcommand{\vh}{{\v^{\meshsize}}}		
\renewcommand{\u}{u}		
\newcommand{\uh}{{\u^{\meshsize}}}		
\newcommand{\B}{{a^{\meshsize}}}	
\newcommand{\F}{{b^{\meshsize}}}	

\newcommand{\bodyheat}{f}
\newcommand{\dirichletdataheat}{g}
\newcommand{\neumanndataheat}{q}

\newcommand{\bodystokes}{\bm{f}}
\newcommand{\neumanndatastokes}{\bm{t}}
\newcommand{\dirichletdatastokes}{\bm{g}}

\newcommand{\rintu}{\bm{r}_{{\rm int},\uu}^\meshsize}
\newcommand{\rintp}{r_{{\rm int},\p}^\meshsize}

\newcommand{\rjumpu}{\bm{r}_{\rm jump}^\meshsize}

\newcommand{\rghostu}{\bm{r}_{\rm ghost}^\meshsize}

\newcommand{\rnitscheu}{\bm{r}_{\rm nitsche}^\meshsize}
\newcommand{\rneumannu}{\bm{r}_{\rm neumann}^\meshsize}
\newcommand{\rskeleton}{r_{\rm skeleton}^\meshsize}

\newcommand{\supp}{\text{supp}}	
\newcommand{\hbasis}{\mathcal{H}} 
\newcommand{\bbasis}{\mathcal{B}} 
\newcommand{\basisfunc}{N} 
\newcommand{\order}{k} 
\newcommand{\regularity}{\alpha} 
\newcommand{\splinespace}{\mathcal{S}^{\order}_{\regularity}} 
\newcommand{\trunc}{\text{trunc}}	
\newcommand{\thblevel}{\ell}	
\newcommand{\thbmaxlevel}{{\thblevel_{\rm max}}}	

\newcommand{\Rh}{r^h}	

\newcommand{\uu}{\bm{u}}			
\newcommand{\uuh}{\uu^{\meshsize}}
\newcommand{\vv}{\bm{v}}			
\newcommand{\vvh}{\vv^{\meshsize}}
\newcommand{\p}{p}			
\newcommand{\ph}{\p^{\meshsize}}			
\newcommand{\q}{q}			
\newcommand{\qh}{\q^{\meshsize}}			
\newcommand{\nn}{\bm{n}}			
\newcommand{\ff}{\bm{f}}			
\newcommand{\g}{\bm{g}}				

\newcommand{\ltrivert}{\left\vert\kern-0.25ex\left\vert\kern-0.25ex\left\vert}
\newcommand{\rtrivert}{\right\vert\kern-0.25ex\right\vert\kern-0.25ex\right\vert}

\newcommand{\dO}{\, {\rm d}V}		
\newcommand{\dG}{\, {\rm d}S}		


%% file: sections/abstract.tex
\begin{abstract}

We propose an adaptive mesh refinement strategy for immersed isogeometric analysis, with application to steady heat conduction and viscous flow problems. The proposed strategy is based on residual-based error estimation, which has been tailored to the immersed setting by the incorporation of appropriately scaled stabilization and boundary terms. Element-wise error indicators are elaborated for the Laplace and Stokes problems, and a THB-spline-based local mesh refinement strategy is proposed. The error estimation and adaptivity procedure is applied to a series of benchmark problems, demonstrating the suitability of the technique for a range of smooth and non-smooth problems. The adaptivity strategy is also integrated in a scan-based analysis workflow,  capable of generating reliable, error-controlled, results from scan data, without the need for extensive user interactions or interventions.
\end{abstract}

%% file: sections/introduction.tex
\section{Introduction}
Immersed finite element methods, such as the finite cell method \cite{parvizian2007, duster2008, schillinger2015} and CutFEM \cite{hansbo2002, burman2012, burman2015}, are a natural companion to isogeometric analysis \cite{hughes2005, hughes2009}. The combination of immersed methods with the spline-based discretization strategy provided by the isogeometric analysis paradigm is recognized as a valuable extension of isogeometric analysis, because the immersed analysis concept provides a cogent framework for the consideration of trimmed CAD objects \cite{schmidt2012, rank2012, schillinger2012, ruess2013, ruess2014, marussig2018}. Moreover, immersed isogeometric analysis enables the construction of spline-based discretization spaces for geometrically and topologically complex volumetric domains \cite{rank2012, schillinger2011, schillinger2012}, a simulation strategy referred to as immersogeometric analysis \cite{hsu2015, kamensky2015}.

In comparison to boundary-fitting isogeometric analysis, the immersed isogeometric analysis strategy requires consideration of three (categories of) non-standard computational aspects. First, the geometry of elements that intersect with the boundary of the computational domain must be resolved by a dedicated integration procedure; see \emph{e.g.}, \cite{kudela2015, kudela2016, joulaian2016, abedian2019, divi2020, antolin2021}. Second, Dirichlet boundary conditions on immersed boundaries can generally not be imposed through basis function constraints. Instead, such boundary conditions are frequently imposed weakly; see, \emph{e.g.}, \cite{nitsche1971, hansbo2002, bazilevs2007weak, embar2010}. Third, unfavorably trimmed elements are notorious for causing ill-conditioning problems and, along Dirichlet boundaries, large unphysical gradients \cite{burman2010, ruess2013, massing2014, schillinger2015, dettmer2016, deprenter2017,deprenter2018}. This  problem is amplified in the higher-order discretization setting of isogeometric analysis \cite{deprenter2017}. Prominent computational remedies to overcome these problems are to supplement the weak formulation with stabilization terms, see, \emph{e.g.}, \cite{burman2010, burman2012, burman2015}, or to constrain, extend, or aggregate basis functions, see \emph{e.g.}, \cite{hollig2001, hollig2005, ruberg2012, ruberg2014, marussig2017, badia2018, badia2018mixed, marussig2018}, or to apply dedicated preconditioning techniques, \emph{e.g.}, \cite{deprenter2017,deprenter2019,jomo2019}.
	
For mixed formulations, such as standard weak forms of the Stokes and Navier-Stokes equations, the immersed isogeometric analysis setting imposes an additional challenge. In order to satisfy the inf-sup condition \cite{babuvska1973, brezzi1974} in boundary-fitting (isogeometric) analyses, generally use is made of stable pairs of basis functions (\emph{e.g.}, Taylor-Hood  \cite{taylor1973, bazilevs2006, buffa2011, bressan2013} or Raviart-Thomas \cite{raviart1977, girault1979, buffa2011, evans2013raviart}). Alternatively, stabilization techniques such as GLS \cite{hughes1989, douglas1989, tezduyar1991}, VMS \cite{hughes1998, hughes2001, bazilevs2007} or projection methods \cite{becker2001, burman2006vms} can be used. Direct utilization of these elements or stabilization techniques in the immersed setting can lead to non-physical spurious oscillations in the solution, even with relatively large and regular cut element configurations \cite{hoang2017, hoang2019}. One remedy for tackling this issue is to employ a skeleton-stabilized immersed isogeometric technique \cite{hoang2019}. The fundamental idea of this stabilization technique is to penalize (high-order) pressure derivative jumps over the edges/faces of the background mesh, resulting in stable discretizations using equal-order spline spaces. The technique proposed in Ref.~\cite{hoang2019} is inspired by the (continuous) interior penalty ((C)IP) and the ghost penalty (GP) methods \cite{burman2010}, extending these techniques to the case of high-regularity isogeometric analysis.

An appraised property of immersed methods in general, and immersed isogeometric analysis in particular, is that the discretization resolution can be controlled independently of the geometry parametrization. The immersed analysis concept avoids the need for geometry-induced mesh refinements in the vicinity of geometric details that are irrelevant in relation to the objective of an analysis. This decoupling of the discretization resolution from the geometry makes it natural to consider immersed finite elements in combination with adaptive discretization strategies. In fact, adaptivity in the form of local $p$- and $hp-$refinements has always been an integral part of the finite cell method \cite{duster2017, zander2015, dangella2016, elhaddad2018}.	

A posteriori error estimation and adaptivity techniques are well-established in the context of finite element methods; see, \emph{e.g.}, the reviews \cite{bank1993, ainsworth1997, gratsch2005}. A variety of error estimation and adaptivity techniques has been studied in isogeometric analysis, such as residual-based error estimators for T-splines \cite{dorfel2010} and hierarchical splines \cite{vuong2011, giannelli2012, buffa2016}, and goal-oriented techniques \cite{kuru2014}. The contemporary overview \cite{bracco2019} is also noteworthy, as is the advanced industrial application considered in  Ref.~\cite{coradello2020}. In the context of Nitsche-based finite element methods (see Refs.~\cite{hansbo2005, chouly2017} for an overview), studies on a posteriori error estimators have been conducted \cite{hansbo2002, hansbo2003, becker2003, juntunen2009, chouly2018}. Local refinement strategies in immersed methods are predominantly feature based, \emph{i.e.}, either based on geometric features such as boundaries, or based on solution features such as sharp gradients in the solution fields; see, \emph{e.g.}, \cite{schillinger2012,  bandara2016, kanduvc2017} for examples of local refinement capabilities in finite cell simulations. Goal-oriented error estimation and adaptivity for immersed methods has also been studied \cite{kuru2014, verhoosel2015, distolfo2019, distolfo2019dual}. In the context of stabilized immersed finite elements, Ref.~\cite{burman2019} considered a posteriori element-wise error estimation and adaptivity to improve boundary approximations.
							
Although the computational setting of immersed isogeometric analysis enables the use of volumetric spline patches, the standard $h$, $p$ and $k$-type refinement strategies in patch-based isogeometric analysis \cite{hughes2005} are not suitable because of the non-local propagation of refinements. Various alternative refinement strategies have been proposed over the last decade to construct local spline refinements, the most prominent of which are T-splines \cite{sederberg2003, bazilevs2010, scott2011, scott2012, schillinger2012, buffa2014, hsu2015}, LRB-splines \cite{kvamsdal2014,kvamsdal2015}, U-splines \cite{thomas2018}, and (Truncated) Hierarchical B-splines \cite{brummelen2020}. In the context of immersed isogeometric analysis on volumetric domains, hierarchical splines are particularly suitable, as they optimally leverage the advantages offered by the geometrically simple background mesh.

In this contribution we propose a computational strategy for the application of residual-based a-posteriori error estimation and mesh adaptivity to stabilized immersed isogeometric analyses. We study various computational aspects of the framework that are non-standard in comparison to error estimation and adaptivity for boundary-fitting analyses, \emph{viz.}:
\begin{enumerate*}[label=\emph{(\roman*)}]
  \item In immersed analyses, the discretization basis is constructed over a mesh comprised of all elements in an ambient mesh that intersect with the computational domain. As a direct consequence of this setting, the support of the computational basis in general changes under refinement operations. The same holds for the mesh skeleton, which is a key ingredient of the considered stabilization methods. The considered computational strategy preserves the geometry of the computational domain under local mesh refinements, despite the change of the background mesh;
  \item Weak formulations in stabilized immersed isogeometric analysis generally involve operators with an explicit dependence on the mesh size. While this mesh size is unambiguously defined in the case of a uniform background mesh, the local mesh refinements considered in the adaptive setting warrant careful consideration of the scaling of the stabilization terms. We herein propose and study a scaling of the stabilization terms based on the local element sizes. 
\end{enumerate*}

We demonstrate the performance of the proposed computational strategy using a series of test cases for steady heat conduction problems (Poisson problem) and steady viscous flow problems (Stokes problem). We consider the application of the proposed adaptivity technique in a scan-based isogeometric analysis setting, and demonstrate that a robust automatic simulation workflow is realized when the methodology presented herein is combined with the topology-preserving image segmentation algorithm presented in Ref.~\cite{divi2020}.

This paper is outlined as follows. Section~\ref{sec:fundamentals} introduces the immersed isogeometric analysis framework, along with a detailed stability analysis for the considered model problems. This analysis focuses particularly on the scaling relations for the stabilization terms. In Section~\ref{sec:eeasection} the residual-based error estimator is introduced, and a mesh-refinement strategy is proposed. Benchmark simulation results are then presented in Section~\ref{sec:results} for both the steady heat conduction problem and the viscous flow problem, after which the developed framework is applied in a scan-based setting in Section~\ref{sec:scanbased}. Conclusions are finally drawn in Section~\ref{sec:conclusions}.

%% file: sections/fundamentals.tex
\section{Stabilized immersogeometric analysis with local mesh refinements}\label{sec:fundamentals}
In this section we introduce the stabilized immersed isogeometric analysis formulations for the steady heat conduction (Laplace) problem and steady viscous flow (Stokes) problem. We commence with presenting the  general setting of the problems in Section~\ref{sec:fcmsetting}, after which the stabilized formulations are presented in Section~\ref{section:immersogeometric}. In preparation of the \emph{a posteriori} error estimation concept discussed in Section~\ref{sec:eeasection}, in Section~\ref{sec:stability} we study the stability of the considered formulations.

\subsection{The finite cell setting}\label{sec:fcmsetting}
\begin{figure}
   \centering
   \begin{subfigure}[b]{\textwidth}
   \centering
   \includegraphics[width=.8\textwidth]{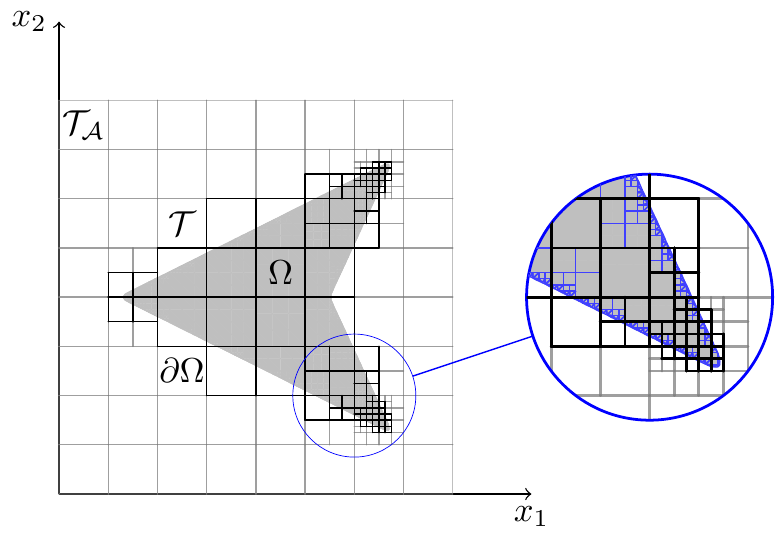}
   \caption{}\label{figure:fcmdomainsa}   
   \end{subfigure}\\[12pt]
   \begin{subfigure}[b]{0.45\textwidth}
   \centering
      \includegraphics[width=\textwidth]{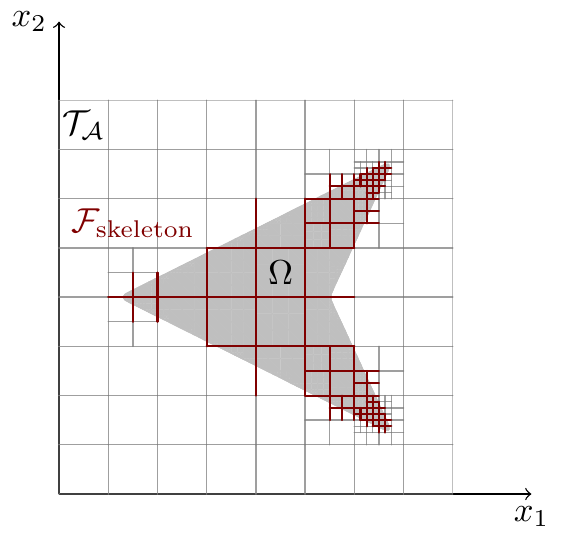}
   \caption{}\label{figure:fcmdomainsb}   
   \end{subfigure}
   \hfill
   \begin{subfigure}[b]{0.45\textwidth}
   \centering
   \includegraphics[width=\textwidth]{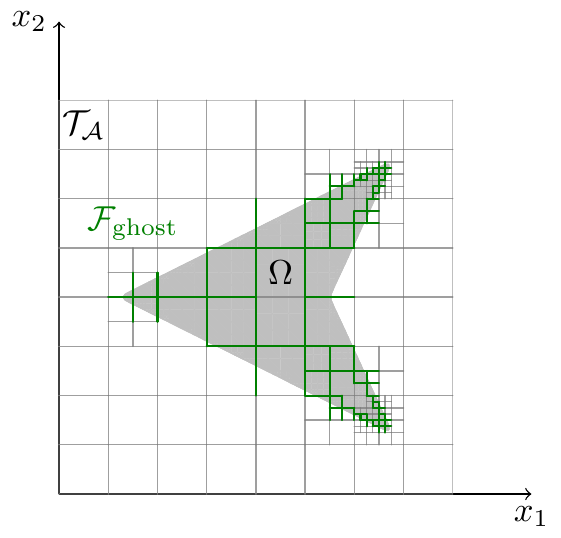}
   \caption{}\label{figure:fcmdomainsc}   
   \end{subfigure}
   
	 \caption{(\subref{figure:fcmdomainsa}) A physical domain $\domain$, with boundary $\boundary$, is embedded in an ambient domain $\ambientdomain$. The background mesh $\mesh$, which consists of all elements  that intersect the physical domain, is constructed by locally refining the ambient domain mesh $\ambientmesh$. The zoom illustrates the employed bisectioning procedure to capture the immersed boundaries. The integration subcells are marked in blue, whereas the background cells are marked in black. The skeleton mesh, $\mathcal{F}_{\rm skeleton}$, and ghost mesh, $\mathcal{F}_{\rm ghost}$, are shown in panels (\subref{figure:fcmdomainsb}) and (\subref{figure:fcmdomainsc}), respectively.}
     \label{figure:fcmdomains}
\end{figure}

We consider a physical domain $\domain \in \R^d$ (with $d\in \{2,3\}$) with boundary $\boundary$, as illustrated in Figure~\ref{figure:fcmdomains}. The boundary is composed of a Neumann part, $\neumannboundary$, and a Dirichlet part, $\dirichletboundary$, such that  $\overline{\neumannboundary} \cup \overline{\dirichletboundary} = \boundary$ and $\neumannboundary\cap\dirichletboundary=\emptyset$. The outward-pointing unit normal vector to the boundary is denoted by  $\nn$.

The physical domain is immersed in a geometrically simple ambient domain, \emph{i.e.}, $\ambientdomain \supset \domain$, on which a locally refined ambient mesh $\ambientmesh$ with elements $\element$ is defined. In this work, the ambient domain is chosen to be rectangular or cuboid, to facilitate simple, tensor-product, spline discretizations. The locally-refined meshes are constructed by sequential bisectioning of (selections of) elements in the mesh, starting from a Cartesian mesh. Truncated hierarchical B-splines can be formed on such meshes, as will be elaborated in Section~\ref{section:immersogeometric}.

Elements that do not intersect with the physical domain can be omitted from the ambient mesh, resulting in the locally refined (active) background mesh
\begin{equation}\label{equation:background}
\mesh := \{\element \, | \, \element \in \ambientmesh, \element \cap \domain \neq \emptyset \}.
\end{equation}
In the remainder, with the abuse of notation, we will use $\mesh$ (and other meshes) to denote both the set of elements in the mesh and the geometry obtained from the union of these elements. The local mesh size of the locally refined background mesh is denoted by 
\begin{align}
\meshsize_{\mathcal{T}}: \element \rightarrow \meshsize_{\element} = \sqrt[d]{\text{meas}(\element)}. \label{eq:meshsize}
\end{align}
By cutting the elements that are intersected by the immersed boundary $\boundary$, a mesh that conforms to the physical domain $\Omega$ is obtained:
\begin{equation}\label{equation:cutmesh}
\cutmesh := \{\element \cap \Omega \, | \, \element \in \mesh \}.
\end{equation}
The collection of elements in the background mesh that are crossed by the immersed boundary $\partial \Omega$ is defined as
\begin{equation}
	\mathcal{G} := \{ K \in \mesh \mid K \cap \partial \Omega \neq \emptyset \}.
\end{equation}

In immersed methods, the geometry of the physical domain is captured by the integration procedure on the cut elements, \emph{i.e.}, elements that are intersected by the immersed boundary $\boundary$. We herein employ an octree integration procedure \cite{verhoosel2015,divi2020}, which we close at the lowest level of bisectioning with a tessellation procedure. The considered integration procedure is illustrated in Figure~\ref{figure:fcmdomains} (in blue) for a typical cut element; see Ref.~\cite{divi2020} for further details. The employed tessellation provides an explicit parametrization of a polygonal approximation of the immersed boundary $\boundary$ through the set of boundary faces 
\begin{equation}\label{equation:boundaryedges}
\boundarymesh := \{ \edge \subset \boundary \, | \, \edge = \partial \element \cap \partial \domain,\,  \element \in \cutmesh \}.
\end{equation}
All faces $E \subset \partial \Omega_{N}$ (respectively $E \subset \partial \Omega_{D}$) are assigned to a set of Neumann faces $\mathcal{T}_{\partial \Omega_{N}}$ (respectively Dirichlet faces $\mathcal{T}_{\partial \Omega_{D}}$). In general, a single polygon face can overlap with both the Neumann and the Dirichlet boundary. Let us note that in an adaptive refinement procedure, the refinements can serve to provide an increasingly accurate approximation of the transition between the Neumann and Dirichlet boundary.

The formulations considered in the remainder of this work incorporate stabilization terms formulated on the edges of the background mesh (see Section \ref{section:immersogeometric}), which we refer to as the skeleton mesh
\begin{equation}\label{equation:skeletoninterfaces}
\skeleton = \{ \partial \element \cap \partial \element' \,|\ \element,\element'\in \mesh, \element\neq \element' \}.
\end{equation}
Note that the boundary of the background mesh is not part of the skeleton mesh. In addition to the skeleton mesh, we define the ghost mesh as the subset of the skeleton mesh that contain a face of an element intersected by the domain boundary
\begin{equation}\label{equation:ghostinterfaces}
\ghost = \{ \partial \element \cap \partial \element' \,|\, \element \in \mathcal{G}, \element'\in \mesh,  \element\neq \element' \}.
\end{equation}
As will be detailed in Section~\ref{sec:stability}, the stabilization terms formed on the skeleton and ghost mesh account for stability and ill-conditioning effects related to unfavorably cut elements, as well as for preventing pressure oscillations in equal-order discretizations of the Stokes problem.

\subsection{Immersogeometric analysis}\label{section:immersogeometric}
\begin{figure}
    \centering
    \includegraphics[width=\textwidth]{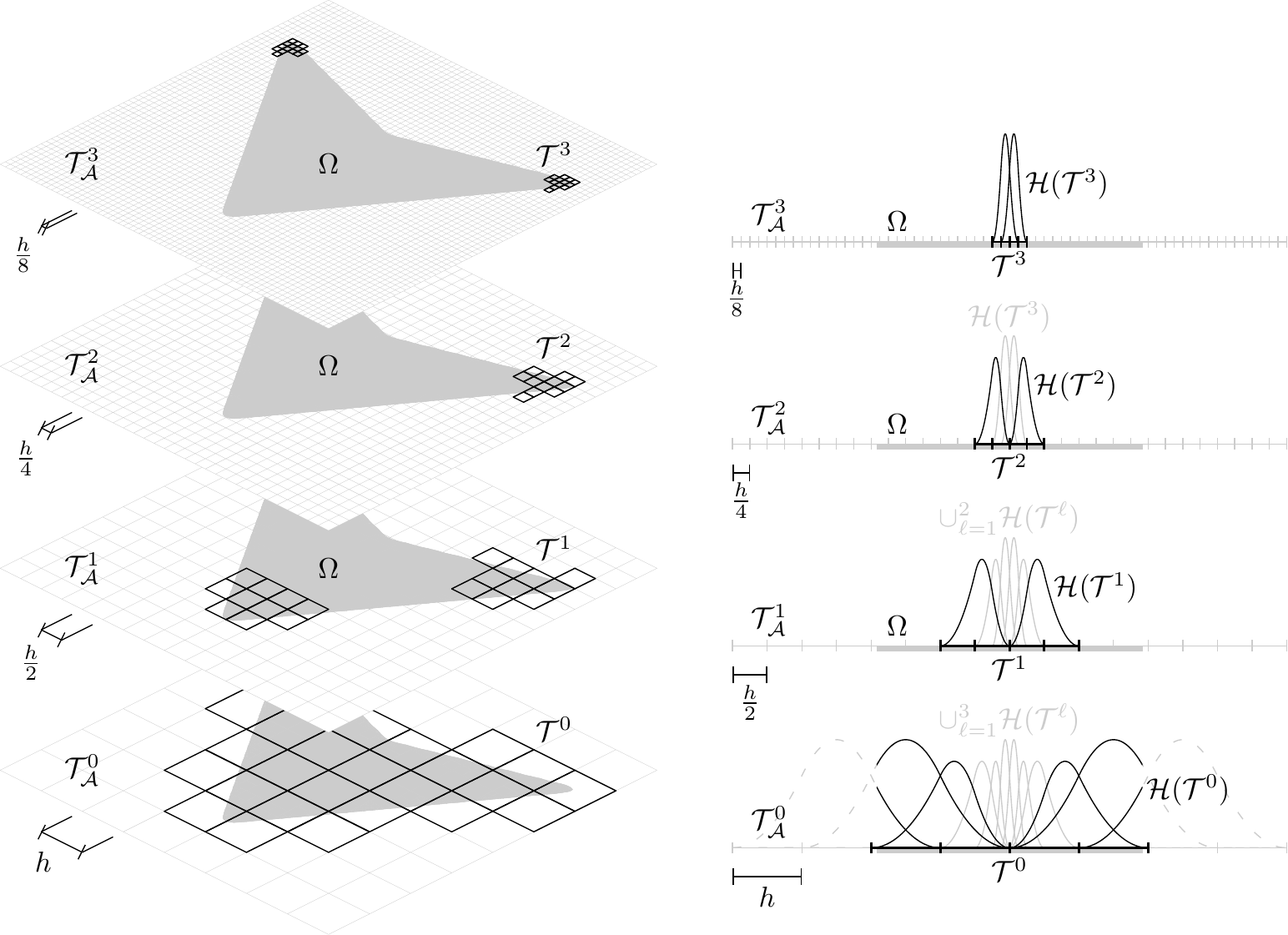}
	  \caption{Illustration of truncated hierarchical B-splines \cite{giannelli2012,brummelen2020} in the immersogeometric analysis setting. The left column shows the hierarchical levels of the mesh $\mesh$ in Figure~\ref{figure:fcmdomains}, while the right column illustrates the concept for a one-dimensional immersed domain $\domain$. The background mesh at the level $\thblevel=0,\cdots,\thbmaxlevel$ (with $\thbmaxlevel=3$ in this illustration) is defined as $\mesh^{\thblevel}=\{ \element \in \ambientmesh^{\thblevel}  \mid \element \cap \Omega \neq \emptyset\}$ where $\ambientmesh^{\thblevel}$ is a regular mesh with mesh size parameter $2^{-\ell} \meshsize$. Note that the meshes are nested, in the sense that the domain covered by the physical mesh at level $\thblevel$, $\mesh^{\thblevel}$, is completely inside that of level $\thblevel-1$, $\mesh^{\thblevel-1}$, \emph{i.e.}, $\mesh^{\thblevel} \subseteq \mesh^{\thblevel-1}$. The THB-spline basis, $\hbasis(\mesh)$, is constructed by selection and truncation of the basis functions in the B-spline basis $\bbasis(\mesh^{\thblevel}) = \{ \basisfunc \in \bbasis(\ambientmesh^{\thblevel}) \mid  \supp{(\basisfunc)} \subseteq \mesh^{\thblevel}  \}$ defined at each level. At the most refined level, \emph{i.e.}, at $\thblevel=\thbmaxlevel$, all basis functions that are completely inside $\mesh^{\thbmaxlevel}$ are selected: $\hbasis(\mesh^{\thbmaxlevel})=\{ \basisfunc \in \bbasis(\mesh^{\thbmaxlevel}) \mid \supp{( N )} \subseteq \mesh^{\thbmaxlevel} \}$. At the coarser levels, \emph{i.e.}, $\thblevel < \thbmaxlevel$, the functions that are completely inside the domain $\mesh^{\thblevel}$ but not completely inside the refined domain $\mesh^{\thblevel+1}$ are selected and truncated: $\hbasis(\mesh^{\thblevel})=\{\trunc{(N)}  \mid \basisfunc \in \bbasis(\mesh^{\thblevel}),\,\supp{( N )} \nsubseteq \mesh^{\thblevel+1}\}$. The truncation operation reduces the support of the B-spline functions by projecting away basis functions retained from the refined levels. The THB-spline basis then follows as $\hbasis(\mesh)=\cup_{\thblevel=0}^{\thbmaxlevel} \hbasis(\mesh^{\thblevel})$. The reader is referred to Ref.~\cite{giannelli2012} for details of THB-spline basis and Ref.~\cite{brummelen2020} for THB-spline basis construction.
}\label{figure:thbsplines}
\end{figure}

We consider the immersogeometric analysis of a single-field steady heat-conduction problem and of a two-field viscous flow problem. Both problems are represented by the abstract Galerkin problem
\begin{equation}\label{equation:galerkin}
\left\{ \begin{array}{l}
\text{Find } \uh \in \Uh\text{ such that:}\\
\B(\uh,\vh) = \F(\vh) \qquad \qquad \forall \vh \in \Vh,
\end{array}\right.
\end{equation}
with mesh-dependent bilinear and linear forms, $\B$ and $\F$, respectively. Note that the superscript $\meshsize$ is used to indicate mesh-dependence. The finite dimensional trial and test spaces, $\Uh$ and $\Vh$, are spanned by truncated hierarchical B-spline (THB-spline) \cite{giannelli2012,brummelen2020} basis functions of degree $\order$ and regularity $\regularity$ constructed over the locally-refined background mesh, \emph{viz.}
\begin{equation}
\splinespace(\mesh) 
= 
\{ \basisfunc\in{}C^{\regularity}(\mesh):\basisfunc|_{\element}\in{}P^{\order}(\element),\,\forall{}\element\in \mesh\},
\end{equation}
with $P^k(K)$ the set of $d$-variate polynomials on the element $K$ constructed by the tensor-product of univariate polynomials of order $\order$. Truncated hierarchical B-splines, which are illustrated in Figure~\ref{figure:thbsplines}, form a partition of unity and have a reduced support compared to their non-truncated counterpart, which is advantageous from the perspective of system matrix sparsity. Our implementation is based on the open source finite element library Nutils \cite{nutils}.

Since the imposition of strong Dirichlet boundary conditions over the immersed boundary $\boundary$ is intractable in the immersogeometric analysis setting, such boundary conditions are imposed weakly through Nitsche's method; see, \emph{e.g.}, Ref.~\cite{embar2010}. A mesh-dependent consistent stabilization term is introduced in order to ensure the well-posedness of the Galerkin problem \eqref{equation:galerkin}.

\subsubsection{Steady heat conduction}\label{section:laplace}
Steady heat conduction is governed by the Poisson problem, which, in dimensionless form, can be formulated as
\begin{equation}\label{equation:stronglaplace}
\left\{
\begin{split}
	-\Delta u 	& = \bodyheat \qquad \text{in } \domain,\\
	u 				& = \dirichletdataheat \qquad \, \text{on } \dirichletboundary,\\
	\gradn u 	& = \neumanndataheat \qquad \text{on } \neumannboundary,
\end{split}
\right.
\end{equation}
where $u$ is the scalar temperature field, $\bodyheat$ is a heat source term, $\neumanndataheat$ represents the prescribed heat flux on the Neumann boundary, and $\dirichletdataheat$ is the prescribed temperature on the Dirichlet boundary. The normal gradient is defined as $\gradn u=  \nabla u \cdot \nn$.

The discretized solution to the strong formulation~\eqref{equation:stronglaplace} with the Dirichlet conditions enforced by Nitsche's method is denoted by $\uh \in \Uh = \mathcal{S}_{\alpha}^k (\mesh) \subset H^1(\mesh)$, with the corresponding test functions given by $\vh \in \Vh = \Uh$. We herein consider maximum regularity B-splines, \emph{i.e.}, $\alpha = \order -1$. The bilinear and linear forms in equation~\eqref{equation:galerkin} are
\begin{subequations}
\begin{align}
\B(\uh,\vh) &=  \int \limits_\domain \nabla \uh \cdot \nabla \vh \dO - \int \limits_{\dirichletboundary} \left(   (\gradn \uh)  \vh + \uh (\gradn v)   \right) \dG   \nonumber \\
&\phantom{=}  + \sum \limits_{\edge \in \dirichletmesh}  \int \limits_{\edge}  \tilde{\beta} \uh \vh \dG +  \sum \limits_{\face \in \ghost} \int \limits_{\face} \tilde{\gamma}_{g} \llbracket \gradn^\order \uh  \rrbracket \llbracket  \gradn^\order v \rrbracket \dG ,\label{eq:ah_laplace}\\
\F(\vh) &= \int \limits_{\domain} \bodyheat \vh \dO + \int \limits_{\neumannboundary} \neumanndataheat \vh \dG  - \int \limits_{\dirichletboundary} \dirichletdataheat  \gradn \vh \dG  \nonumber \\
&\phantom{=}+\sum \limits_{\edge \in \dirichletmesh} \int \limits_\edge \tilde{\beta} \dirichletdataheat  \vh  \dG , \label{eq:bh_laplace}
\end{align}%
\end{subequations}
where $\tilde{\beta}$ is the Nitsche stabilization parameter. This parameter should be selected and scaled (with the mesh size) appropriately, being large enough to ensure stability, while not being too large to cause a reduction in accuracy (see, \emph{e.g.}, Refs.~\cite{deprenter2018,badia2018}). The ghost-penalty operator in \eqref{eq:ah_laplace} controls the $k^{\rm th}$-order normal derivative jumps, indicated by $\llbracket \cdot \rrbracket$, over the interfaces of the elements which are intersected by the domain boundary $\boundary$. Since in this contribution B-splines of degree $\order$ with $C^{\order-1}$-continuity are considered, only the $\order^{\rm th}$ normal derivative is non-vanishing at the ghost mesh. As will be discussed in detail in Section~\ref{sec:stability}, upon approriate selection and scaling (with the mesh size) of $\tilde{\gamma}_{g}$, a Nitsche stabilization parameter, $\tilde{\beta}$, can be selected in such a way that stability of the formulation can be assured independent of the cut-cell configurations. To avoid loss of accuracy, the ghost-penalty parameter $\tilde{\gamma}_{g}$ should also not be too large \cite{badia2022}.
%
\subsubsection{Steady viscous flow}\label{section:stokes}
Steady viscous flow can be modeled by the Stokes equations,
\begin{equation}\label{equation:stokesequations}
\left\{
	\begin{split}
		-\nabla \cdot (2\mu \nabla^s \uu ) + \nabla p 	& = \bodystokes \qquad \text{in } \domain, \\
		\nabla \cdot \uu 								& = 0 \qquad \text{in } \domain, \\
		\uu  		& = \dirichletdatastokes \qquad  \text{on } \dirichletboundary, \\
		2 \mu \left( \nabla^s \uu  \right) \nn - p \nn  			& = \neumanndatastokes \qquad \text{on } \neumannboundary, 
	\end{split}
\right. 
\end{equation}
with velocity $\uu$, pressure $p$, constant viscosity $\mu$, body force $\bodystokes$,  Dirichlet data $\dirichletdatastokes$ and Neumann data $\neumanndatastokes$. By consideration of the solution in the abstract Galerkin problem \eqref{equation:galerkin} as a velocity-pressure pair, \emph{i.e.}, $\uh = (\uuh,\ph) \in \Uh = U^h_u \times U^h_p = [\mathcal{S}_{k-1}^k]^d \times \mathcal{S}_{k-1}^k \subset [H^1(\domain)]^d \times L^2(\domain)$ and the corresponding test functions as $\vh=(\vvh,\qh) \in \Vh = V^h_u \times V^h_p = \Uh$, the aggregate bilinear and linear forms corresponding to \eqref{equation:stokesequations} follow~as
\begin{subequations}
\label{eq:stokesforms}
\begin{align}
\B(\uh,\vh) &=  a^{\meshsize}_1(\uuh,\vvh)+a^{\meshsize}_2(\ph,\vvh)+a^{\meshsize}_2(\qh,\uuh)-a^{\meshsize}_3(\ph,\qh),
\label{equation:ah_stokes}
\\
\F(\vh) &= \int \limits_{\domain} \bodystokes \cdot  \vvh \dO + \int \limits_{\neumannboundary} \neumanndatastokes \cdot  \vvh \dG 
+ \int \limits_{\dirichletboundary} \left( \qh  \nn -2\mu  (\nabla^s \vvh)\nn \right)\cdot \dirichletdatastokes   \dG  
+\sum \limits_{\edge \in \dirichletmesh} \int \limits_\edge \tilde{\beta} \dirichletdatastokes \cdot \vvh \dG,
\end{align}
\end{subequations}
where
\begin{subequations}
\label{eq:ahcomponents}
\begin{align}
a^{\meshsize}_1(\uuh,\vvh)&=\int \limits_\domain  2 \mu \nabla^s \uuh : \nabla^s \vvh \dO  - \int \limits_{\dirichletboundary} 2 \mu \left(   (\nabla^s \uuh)\nn \cdot  \vvh + (\nabla^s \vvh)\nn \cdot \uuh \right) \dG   \nonumber \\
&\phantom{=}  + \sum \limits_{\edge \in \dirichletmesh}  \int \limits_{\edge} \tilde{\beta} \mu \uuh \cdot \vvh \dG 
+  \sum \limits_{F \in \ghost} \int \limits_F \tilde{\gamma}_{g}  \mu \llbracket \gradn^\order \uu^\meshsize  \rrbracket \cdot \llbracket  \gradn^\order \vvh \rrbracket \dG, 
\\
a^{\meshsize}_2(\ph,\vvh)&=-\int \limits_\domain \ph \nabla \cdot \vvh  \dO + \int \limits_{\dirichletboundary}\ph \vvh \cdot \nn \dG, 
\label{eq:ah2}
\\
a^{\meshsize}_3(\ph,\qh)&= \sum \limits_{F \in \skeleton} \int \limits_F \tilde{\gamma}_{s}\mu^{-1} \llbracket \gradn^\order \ph  \rrbracket  \llbracket  \gradn^\order \qh \rrbracket \dG. \label{equation:ah3_skeleton}
\end{align}	
\end{subequations}
For the selection of the Nitsche parameter, $\tilde{\beta}$, and ghost stabilization constant, $\tilde{\gamma}_{g}$, the same arguments apply as for the steady heat conduction problem discussed above. A discussion on the selection and scaling of these parameters for the Stokes problem is presented in Section~\ref{sec:stabilitystokes}.

An additional stability issue is encountered for the immersed Stokes flow problem \eqref{eq:stokesforms} on account of the selected equal-order optimal regularity spline spaces of degree $k$. In the conforming setting, inf-sup stability is achieved by adopting a suitable velocity-pressure pair, \emph{e.g.}, Taylor-Hood \cite{taylor1973, bazilevs2006, buffa2011, bressan2013} or Raviart-Thomas \cite{raviart1977, girault1979, buffa2011, evans2013raviart}. In the immersed setting, such pairs can still lead to pressure oscillations in the vicinity of cut elements \cite{hoang2017}. To resolve these pressure oscillations, the immersogeometric skeleton stabilization technique developed in Ref.~\cite{hoang2019} is applied. This stabilization technique can be regarded as the higher-order continuous version of the method proposed in Ref.~\cite{burman2006}, which has also been applied in the conforming isogeometric analysis setting \cite{hoang2017}. 

From equation~\eqref{equation:ah3_skeleton} it is seen that the skeleton stabilization term penalizes jumps in higher-order pressure gradients, where the parameter $\tilde{\gamma}_{s}$ should be selected such that oscillations are suppressed, while the influence of the additional term on the accuracy of the solution remains limited. The purpose of the skeleton stabilization method is to avoid pressure oscillations induced by inf-sup stability problems, allowing for the utilization of identical spaces for the velocity components and the pressure. Since the inf-sup stability problem is not restricted to the immersed boundary, the skeleton stabilization pertains to all interfaces of the background mesh. The appropriate selection and scaling of the skeleton stability parameter is discussed in detail in Section~\ref{sec:stabilitystokes}.


\subsection{Selection of the stabilization parameters: continuity and coercivity of the formulation}
\label{sec:stability}
Before considering a-posteriori error estimation in Section~\ref{sec:eeasection}, we first study the continuity and coercivity of the immersed formulations introduced above. We commence with the introduction of the following inequalities:
\begin{itemize}
  \item Using Young's inequality, it follows that for any constant $\varepsilon >0$ it holds that  
			\begin{align} 
				2 \lVert u^{h} \rVert_{L^{2}} \lVert {\tilde{u}}^{h} \rVert_{L^{2}} &\leq \varepsilon \lVert u^{h} \rVert_{L^{2}}^{2} + \cfrac{1}{\varepsilon} \lVert {\tilde{u}}^{h} \rVert_{L^{2}}^{2}  &  &\forall u^{h} \in U^{h},~\forall \tilde{u}^{h} \in U^{h}. \label{ineq:peterpaul}
			\end{align}
	In combination with the Cauchy-Schwarz inequality, this inequality can be applied to obtain
	\begin{align} \label{ineq:nitscheinequality}
			2 \int \limits_{\partial \Omega_{D}} (\partial_{n} u^{h}) u^{h} \, {\rm d}S &\leq \varepsilon  \lVert \partial_{n} u^{h} \rVert_{L^{2}(\partial \Omega_D)}^{2} + \cfrac{1}{\varepsilon} \lVert  u^{h} \rVert_{L^{2}(\partial \Omega_{D})}^{2} & &\forall u^h \in U^h. 
		\end{align}
	\item For any background element $K$ crossed by the boundary $\partial \Omega$, with $E = K \cap \partial \Omega$, under the assumption of shape regularity (\emph{i.e.}, provided with an upper bound on the length of the intersection of the boundary within one single element ${\rm meas}\,(K \cap \partial\Omega)$), it holds that (see, \emph{e.g.}, Ref.~\cite[Lemma 4.2]{evans2013})
			\begin{align} 
				\lVert \phi \rVert_{L^{2}(E)}^{2} 
				&\leq C_{T} \lVert {\meshsize_{\element}}^{-1/2} \phi \rVert_{L^{2}(\element)}^{2}  &  & \forall \phi \in P^{\order} \label{ineq:trace},
			\end{align}
where it is noted that this inequality holds for the finite-dimensional space $P^{\order}$ of tensor-product polynomials of order $k$ (not for functions in $H^1$ in general). The constant $C_{T}>0$,  referred to as the trace inequality constant, is independent of the size of the element, but dependent on the order $\order$. Note that the right part of the inequality contains the norm over the full background element $K$, and not just its intersection with the physical domain. 

Using inequality \eqref{ineq:trace}, the following bound for the normal gradient of $u^h$ on the immersed boundary is obtained:
\begin{equation}
\begin{aligned}
\lVert \partial_{n} u^{h} \rVert_{L^{2}(\partial \Omega)}^{2}
& \leq \lVert \nabla u^{h} \rVert_{L^{2}(\partial \Omega)}^{2}  
= \sum \limits_{ \edge \in \mathcal{T}_{\partial \Omega}} \lVert \nabla u^{h} \rVert_{L^{2}(\edge)}^{2}
\\
& \leq \sum \limits_{\element \in\mathcal{G}} C_{T} \lVert {\meshsize_{\element}}^{-1/2} \nabla u^{h} \rVert_{L^{2}(\element)}^{2}
\leq  {C}_{T} \lVert {\meshsize_{\mathcal{T}}}^{-1/2} \nabla u^{h} \rVert_{L^{2}(\mathcal{T})}^{2}\qquad\forall u^h \in U^h,
\end{aligned} \label{eq:trace_HvB}
\end{equation}
		with $h_{\mesh}$ defined in Eq.~\eqref{eq:meshsize} and where, with abuse of notation, the constant $C_T$ is used to both represent the local trace inequality constant (second line) and its global maximum (third line).
		\item Norms of functions over the entire background domain $\mesh$ can be bounded by norms over the physical domain $\Omega$ and the ghost penalty. Using the ghost-penalty, the gradients on the background mesh are bounded  by those in the physical domain. To demonstrate this bound, we split the norm over the background mesh as
\begin{equation}
\begin{aligned}
  \lVert \nabla u^{h} \rVert_{L^{2}(\mathcal{T})}^{2} &= \lVert \nabla u^{h} \rVert_{L^{2}(\mathcal{T} \setminus \mathcal{G})}^{2} + \lVert \nabla u^{h} \rVert_{L^{2}(\mathcal{G})}^{2} & &\\
  & \leq \lVert \nabla u^{h} \rVert_{L^{2}(\Omega)}^{2} + \lVert \nabla u^{h} \rVert_{L^{2}(\mathcal{G})}^{2} & &\\
  &\leq \lVert \nabla u^{h} \rVert_{L^{2}(\Omega)}^{2} + \sum \limits_{K \in \mathcal{G}} \lVert \nabla u^{h} \rVert_{L^{2}(K)}^{2} & & \forall u^h \in U^h.
\end{aligned}
  \label{eq:integralsplit}
\end{equation}
To show the last inequality, we consider an element $K \in \mathcal{G}$ which shares the interface $F$ with an element $K' \notin \mathcal{G}$ that completely lies inside $\Omega$, such that the volume integral over the background element $K'$ is included in the norm over $\Omega$. We will first demonstrate that the gradients on $K$ are controlled by the ghost penalty and the norms on the physical domain. Later on, elements in $\mathcal{G}$ that do not share an interface with an element in $\mesh\setminus\mathcal{G}$ will be considered by means of recursion. To demonstrate that the gradients on $K$ are bound by those in the physical domain, we define the polynomial extension of $\left. u^h \right|_{K'}$ as the global polynomial $\bar{u}^h_{K'} \in P^{\order}$ (see Figure~\ref{fig:integralsplit}). Using this extension, the spline function $u^h$ on the element $K$ can be decomposed as
\begin{align}
 \left. u^h \right|_K = \bar{u}^h_{K'} + \tilde{u}^h_{K'}.
 \label{eq:decomposition}
\end{align}
Let us consider $\boldsymbol{x}_F$ as a projection of $\boldsymbol{x}$ on the straight or flat interface $F$, such that $\boldsymbol{x}$ can be written as $\boldsymbol{x}_F + x_n\boldsymbol{n}_F$, where $x_n = (\boldsymbol{x} - \boldsymbol{x}_F ) \cdot \boldsymbol{n}_F$. Here, the interface coordinate $\boldsymbol{x}_F \in F$ is interpreted to be on the side of the element $K$, and related to the coordinate $\boldsymbol{x} \in K$. The function $\tilde{u}^h_{K'}$ has no support on $K'$ and has vanishing normal derivatives up to order $k$ at the interface $F$. By Taylor-series expansion one can infer
\begin{align}
 \tilde{u}^h_{K'}(\boldsymbol{x}) &= \frac{1}{k!} \partial_n^k (u^h(\boldsymbol{x}_F) - \bar{u}^h_{K'}(\boldsymbol{x}_F)) x_n^k = \frac{1}{k!} \partial_n^k \llbracket u^h(\boldsymbol{x}_F) \rrbracket x_n^k  & &\forall \boldsymbol{x} \in K, \forall u^h \in U^h.
 \label{eq:expansion}
\end{align}
This splitting is very natural through the use of maximum regularity splines (\emph{i.e.}, $\tilde{u}^h_{K'}$ contains all degrees of freedom of $K$ that are independent of $K'$). 

\begin{figure}[htp]
  \centering
  \begin{subfigure}[b]{0.45\textwidth}
  	 \centering
  	 \includegraphics[width=\textwidth]{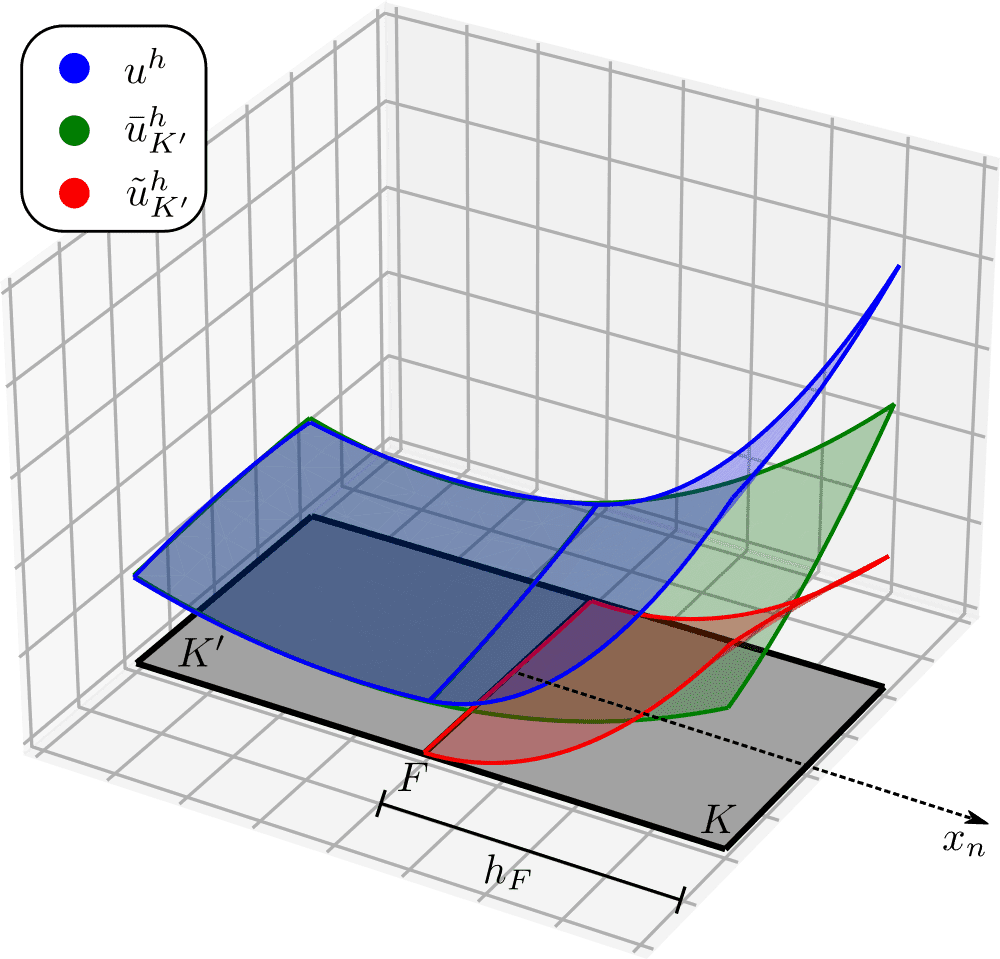}
  	 \caption{ }
  	 \label{fig:funcKK'}
  \end{subfigure}\hspace{0.08\textwidth}%
  \begin{subfigure}[b]{0.45\textwidth}
	 \centering
	 \includegraphics[width=\textwidth]{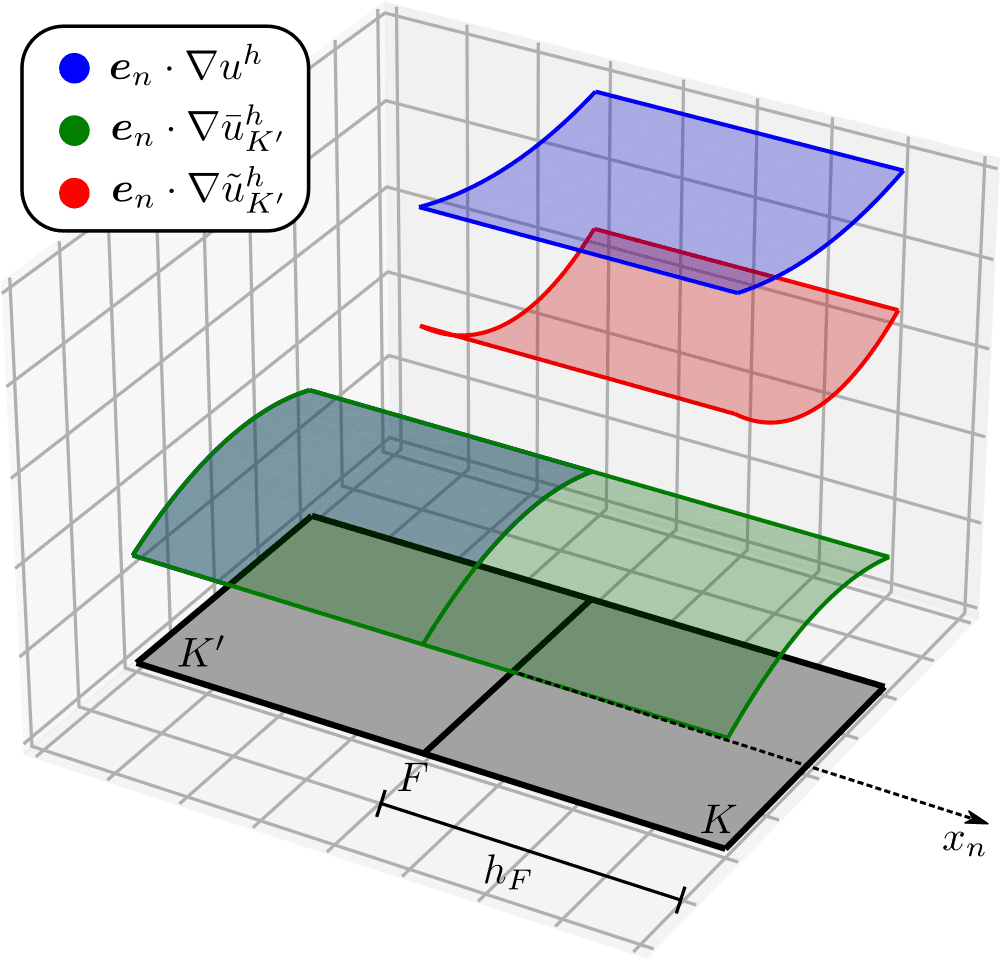}
	 \caption{ }
	 \label{fig:d2KK'}
  \end{subfigure}\\[24pt]
  \begin{subfigure}[b]{0.45\textwidth}
  	\centering
  	 \begin{tikzpicture}
  	 	\node (1) [black, fill=black!20, draw, rounded corners=5pt] {$\lVert \nabla  \bar{u}^h_{K'}  \rVert_{L^{2}(K)}^{2} \leq {\color{blue} C_Q} \lVert \nabla  \bar{u}^h_{K'}  \rVert_{L^{2}(K')}^{2}$ };
  	 	\node[below=-.3em of 1] (2) {\hspace{-3.5em} \includegraphics[width=\textwidth]{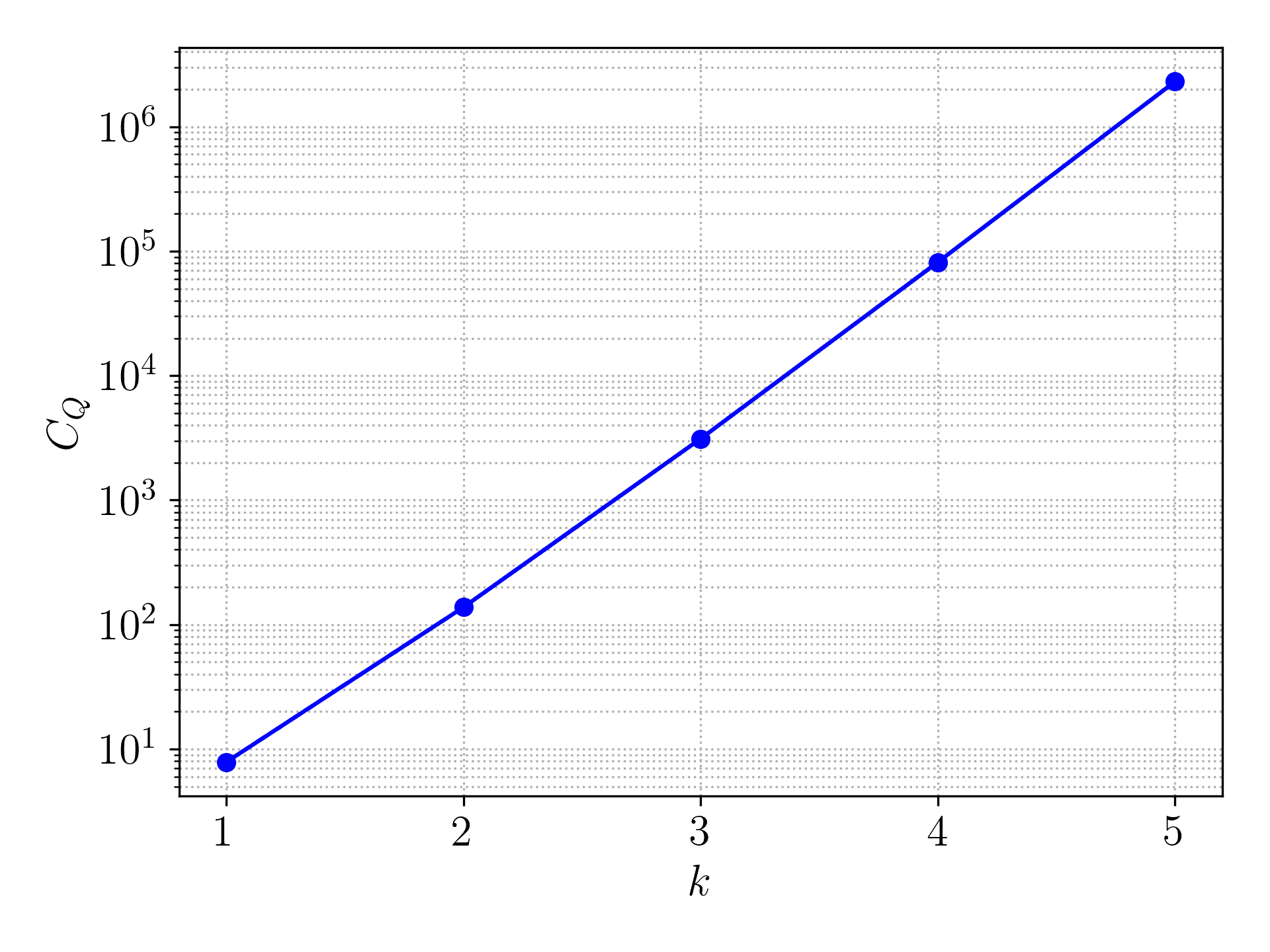}};
  	 \end{tikzpicture}
  	 \caption{ }
  	 \label{fig:C_Q}
  \end{subfigure}\hspace{0.08\textwidth}%
  \begin{subfigure}[b]{0.45\textwidth}
    \centering
      	 \begin{tikzpicture}
    	\node (1) [black, fill=black!20, draw, rounded corners=5pt] {$\lVert \nabla \tilde{u}^h_{K'} \rVert^2_{L^2(K)} \leq {\color{red} C_F} h_F^{2k-1} \lVert  \llbracket \partial_n^k  u^h \rrbracket \rVert^2_{L^2(F)} $ };
    	\node[below=-.3em of 1] (2) {\hspace{-2em} \includegraphics[width=\textwidth]{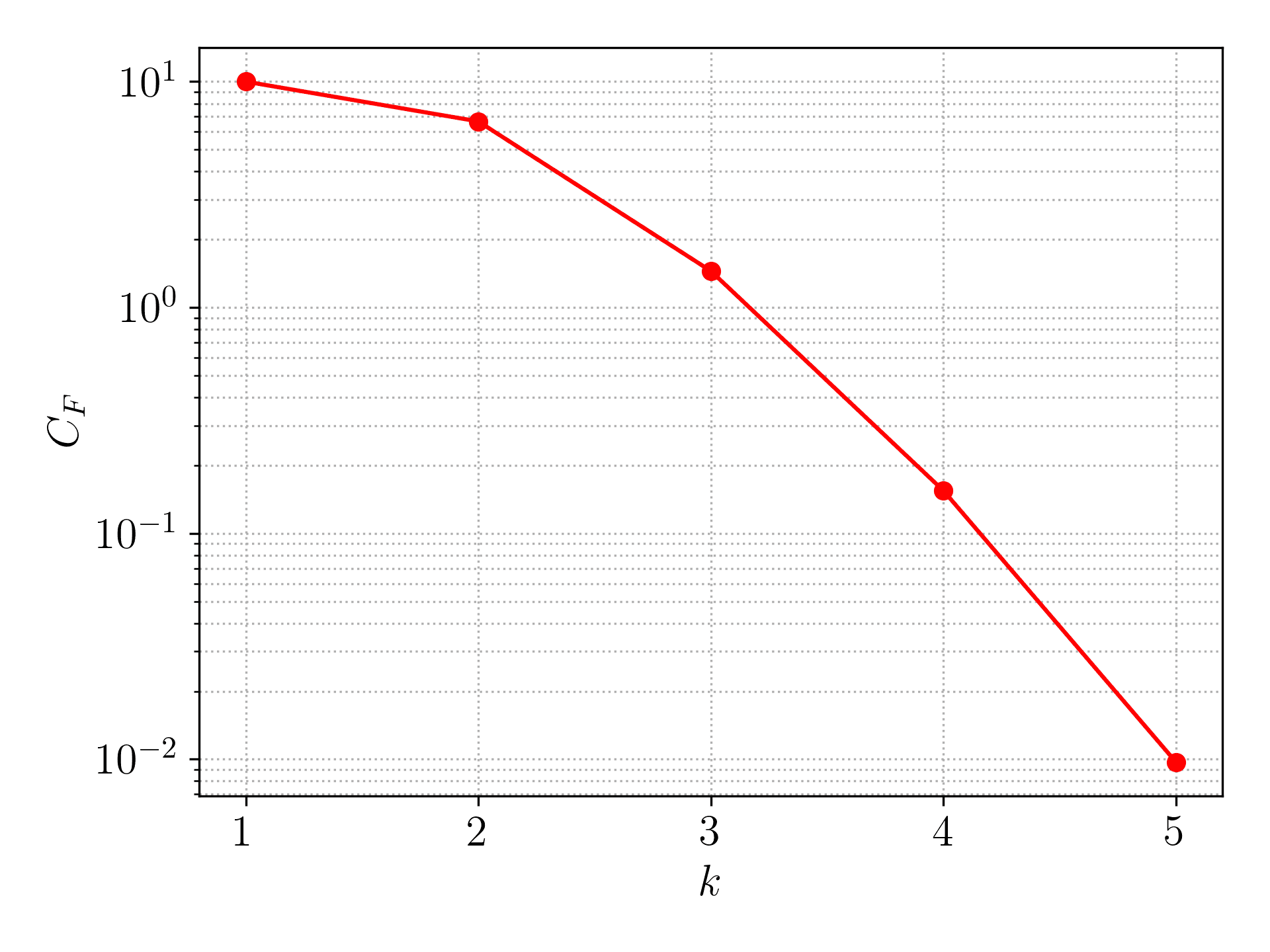}};
    \end{tikzpicture}
  	 \caption{ }
  	 \label{fig:C_F}
  \end{subfigure}%
  \caption{Panel (\subref{fig:funcKK'}) is an illustration of a second order B-spline on an element $K \in  \mathcal{G}$ and its adjacent element $K'$ with an interface $F$. Panel (\subref{fig:d2KK'}) is its second order gradient in the direction normal to the interface (with $\boldsymbol{e}_n$ the unit vector in the normal direction). Panels (\subref{fig:C_Q}) and (\subref{fig:C_F}) show the dependence of the constants in \eqref{eq:polynomialinequality} and \eqref{eq:expansionintegrated} on the order $\order$.}
  \label{fig:integralsplit}
\end{figure}

For the polynomial extension $\bar{u}^h_{K'}$ it holds that 
\begin{align}
\lVert \nabla  \bar{u}^h_{K'}  \rVert_{L^{2}(K)}^{2} &\leq C_Q \lVert \nabla  \bar{u}^h_{K'}  \rVert_{L^{2}(K')}^{2}  & \forall u^h \in U^h,
\label{eq:polynomialinequality}
\end{align}
where the constant $C_Q$ is independent of the mesh size, but dependent on the order of the approximation and the ratio of the size of the elements at either side of the interface. The order-dependence of this constant is illustrated in Figure~\ref{fig:C_Q}. The presented results have been computed by solving the generalized eigenvalue problem corresponding to Eq.~\eqref{eq:polynomialinequality}.

From the definition of the expansion $\tilde{u}^h_{K'}$ in equation~\eqref{eq:expansion} it follows that
\begin{equation}
\begin{aligned}
\lVert \nabla \tilde{u}^h_{K'} \rVert^2_{L^2(K)}
&=  \lVert  \frac{1}{k!}\left( \nabla_F \partial_n^k \tilde{u}^h \right) x_n^k +   \frac{k}{k!}\left(  \partial_n^k \tilde{u}^h \right) x_n^{k-1}    \rVert^2_{L^2(K)} \\
&\leq  \frac{2}{(k!)^2} \left(\lVert  \left( \nabla_F \partial_n^k \tilde{u}^h \right) x_n^k  \rVert^2_{L^2(K)} + \lVert  k \left(  \partial_n^k \tilde{u}^h \right) x_n^{k-1}    \rVert^2_{L^2(K)}\right)\\
&\leq \frac{2}{(k!)^2}  \left( \lVert\nabla_F \partial_n^k \tilde{u}^h \rVert^2_{L^2(F)} (2k+1)^{-1} h_F^{2k+1} \right. \\
&\phantom{\leq \frac{2}{(k!)^2}  \lVert\nabla_F} \left. + ~ k^2 \lVert \partial_n^k  \tilde{u}^h \rVert^2_{L^2(F)} (2k-1)^{-1} h_F^{2k-1}
\right)\\
&\leq \frac{C_F}{2} h_F^{2k-1} \lVert \partial_n^k  \tilde{u}^h \rVert^2_{L^2(F)} = \frac{C_F}{2} h_F^{2k-1} \lVert  \llbracket \partial_n^k  u^h \rrbracket \rVert^2_{L^2(F)},  \label{eq:expansionintegrated}
\end{aligned}
\end{equation}
with $h_F$ the size of $K$ in the direction normal to the interface and where $\nabla_F$ denotes the surface gradient in the interface $F$ and where use has been made of the polynomial inequality $\lVert \nabla_F f^h \rVert^2_{L^2(F)} \lesssim h_F^{-2} \lVert f^h \rVert^2_{L^2(F)}$ for all $f^{h} \in P^{\order}$ \cite{kroo2009}. The dependence of the constant $C_F$ in the inequality \eqref{eq:expansionintegrated} on the order is illustrated in Figure~\ref{fig:C_F}. This constant is independent of the mesh size.

Substituting the decomposition \eqref{eq:decomposition} in equation \eqref{eq:integralsplit} yields
\begin{align}
 \lVert \nabla u^{h} \rVert_{L^{2}(\mathcal{G})}^{2} &= \sum \limits_{K \in \mathcal{G}} \lVert \nabla \bar{u}^h_{K'} + \nabla  \tilde{u}^h_{K'} \rVert_{L^{2}(K)}^{2} \leq 2 \sum \limits_{K \in \mathcal{G}} \left( \lVert \nabla  \bar{u}^h_{K'}  \rVert_{L^{2}(K)}^{2} + \lVert  \nabla \tilde{u}^h_{K'} \rVert_{L^{2}(K)}^{2}  \right) & & \forall u^h \in U^h.
\end{align}
Using the inequalities \eqref{eq:polynomialinequality} and \eqref{eq:expansionintegrated}, and noting that since $K' \in \mathcal{T} \setminus \mathcal{G}$ it follows that $\lVert \nabla  \bar{u}^h_{K'}  \rVert_{L^{2}(K')}^{2} =\lVert \nabla  u^h  \rVert_{L^{2}(K')}^{2}$, then results in
\begin{equation}
\begin{aligned}
 \lVert \nabla u^{h} \rVert_{L^{2}(\mathcal{G})}^{2} &\leq (C_G - 1) \lVert \nabla  u^h  \rVert_{L^{2}(\Omega)}^{2} + \sum \limits_{F \in \mathcal{F}} C_F h_F^{2k-1} \lVert  \llbracket \partial_n^k  u^h \rrbracket \rVert^2_{L^2(F)}   & & \forall u^h \in U^h,
\end{aligned}
\label{eq:gradientboundedness}
\end{equation}
where $C_G = 1 + 2 \max_{K \in \mathcal{G}}(C_{Q})$. To obtain this result, the inequality is first applied to the layer of elements in $\mathcal{G}$ that share an interface with the interior mesh $\mathcal{T} \setminus \mathcal{G}$. With control over the gradients in this first layer, the inequality is then applied to a second layer of elements. This recursive application is repeated until all elements in $\mathcal{G}$ have been considered.  As a result of this recursive application of the ghost inequality, the constant $C_G$ depends on the number of layers, which in turn depends on the mesh size.

The boundedness of the gradients on the background mesh finally follows by substitution of \eqref{eq:gradientboundedness} in \eqref{eq:integralsplit}:
\begin{align}
  \lVert \nabla u^{h} \rVert_{L^{2}(\mathcal{T})}^{2} &\leq C_G \lVert \nabla u^{h} \rVert_{L^{2}(\Omega)}^{2} + \sum \limits_{F \in \mathcal{F}} C_F h_F^{2k-1} \lVert  \llbracket \partial_n^k  u^h \rrbracket \rVert^2_{L^2(F)}   & & \forall u^h \in U^h.
  \label{ineq:ghostinequality}
\end{align}

\item For any $p^h \in U^h_p$, there exists a $\mathbf{w}^h \in  V^h_u$ such that
\begin{subequations}
	\begin{align}
		-a_2^h(\ph, \mathbf{w}^h) &\geq \Big. C_{1} \lVert \mu^{-1/2} \ph \rVert_{L^{2}(\mathcal{T})}^{2} - C_{2} \sum \limits_{F \in \skeleton} \int \limits_F \meshsize_\face^{2\order+1}\mu^{-1} \llbracket \gradn^\order \ph  \rrbracket  \llbracket  \gradn^\order \ph \rrbracket \dG, \Big.
		\label{eq:skeletoninequality}\\
		\ltrivert \mathbf{w}^h \rtrivert_u &\leq C_3 \lVert \mu^{-1/2} \ph \rVert_{L^{2}(\mathcal{T})},
	\end{align}
	\label{eq:skeletoninequalities}%
\end{subequations}
with $a^{\meshsize}_2$ according to~(\ref{eq:ah2}) and $\ltrivert \,\cdot\, \rtrivert_u$ the velocity energy norm (see Section~\ref{sec:stabilitystokes}), for certain positive constants $C_1,C_3>0$ and a non-negative constant $C_2 \geq 0$. Existence of a velocity field~$\mathbf{w}^h$ in accordance with~\eqref{eq:skeletoninequalities} is established in~\cite[Lemma 3.11]{hansbo2014} for piecewise linear ($k=1$) polynomials. However, this result generalizes to polynomial orders $k\geq{}1$ and increased continuity of the pressure and velocity spaces. Proof of~\eqref{eq:skeletoninequalities} in the general case is however extensive, and beyond the scope of the present manuscript.
\end{itemize}

\subsubsection{Steady heat conduction} \label{sec:stabilitylaplace}


Continuity of the bilinear form \eqref{eq:ah_laplace} cannot be shown in the $H^1$-norm on account of the immersed boundary terms, and coercivity cannot be shown on the infinite-dimensional space. However, with an appropriate selection of the stabilization parameters, continuity and coercivity can be established with respect to the mesh-dependent norm
\begin{align}
	\ltrivert u^{h} \rtrivert^{2} :&= \lVert \nabla u^{h} \rVert_{L^{2}(\mathcal{T})}^{2} +  \lVert \tilde{\beta}^{-1/2} \partial_{n} u^{h} \rVert_{L^{2}(\partial \Omega_{D})}^{2} + \lVert \tilde{\beta}^{1/2} u^{h} \rVert_{L^{2}(\partial \Omega_{D})}^{2} && \nonumber \\
&\phantom{=} + \sum \limits_{F \in \mathcal{F}_{\rm ghost}} \lVert \tilde{\gamma}_{g}^{1/2} \llbracket \partial_{n}^{k} u^h \rrbracket \rVert_{L^{2}(F)}^{2} && \forall u^{h} \in \Uh, \label{eq:energynorm}
\end{align}
which we refer to as the \emph{energy norm}.

The bilinear form \eqref{eq:ah_laplace} is continuous on $U^h \times V^h$ if there exists a constant, $C>0$, independent of the mesh size, such that
	\begin{align}
		a^{h}(u^{h}, v^{h}) &\leq C \ltrivert u^{h} \rtrivert \ltrivert v^{h} \rtrivert & \forall u^h \in U^h,~\forall v^h \in V^h.\label{eq:continuity}
	\end{align}
Using the Cauchy-Schwarz inequality, for all $u^h \in U^h,~v^h \in V^h$ one obtains
\begin{align*}
		a^{h}(u^h, v^h)
		&\leq \Big. \lVert \nabla u^{h} \rVert_{L^{2}(\Omega)} \lVert \nabla v^{h} \rVert_{L^{2}(\Omega)} + \lVert \tilde{\beta}^{-1/2} \partial_{n} u^{h} \rVert_{L^{2}(\partial \Omega_{D})}  \lVert \tilde{\beta}^{1/2} v^{h} \rVert_{L^{2}(\partial \Omega_{D})} \Big. \\
		&\phantom{=} \Big. + \lVert \tilde{\beta}^{1/2} u^{h} \rVert_{L^{2}(\partial \Omega_{D})} \lVert \tilde{\beta}^{-1/2} \partial_{n} v^{h} \rVert_{L^{2}(\partial \Omega_{D})} \Big. \\
		&\phantom{=} \Big. + \sum_{E \in \mathcal{T}_{\partial \Omega_{D}}} \lVert \tilde{\beta}^{1/2} u^{h} \rVert_{L^{2}(E)} \lVert \tilde{\beta}^{1/2} v^{h} \rVert_{L^{2}(E)} \Big. \\
		&\phantom{=} \Big. + \sum_{F \in \mathcal{F}_{\rm ghost}} \lVert \tilde{\gamma}_{g}^{1/2} \llbracket \partial_{n}^{\order} u^{h}  \rrbracket \rVert_{L^{2}(F)}  \lVert \tilde{\gamma}_{g}^{1/2} \llbracket \partial_{n}^{\order} v^{h}  \rrbracket \rVert_{L^{2}(F)} \Big. .
	\end{align*}
Since each of the norms in this expression is bounded from above by the energy norm \eqref{eq:energynorm}, it follows that $a^{h}(u^h, v^h) \leq 5 \ltrivert u^{h} \rtrivert \ltrivert v^{h} \rtrivert$. Hence the bilinear form is continuous.

The bilinear form \eqref{eq:ah_laplace} is uniformly (\emph{i.e.}, independent of $h$) coercive on $U^h$ if there exists a constant, $c>0$, such that
	\begin{align} \label{eq:coercivitylaplace}
		&a^{h}(u^{h}, u^{h}) \geq c \ltrivert u^{h} \rtrivert^{2} &  & u^{h} \in U^{h}.
	\end{align}
To demonstrate that this is indeed the case, we apply the inequalities \eqref{ineq:nitscheinequality} and \eqref{ineq:ghostinequality} to obtain
\begin{align*}
	a^{h}(u^h, u^h)  &\geq \Bigg( \frac{1}{C_G} \lVert \nabla u^{h} \rVert_{L^2(\mathcal{T})}^{2} - \sum \limits_{F \in \mathcal{F}_{\rm ghost}} \frac{C_F}{C_G} \int \limits_{F} h_F^{2k-1}\llbracket \partial_{n}^{\order} u^{h}  \rrbracket^{2} \, {\rm d}S \Bigg) \nonumber \\
	&\phantom{\geq} \Bigg. - \Bigg( \varepsilon  \lVert \partial_{n} u^{h} \rVert_{L^{2}(\partial \Omega_D)}^{2} + \varepsilon^{-1} \lVert u^h \rVert_{L^{2}(\partial \Omega_{D})}^{2} \Bigg) \Bigg. \nonumber \\
	&\phantom{\geq} \Bigg. + \sum \limits_{ E \in \mathcal{T}_{\partial \Omega_{D}}}  \lVert {\tilde{\beta}}^{\frac{1}{2}} u^{h} \rVert_{L^{2}(E)}^{2} + 	\sum \limits_{F \in \mathcal{F}_{\rm ghost}} \int \limits_{F} \tilde{\gamma}_{g}\llbracket \partial_{n}^{\order} u^{h}  \rrbracket^{2} \, {\rm d}S. \Bigg.
\end{align*}
Application of the trace inequality \eqref{eq:trace_HvB} and collecting terms then yields
\begin{align*}
	a^{h}(u^h, u^h)  &\geq  \bigg\lVert \bigg( \frac{1}{C_G} -\varphi C_T h_{\mathcal{T}}^{-1}  \bigg)^{\frac{1}{2}} \nabla u^{h} \bigg\rVert_{L^2(\mathcal{T})}^{2} +   \Big\lVert \left(\varphi-\varepsilon \right)^{\frac{1}{2}} \partial_{n} u^{h} \Big\rVert_{L^{2}(\partial \Omega_D)}^{2}  \nonumber \\
	&\phantom{\geq} + \Big\lVert  \left( \tilde{\beta} - \varepsilon^{-1}\right)^{\frac{1}{2}} u^{h} \Big\rVert_{L^{2}(\partial \Omega_{D})}^{2} + \sum \limits_{F \in \mathcal{F}_{\rm ghost}} \int \limits_{F} \bigg( \tilde{\gamma}_{g} - \frac{C_F}{C_G} h_F^{2k-1} \bigg) \llbracket \partial_{n}^{\order} u^{h}  \rrbracket^{2} \, {\rm d}S,
\end{align*}
for arbitrary $\varphi > 0$. By selecting element-wise constants $0 < \varepsilon \leq \varphi$ and $0 < \varphi \leq \frac{C_G^{-1}}{C_T } h_{\mathcal{T}}$, one can infer that coercivity is ensured provided that
\begin{align}
\tilde{\beta} &= \beta h_K^{-1} \geq C_T C_G h_{K}^{-1}
 &
 \tilde{\gamma}_{g} &= \gamma_{g} h_F^{2k-1} \geq \frac{C_F }{C_G} h_F^{2k-1}.
\end{align}
for all elements $K$ and interfaces $F$, where the positive constants $\beta \geq C_T C_G$ and $\gamma_{g} \geq C_F C_G^{-1}$ are independent of the mesh size. The interface length scale is defined as $h_F=\max{(h_K,h_{K'})}$ with $K$ and $K'$ being the elements on either side of the interface $F$. The rational behind this choice is that the ghost stabilization term scales with $h_F^{2k-1}$ ($k\geq1$) and that hence the larger element size ensures that the stability constant is sufficiently large.
\subsubsection{Steady viscous flow}\label{sec:stabilitystokes}
Recalling that for the Stokes problem $u^h=(\boldsymbol{u}^h,p^h)$, we define the mesh-dependent energy norm as
\begin{equation} \label{eq:energynormstokes}
	\ltrivert  u^h \rtrivert^2 = \ltrivert (\boldsymbol{u}^h,p^h) \rtrivert^2 =  \ltrivert \uuh \rtrivert^{2}_u + \ltrivert \ph \rtrivert^{2}_p,
\end{equation}
with 
\begin{subequations}	
	\begin{align}
		\ltrivert \uuh \rtrivert^{2}_u &:= \lVert \mu^{1/2} \nabla^s \boldsymbol{u}^{h} \rVert_{L^{2}(\mathcal{T})}^{2} +  \lVert \tilde{\beta}^{-1/2}  \mu^{1/2} \partial_{n} \boldsymbol{u}^{h} \rVert_{L^{2}(\partial \Omega_{D})}^{2} + \lVert \tilde{\beta}^{1/2}  \mu^{1/2}\boldsymbol{u}^{h} \rVert_{L^{2}(\partial \Omega_{D})}^{2} \nonumber \\ 
		&\phantom{:=}+ \sum \limits_{F \in \mathcal{F}_{\rm ghost}} \lVert \tilde{\gamma}_{g}^{1/2}  \mu^{1/2} \llbracket \partial_{n}^{k} \boldsymbol{u}^h \rrbracket \rVert_{L^{2}(F)}^{2}, \label{eq:velocitynorm} \\
		\ltrivert \ph \rtrivert^{2}_p &:=  \lVert \mu^{-1/2} p^{h} \rVert_{L^{2}(\mathcal{T})}^{2} + \sum \limits_{F \in \mathcal{F}_{\rm skeleton}} \lVert \tilde{\gamma}_{s}^{1/2}  \mu^{-1/2} \llbracket \partial_{n}^{k} p^h \rrbracket \rVert_{L^{2}(F)}^{2}. \label{eq:pressurenorm}
	\end{align}
\end{subequations}
Continuity of the bilinear form \eqref{equation:ah_stokes} with respect to this energy norm in the sense of \eqref{eq:continuity} follows directly by application of the Cauchy-Schwarz inequality to all terms in \eqref{equation:ah_stokes}.

With an appropriate selection of the stability parameters for the Stokes problem, it holds that the bilinear form \eqref{equation:ah_stokes} is inf-sup stable in accordance with
\begin{align}
	&\sup_{v^h \in V^h \setminus \{0\}} \frac{a^h(u^h,v^h)}{\ltrivert v^h \rtrivert} \geq c^\star \ltrivert u^h \rtrivert & &\forall u^h \in U^h,  
	\label{eq:infsupstabilitystokes}
\end{align}	  
where $c^\star>0$ is referred to as the inf-sup stability constant. To demonstrate this stability property, we recall the splitting of the bilinear form $\B$ according to~\eqref{equation:ah_stokes} and~\eqref{eq:ahcomponents}. We now take a function $\varphi^{h} = (\vvh,\qh) = (\uuh - \alpha \mathbf{w}^h, -\ph) \in \Vh$ where $\mathbf{w}^h$ depends on $\ph$ in accordance with \eqref{eq:skeletoninequalities}, and with some constant $\alpha > 0$, such that
\begin{align}
	\sup_{v^h \in V^h \setminus \{0\}} \frac{a^h(u^h,v^h)}{\ltrivert v^h \rtrivert}
	&\geq  \frac{a^h(u^h,\varphi^h)}{\ltrivert \varphi^h \rtrivert}  & &  \nonumber \\
	& \geq  \frac{a_{1}^{h}(\uuh, \uuh) - \alpha a_{1}^{h}(\uuh, \mathbf{w}^h) - \alpha a^{h}_2(\ph, \mathbf{w}^h) + a_3^h(\ph, \ph)}{\ltrivert \varphi^h \rtrivert}  & & \forall u^h \in U^h.
\end{align}	 
Following Section~\ref{sec:stabilitylaplace}, $a_1^{h} (\uuh, \uuh)$ is coercive (with constant $c^{\star}_{u}$) and $a_1^{h}(\uuh, \mathbf{w}^{h})$ is continuous (with constant $C^{\star}_{u}$) with respect to the velocity energy norm \eqref{eq:velocitynorm} in accordance with Eqs.~\eqref{eq:continuity} and \eqref{eq:coercivitylaplace}, respectively. Hence, 
\begin{align}
	\sup_{v^h \in V^h \setminus \{0\}} \frac{a^h(u^h,v^h)}{\ltrivert v^h \rtrivert} \geq  \frac{c^{\star}_{u} \ltrivert \uuh \rtrivert^{2}_u - \alpha C^{\star}_{u} \ltrivert \uuh \rtrivert_u \ltrivert \mathbf{w}^{h} \rtrivert_u - \alpha a_2^{h}(\ph, \mathbf{w}^h) + a_3^h(\ph, \ph)}{\sqrt{2 \ltrivert \uuh \rtrivert_u^2 + 2 \alpha^2 \ltrivert  \mathbf{w}^h \rtrivert_u^2 +  \ltrivert p^h \rtrivert_p^2}}  & &\forall u^h \in U^h,
\end{align}	 
where use has been made of
\begin{align}
\ltrivert \varphi^h \rtrivert = \sqrt{\ltrivert \uuh -\alpha \mathbf{w}^h \rtrivert_u^2 + \ltrivert p^h \rtrivert_p^2 } \leq \sqrt{2 \ltrivert \uuh \rtrivert_u^2  + 2 \alpha^2 \ltrivert  \mathbf{w}^h \rtrivert_u^2 +  \ltrivert p^h \rtrivert_p^2}.
\end{align}
From the inequalities \eqref{eq:skeletoninequalities} it follows that
\begin{align}
	& \sup_{v^h \in V^h \setminus \{0\}} \frac{a^h(u^h,v^h)}{\ltrivert v^h \rtrivert} \nonumber \\
	&\phantom{\sup_{v^h \in V^h \setminus \{0\}}} \geq  \frac{c^{\star}_{u} \ltrivert \uuh \rtrivert^{2}_u - \alpha C^{\star}_{u} C_3 \ltrivert \uuh \rtrivert_u \lVert \mu^{-1/2} \ph \rVert_{L^{2}(\mathcal{T})} +  \alpha C_{1} \lVert \mu^{-1/2} \ph \rVert_{L^{2}(\mathcal{T})}^{2} }{\sqrt{2 \ltrivert \uuh \rtrivert_u^2 + (1+2\alpha^2 C_3^2) \ltrivert p^h \rtrivert_p^2}}  & & \nonumber\\
	&\phantom{\sup_{v^h \in V^h \setminus \{0\}}} +  \frac{ \sum \limits_{F \in \skeleton} \left( 1 - \alpha C_{2} \tilde{\gamma}_{s}^{-1} \meshsize_{F}^{2\order+1} \right) \lVert \tilde{\gamma}_{s}^{1/2} \mu^{-1/2} \llbracket \gradn^\order \ph  \rrbracket \rVert_{F}^{2} }{\sqrt{2 \ltrivert \uuh \rtrivert_u^2 + (1+2\alpha^2 C_3^2) \ltrivert p^h \rtrivert_p^2}}  & &\forall u^h \in U^h,
\end{align}	 
which, using Young's inequality \eqref{ineq:peterpaul} with $\varepsilon=1$, can be reformulated as
\begin{align}
	& \sup_{v^h \in V^h \setminus \{0\}} \frac{a^h(u^h,v^h)}{\ltrivert v^h \rtrivert} \nonumber \\
	&\phantom{\sup_{v^h \in V^h \setminus \{0\}}} \geq  \frac{
	\left( c^{\star}_{u} -  \cfrac{\alpha C^{\star}_{u} C_3}{2} \right) \ltrivert \uuh \rtrivert^{2}_u +  \alpha \left( C_{1} - \cfrac{C^{\star}_{u} C_3}{2} \right) \lVert \mu^{-1/2} \ph \rVert_{L^{2}(\mathcal{T})}^{2}
	 }{C_{4} \ltrivert \uh \rtrivert}  & & \nonumber\\
	&\phantom{\sup_{v^h \in V^h \setminus \{0\}}} +  \frac{ \sum \limits_{F \in \skeleton} \left( 1 - \alpha C_{2} \tilde{\gamma}_{s}^{-1} \meshsize_{F}^{2\order+1} \right) \lVert \tilde{\gamma}_{s}^{1/2} \mu^{-1/2} \llbracket \gradn^\order \ph  \rrbracket \rVert_{F}^{2}  }{C_{4} \ltrivert \uh \rtrivert}  & &\forall u^h \in U^h.
\end{align}	 
Inf-sup stability as in \eqref{eq:infsupstabilitystokes} then holds, provided that 
\begin{align}
 \alpha &< \frac{2 c^{\star}_{u}}{C^{\star}_{u} C_3},   &  \tilde{\gamma}_s = \gamma_s h_F^{2k+1} &\geq \alpha C_2 h_F^{2k+1} .
\end{align}

We note that the skeleton penalty has two purposes: \emph{i)} It extends the stability of the pressure field to the background grid $\mathcal{T}$ as in \eqref{eq:skeletoninequalities}, in the same way as for the ghost penalty discussed in Section~\ref{sec:stabilitylaplace}. Since stability is here defined with respect to the $L^{2}$-norm of the pressure field, the skeleton stability constant $\tilde{\gamma}_{s}$ scales with $\meshsize^{2\order+1}$, following the same reasoning as in Eq.~\eqref{ineq:ghostinequality}; \emph{ii)} It ensures the inf-sup stability for equal-order discretizations, essentially meaning that pressure oscillations in the interior are penalized. This is the reason why this term is applied over the complete skeleton and not only the ghost interfaces.

%% file: sections/eea.tex
\section{Error estimation and adaptivity}\label{sec:eeasection}
We study a posteriori error estimation and adaptivity for immersogeometric analysis. In Section~\ref{section:errorestimation} we first introduce a residual-based error indicator, and elaborate it for the heat conduction problem and viscous flow problem introduced in the previous section. In Section~\ref{section:eeaprocedure} the refinement strategy is discussed. 

\subsection{Residual-based error estimation}\label{section:errorestimation}
We propose an error estimator pertaining to the background mesh, $\mesh$, of the form
\begin{align}
 \mathcal{E} = \sqrt{\sum \limits_{\element \in \mesh} \eta_\element^2},
 \label{equation:estimator}
\end{align}
where the element-wise error indicators, $\eta_\element$, will serve to guide an adaptive refinement procedure. The derivations of the indicators for the heat conduction problem and viscous flow problem will be elaborated in the following sections.

From an abstract perspective, the element-wise error indicators are defined in such a way that the estimator \eqref{equation:estimator} bounds the residual from above as
\begin{equation}
  \mathcal{E} \gtrsim \| r^h \|_{\widehat{V}^{h\ast}}.
  \label{equation:errorresidual}
\end{equation}
In this expression, the residual and its (dual) norm are defined as
\begin{subequations}
\begin{align}
r^h(\widehat{v}^{h}) &:= r^h(u^h)(\widehat{v}^h) :=b^h(\widehat{v}^h) -a^h(u^h,\widehat{v}^h),  \\
\lVert r^{h} \rVert_{\widehat{V}^{h\ast}} &:=\sup_{\widehat{v}^h \in \widehat{V}^h \setminus \{0\}} \frac{r^h(\widehat{v}^h)}{\ltrivert \widehat{v}^h \rtrivert},
\label{equation:residualdefinition}
\end{align}
\end{subequations}
The function space $\widehat{V}^h\supset{}V^h$ corresponds to a suitable extension of~$V^h$ in such a manner that 
$\widehat{V}^h$ contains an approximation of the solution (possibly the solution itself) that is sufficiently accurate to estimate the error in the approximation~$u^h\in{}V^h$. An example of such an extended space is an order elevated approximation space on the same mesh and with the same regularity as the space $V^h$, or an approximation space with the same order and the same regularity on a hierarchically refined mesh. The Galerkin approximation problem in $\widehat{V}^h$ writes
\begin{equation}
\label{eq:extendedGalerkin}
\widehat{u}^h\in\widehat{V}^h:\qquad \widehat{a}^h(\widehat{u}^h,\widehat{v}^h)=b^h(\widehat{v}^h)\qquad\forall{}\widehat{v}^h\in\widehat{V}^h.
\end{equation}
The bilinear form~$\widehat{a}^h:\widehat{V}^h\times{}\widehat{V}^h\to\mathbb{R}$ is an extension of the original bilinear form~$a^h$ 
according to
\begin{equation}
\widehat{a}^h(\widehat{u}^h,\widehat{v}^h)=a^h(\widehat{u}^h,\widehat{v}^h)+s^h(\widehat{u}^h,\widehat{v}^h),
\end{equation}
where the auxiliary symmetric bilinear form~$s^h$ contains additional stabilization terms. For instance, for an order-elevated space the bilinear form $s^h$ contains jumps of higher-order normal derivatives (cf.~\eqref{eq:ah_laplace} and see~\cite{larson2020,hoang2018}), and
for a hierarchically refined mesh $s^h$ contains the stabilization terms on the supplementary faces of the mesh. The additional stabilization terms vanish on the original approximation space~$V^h$, i.e.
\begin{equation}
\label{eq:shvanishes}
s^h:V^h\times\widehat{V}^h\to\{0\}.
\end{equation}
We equip $\widehat{V}^h$ with the extended energy norm according to
\begin{equation}
\ltrivert\,\cdot\,\rtrivert_{\widehat{V}^h}^2=  \ltrivert\,\cdot\,\rtrivert^2+s^h(\cdot,\cdot),
\end{equation}
%
With a suitable choice of the stabilization parameters in~$a^h$ and $s^h$, the bilinear
form~$\widehat{a}^h$ is weakly coercive and continuous. It is to be noted that this may require that the stabilization parameters in~$a^h$ are larger than would be required for weak coercivity of~$a^h$ on $V^h\times{}V^h$. By virtue of~\eqref{eq:extendedGalerkin}\nobreakdash--\eqref{eq:shvanishes} and the weak coercivity and linearity of~$\widehat{a}^h$, the following chain of inequalities holds:
\begin{equation}
\begin{aligned}
 \mathcal{E} \gtrsim \ltrivert r^h(\widehat{v}^h)\rtrivert_{\widehat{V}^{h\ast}}
 &= \sup_{\widehat{v}^h \in \widehat{V}^h \setminus \{0\}} \frac{b^h(\widehat{v}^h)-a^h(u^h,\widehat{v}^h)}{\ltrivert \widehat{v}^h \rtrivert} \\
 &= \sup_{\widehat{v}^h \in \widehat{V}^h \setminus \{0\}} \frac{\widehat{a}^h(\widehat{u}^h,\widehat{v}^h)-\widehat{a}^h(u^h,\widehat{v}^h)}{\ltrivert \widehat{v}^h \rtrivert} \\ 
 &=\sup_{\widehat{v}^h \in \widehat{V}^h \setminus \{0\}} \frac{\widehat{a}^h(\widehat{u}^h-u^h,\widehat{v}^h)}{\ltrivert \widehat{v}^h \rtrivert} 
  \gtrsim \ltrivert e \rtrivert_{\widehat{V}^h},
\end{aligned} \label{equation:boundonenergy}
\end{equation}
with the error in the ultimate expression according to~$e\coloneqq{}\widehat{u}^h-u^h$. The chain of inequalities in~\eqref{equation:boundonenergy} implies that the error estimator~$\mathcal{E}$ controls the error in the extended energy norm,~$\ltrivert e \rtrivert_{\widehat{V}^h}$.

The reason for defining the residual as a map from $V^h$ to~$\widehat{V}^{h\ast}$ is that the stabilization terms in the residual are generally unbounded in the ambient space of the continuum problem, \emph{viz.} $H^1(\Omega)$ for the steady heat equation and $H^1(\Omega,\mathbb{R}^d)\times{}L^2(\Omega)$ for the steady viscous-flow equation. As we will elaborate in Sections~\ref{section:laplaceerror} and~\ref{section:stokeserror}, 
the refined approximation $\widehat{u}^h$ is not required for the calculation of the residual-based estimator \eqref{equation:estimator}. The 
extended space $\widehat{V}^h$ merely serves to establish the error-control relation~\eqref{equation:boundonenergy}.

\subsubsection{Steady heat conduction}\label{section:laplaceerror}
To derive the error indicators for the steady heat conduction problem introduced in Section~\ref{section:laplace}, it is first noted that because of Galerkin orthogonality
\begin{align}
  r^h(\widehat{v}^h) &=  r^h(\widehat{v}^h - \Pi^h \widehat{v}^h) = r^h(\tilde{v}),
\end{align}
where $\tilde{v}=\widehat{v}^h - \Pi^h \widehat{v}^h \in \widehat{V}^h$ and $\Pi^h: \widehat{V}^h \to \Vh$ is an interpolation operator \cite{schumaker2007,veiga2014}. Note that, for notational convenience, we will drop the diacritic and superscript from $\widehat{v}^h \in \widehat{V}^h$ in the remainder of this section, \emph{i.e.}, $\widehat{v}^h=v$. 

Using the definition of the residual \eqref{equation:residualdefinition} in combination with the definitions of the bilinear and linear forms \eqref{eq:ah_laplace} and \eqref{eq:bh_laplace}, (reverse) integration by parts yields
\begin{equation}
\begin{aligned}
	r^{h}(\tilde{v}) &= \sum \limits_{K \in \mesh} \Bigg\{ \int \limits_{K \cap \Omega} r^{h}_{\rm volume}  \tilde{v} \, {\rm d}V +  \int \limits_{K \cap \partial \Omega_{N}} r_{\rm neumann}^{h} \tilde{v} \, {\rm d}S \Bigg. \\
	&\phantom{\sum \limits_{K \in \mesh}} \Bigg. + \int \limits_{K \cap \partial \Omega_{D}} \left( - r_{\rm nitsche}^{h} \right) \partial_{n} \tilde{v} \, {\rm d}S  +  \int \limits_{K \cap \partial \Omega_D} \frac{\beta}{h_{K}} r_{\rm nitsche}^{h} \tilde{v} \, {\rm d}S \Bigg. \\
	&\phantom{\sum \limits_{K \in \mesh}} \Bigg. + \int \limits_{\partial K \cap \mathcal{F}_{\rm skeleton}} \left( - r^{h}_{\rm jump} \right) \tilde{\v} \, {\rm d}S +  \int \limits_{\partial K \cap \mathcal{F}_{\rm ghost}} \gamma_{g} h_{F}^{2\order-1} \left(-r_{\rm ghost}^{h} \right) \llbracket  \partial_{n}^{k} \tilde{v} \rrbracket \, {\rm d}S \Bigg\},
	\label{equation:residualheat}
\end{aligned}
\end{equation}
where
\begin{subequations}
	\begin{align}
		r^{h}_{\rm volume}    &:= f + \Delta u^{h}, \\
		r^{h}_{\rm neumann} &:= q - \partial_{n} u^{h}, \\
		r^{h}_{\rm nitsche}    &:= g - u^{h}, \\
		r^{h}_{\rm jump}        &:= \tfrac{1}{2} \llbracket \partial_{n} u^h \rrbracket ,\\
		r^{h}_{\rm ghost}      &:= \tfrac{1}{2}  \llbracket \partial_{n}^{k} u^{h} \rrbracket.
	\end{align}
\end{subequations}
The factor $\frac{1}{2}$ in the jump and ghost terms accounts for the presence of the associated interfaces in two elements. Using the Cauchy-Schwarz inequality it then follows that
\begin{equation}
\begin{aligned}
	\lvert r^{h}(\tilde{v}) \rvert
	&\leq \sum \limits_{K \in \mesh} \Bigg\{ \lVert r^{h}_{\rm volume} \rVert_{L^2(K \cap \Omega)} \lVert \tilde{v} \rVert_{L^{2}(K \cap \Omega)} +  \lVert r_{\rm neumann}^{h} \rVert_{L^{2}(K \cap \partial \Omega_{N})} \lVert \tilde{v}\rVert_{L^{2}(K \cap \partial \Omega_{N})} \Bigg. \\
	&\phantom{\leq \sum \limits_{K \in \mesh}} \Bigg. + \lVert r_{\rm nitsche}^{h} \rVert_{L^{2}(K \cap \partial \Omega_{D})} \lVert \partial_{n} \tilde{v} \rVert_{L^{2}(K \cap \partial \Omega_{D})} + \beta {\meshsize_{K}^{-1}}  \lVert r_{\rm nitsche}^{h} \rVert_{L^{2}(K \cap \partial \Omega_D)} \lVert \tilde{v} \rVert_{L^{2}(K \cap \partial \Omega_D)} \Bigg. \\
	&\phantom{\leq \sum \limits_{K \in \mesh}} \Bigg. + \sum_{F \in \mathcal{F}_{\rm skeleton}} \lVert r^{h}_{\rm jump} \rVert_{L^2(\partial K \cap F)} \lVert \tilde{v} \rVert_{L^2(\partial K \cap F)} \Bigg. \\
	&\phantom{\leq \sum \limits_{K \in \mesh}} + \sum_{F \in \mathcal{F}_{\rm ghost}}  \gamma_{g} \meshsize_{F}^{2\order-1}  \lVert r_{\rm ghost}^{h} \rVert_{L^{2}(\partial K \cap F)} \lVert \llbracket  \partial_{n}^{k} \tilde{v} \rrbracket \rVert_{L^{2}(\partial K \cap F)} \Bigg\}. \label{ineq:res2}
\end{aligned}
\end{equation}

Using standard interpolation inequalities \cite{ern2013,bazilevs2006} and the definition of the norm \eqref{eq:energynorm}, and noting that we consider the functions $\tilde{v}$ and $v$ to be piecewise polynomials, it follows that
\begin{subequations} \label{eq:interpolations}
\begin{align}
&\lVert \tilde{v}\rVert_{L^{2}(K \cap \Omega)} \lesssim h_K \lVert \nabla v \rVert_{L^2(\widetilde{K})} \lesssim h_K \ltrivert {v} \rtrivert_{\widetilde{K}},\\
&\lVert \tilde{v}\rVert_{L^{2}(K \cap \partial \Omega)} \lesssim h_K^{-\frac{1}{2}}\lVert \tilde{v}\rVert_{L^{2}(K)}  \lesssim h_K^{\frac{1}{2}} \ltrivert v \rtrivert_{\widetilde{K}}, \\ 
&\lVert \partial_n \tilde{v}\rVert_{L^{2}(K \cap \partial \Omega_{D})}  \lesssim h_K^{-\frac{1}{2}}\lVert \nabla \tilde{v}\rVert_{L^{2}(K)}  \lesssim h_K^{-\frac{1}{2}} \ltrivert v \rtrivert_{\widetilde{K}}, \\
&\lVert \tilde{v} \rVert_{L^2(\partial K \cap F)} 
\lesssim h_K^{-\frac{1}{2}}\lVert \tilde{v}\rVert_{L^{2}(K)}  \lesssim h_K^{\frac{1}{2}} \ltrivert v \rtrivert_{\widetilde{K}}
,\\
&\lVert \llbracket \partial_n^k \tilde{v} \rrbracket \rVert_{L^2(\partial K \cap F)} \lesssim \lVert \partial_n^k \tilde{v}_K \rVert_{L^2(\partial K \cap F)} + \lVert \partial_n^k \tilde{v}_{K'} \rVert_{L^2(\partial K' \cap F)}  \lesssim h_F^{\frac{1}{2}-k} \ltrivert v \rtrivert_{\widetilde{K}\cup \widetilde{K}'},
\end{align}
\label{eq:interpolationinequalities}%
\end{subequations}
where $\widetilde{K}$ is the support extension \cite{bazilevs2006} of the element $K$ and $K'$ is the element that shares the interface $F$ with element $K$. The residual can then be bounded as
\begin{equation}
\begin{aligned}
	\lvert r^{h}(\tilde{v}) \rvert
	&\lesssim \sum \limits_{K \in \mesh} \Bigg\{ h_K \lVert  r^{h}_{\rm volume} \rVert_{L^2(K \cap \Omega)} +  h_K^{\frac{1}{2}}\lVert r_{\rm neumann}^{h} \rVert_{L^{2}(K \cap \partial \Omega_{N})} \Bigg. \\
	&\phantom{\lesssim \sum \limits_{K \in \mesh} } \Bigg. + h_K^{-\frac{1}{2}} \lVert r_{\rm nitsche}^{h} \rVert_{L^{2}(K \cap \partial \Omega_{D})} + \beta \meshsize_{K}^{-\frac{1}{2}}  \lVert r_{\rm nitsche}^{h} \rVert_{L^{2}(K \cap \partial \Omega_D)} \Bigg. \\
	&\phantom{\lesssim \sum \limits_{K \in \mesh} } \Bigg. + \sum_{F \in \mathcal{F}_{\rm skeleton}} h_F^{\frac{1}{2}} \lVert r^{h}_{\rm jump} \rVert_{L^2(\partial K \cap F)} + \sum_{F\in \mathcal{F}_{\rm ghost}}\gamma_{g} \meshsize_{F}^{\order-\frac{1}{2}} \lVert r_{\rm ghost}^{h} \rVert_{L^{2}(\partial K \cap F)}  \Bigg\} \ltrivert {v} \rtrivert_{\widetilde{K}\cup \widetilde{K}'},
 \label{ineq:res3}
\end{aligned}
\end{equation}
which, using the discrete Cauchy-Schwarz inequality can be rewritten as
\begin{equation}
\begin{aligned}
 \frac{\lvert r^{h}({v}) \rvert^{2}}{\ltrivert {v} \rtrivert^2}
	&\lesssim \sum \limits_{K \in \mesh} \Bigg\{ h_K^2 \lVert  r^{h}_{\rm volume} \rVert_{L^2(K \cap \Omega)}^{2}  +  h_K\lVert r_{\rm neumann}^{h} \rVert_{L^{2}(K \cap \partial \Omega_{N})}^{2} \Bigg. \\
	&\phantom{\lesssim \sum \limits_{K \in \mesh}} \Bigg. + h_K^{-1} \lVert r_{\rm nitsche}^{h} \rVert_{L^{2}(K \cap \partial \Omega_{D})}^{2}  + \beta^2 \meshsize_{K}^{-1}  \lVert r_{\rm nitsche}^{h} \rVert_{L^{2}(K \cap \partial \Omega_D)}^{2} \Bigg. \\
	&\phantom{\lesssim \sum \limits_{K \in \mesh}} \Bigg. + \sum_{F\in \mathcal{F}_{\rm skeleton}}h_F \lVert r^{h}_{\rm jump} \rVert_{L^2(\partial K \cap F)}^{2} + \sum_{F\in \mathcal{F}_{\rm ghost}}\gamma_{g}^{2} \meshsize_{F}^{2\order-1} \lVert r_{\rm ghost}^{h} \rVert_{L^{2}(\partial K \cap F)}^{2}  \Bigg\} .
 \label{ineq:res}
\end{aligned}
\end{equation}
Using the definition of the residual norm \eqref{equation:residualdefinition} it follows that
\begin{align}
 \| r^h \|_{\widehat{V}^h} \lesssim \mathcal{E} = \sqrt{ \sum \limits_{K \in \mesh}\eta_{K}^{2}},
\end{align}
with the element error indicators defined as
\begin{equation}
\begin{aligned}
 \eta_K^2
	&:= \Big. h_K^2 \lVert  r^{h}_{\rm volume} \rVert_{L^2(K \cap \Omega)}^{2}  +  h_K\lVert r_{\rm neumann}^{h} \rVert_{L^{2}(K \cap \partial \Omega_{N})}^{2} \Big. \\
	&\phantom{+} \Big. + h_K^{-1} \lVert r_{\rm nitsche}^{h} \rVert_{L^{2}(K \cap \partial \Omega_{D})}^{2}  + \beta^2 \meshsize_{K}^{-1}  \lVert r_{\rm nitsche}^{h} \rVert_{L^{2}(K \cap \partial \Omega_D)}^{2} \Big. \\
	&\phantom{+} \Big. + \sum_{F\in \mathcal{F}_{\rm skeleton}}h_F \lVert r^{h}_{\rm jump} \rVert_{L^2(\partial K \cap F)}^{2} + \sum_{F\in \mathcal{F}_{\rm ghost}}\gamma_{g}^{2} \meshsize_{F}^{2\order-1} \lVert r_{\rm ghost}^{h} \rVert_{L^{2}(\partial K \cap F)}^{2}. \Big.
 \label{equation:laplaceindicators}
\end{aligned}
\end{equation}
This error indicator reflects that the total element error for all elements that do not intersect the boundary of the domain is composed of the interior residual and the residual term for the jump in the solution normal derivative across the element interfaces. It is noted that for higher-order continuous discretizations, \emph{i.e.}, $\regularity > 0$, the jump contribution vanishes. For elements that intersect the Neumann boundary, additional error contributions are obtained from the Neumann residual and the ghost penalty residual, while additional Nitsche-related contributions appear for elements intersecting the Dirichlet boundary.

\subsubsection{Steady viscous flow}\label{section:stokeserror}
For the Stokes problem introduced in Section~\ref{section:stokes}, using (reverse) integration by parts, the error indicators in equation \eqref{equation:estimator} are obtained by considering the residual \eqref{equation:residualdefinition} as
\begin{equation}
\begin{aligned}
 \Rh(\tilde{\v})
 &= \sum \limits_{\element \in \mesh} \Bigg\{ \int \limits_{\element \cap \Omega} \rintu 
  \cdot \tilde{\vv} + \rintp \tilde{\q} \dO +  \int \limits_{ \element  \cap \neumannboundary} \rneumannu \cdot \tilde{\vv} \dG \Bigg. \\
&\phantom{=\sum \limits_{\element \in \cutmesh}} \Bigg. +
\int \limits_{ \element \cap \dirichletboundary} (-\rnitscheu) \cdot \left( \left( 2 \mu \nabla^s \tilde{\vv} \right) \nn + \tilde{\q} \nn \right) \dG + \int \limits_{\element \cap \partial \Omega_D} \frac{\beta \mu}{\meshsize_{K}} \rnitscheu \cdot \tilde{\vv} \dG \Bigg. \\
&\phantom{=\sum \limits_{\element \in \cutmesh}} \Bigg. + \sum_{F \in \mathcal{F}_{\rm skeleton}}\int \limits_{\partial \element \cap F} (-\rjumpu) \cdot \tilde{\vv} \dG \Bigg. \\
&\phantom{=\sum \limits_{\element \in \cutmesh}} \Bigg. + \sum \limits_{\face \in \ghost} \int \limits_{\partial \element \cap \face} \gamma_{g} \mu  \meshsize_{\face}^{2\order-1} (-\rghostu) \cdot \llbracket  \gradn^\order \tilde{\vv} \rrbracket \dG \Bigg. \\
&\phantom{=\sum \limits_{\element \in \cutmesh}} \bigg.  + \sum \limits_{\face \in \skeleton} \int \limits_{\partial \element \cap \face}  \frac{ {\gamma_{s} }\meshsize_{\face}^{2\order+1}}{\mu} \rskeleton \llbracket  \gradn^\order \tilde{\q} \rrbracket \dG  \Bigg\},
\label{equation:residualstokes}
\end{aligned}
\end{equation}
where $\tilde{\v}=\v-\Pi^h \v=(\tilde{\vv},\tilde{\q})$ and
{\allowdisplaybreaks
\begin{subequations}
\begin{align}
 \rintu &:= \bodystokes + \nabla \cdot \left( 2 \mu \nabla^s \uuh \right) - \nabla \ph,\\
 \rintp &:= \nabla \cdot \uuh,\\
 \rneumannu &:= \neumanndatastokes - \left( 2 \mu \nabla^s \uuh \right) \nn + \ph \nn ,\\
 \rnitscheu &:= \dirichletdatastokes - \uuh ,\\
 \rjumpu &:= \tfrac{1}{2}  \llbracket \left( 2 \mu \nabla^s \uuh \right) \nn \rrbracket,\\
 \rghostu &:= \tfrac{1}{2}  \llbracket \gradn^k \uuh \rrbracket,\\
 \rskeleton &:= \tfrac{1}{2} \llbracket \gradn^k \ph \rrbracket.
\end{align}
\end{subequations}
}
Application of the Cauchy-Schwarz inequality gives
\begin{equation}
\begin{aligned}
	\left| \Rh(\tilde{\v}) \right| &\leq \sum \limits_{\element \in \mesh} \Bigg\{ \| \rintu \|_{L^2(\element \cap \domain)}
	\| \tilde{\vv} \|_{L^2(\element \cap \domain)} + \| \rintp \|_{L^2(\element\cap \domain)} \| \tilde{\q} \|_{L^2(\element\cap \domain)} \Bigg. \\
	 &\phantom{\leq  \sum \limits_{\element \in \mesh}} \Bigg. + \| \rneumannu \|_{L^2({ \element \cap \neumannboundary})} \| \tilde{\vv} \|_{L^2({ \element \cap \neumannboundary})} \Bigg. \\
	&\phantom{\leq  \sum \limits_{\element \in \mesh}}\Bigg. + \| \rnitscheu \|_{L^2({\element \cap \dirichletboundary})} \left( 2 \mu \| \left( \nabla^s \tilde{\vv} \right) \nn \|_{L^2({ \element \cap \dirichletboundary})} + \| \tilde{\q} \|_{L^2({\element \cap \dirichletboundary})}\right) \Bigg. \\
	&\phantom{\leq  \sum \limits_{\element \in \mesh}} \Bigg. +  \mu \beta \meshsize_{\element}^{-1}  \| \rnitscheu \|_{L^2(\element \cap \partial \domain_D)} \| \tilde{\vv} \|_{L^2(\element \cap \partial \domain_D)} \Bigg. \\
	&\phantom{\leq  \sum \limits_{\element \in \mesh}} \Bigg. + \sum_{\face \in \mathcal{F}_{\rm skeleton}}\| \rjumpu \|_{L^2({\partial \element \cap F})}  \| \tilde{\vv} \|_{L^2({\partial \element \cap \face})} \Bigg. \\
	&\phantom{\leq  \sum \limits_{\element \in \mesh}} \Bigg. + \sum \limits_{\face \in \ghost^\element} \mu \gamma_g \meshsize_{\face}^{2\order-1} \| \rghostu \|_{L^2(\partial \element \cap \face)} \| \llbracket  \gradn^\order \tilde{\vv} \rrbracket \|_{L^2(\partial \element \cap \face)}  \Bigg. \\
	&\phantom{\leq  \sum \limits_{\element \in \mesh}}+\sum \limits_{\face \in \skeleton^\element} \mu^{-1} \gamma_s\meshsize_{\face}^{2\order+1}  \| \rskeleton \|_{L^2(\partial \element \cap \face)}  \|  \llbracket  \gradn^\order \tilde{\q} \rrbracket \|_{L^2(\partial \element \cap \face)} \Bigg\},
\end{aligned}
\end{equation}
which, using the inequalities \eqref{eq:interpolationinequalities} and
{\allowdisplaybreaks
\begin{subequations}
\begin{align}
 &\| \tilde{q} \|_{L^2(K \cap \Omega)} \lesssim  \| q \|_{L^2(K)},\\
 &\| \tilde{q} \|_{L^2(K \cap \partial \Omega)} \lesssim h_K^{-\frac{1}{2}} \| \tilde{q} \|_{L^2(K)} \lesssim h_K^{-\frac{1}{2}} \| {q} \|_{L^2(K)},\\
 &\lVert \llbracket \partial_n^k \tilde{q} \rrbracket \rVert_{L^2(\partial K \cap F)} \lesssim \lVert \partial_n^k \tilde{q}_K \rVert_{L^2(\partial K \cap F)} + \lVert \partial_n^k \tilde{q}_{K'} \rVert_{L^2(\partial K' \cap F)}  \lesssim h_F^{-\frac{1}{2}-k} \| q \|_{L^2(\widetilde{K}\cup \widetilde{K}')},
\end{align}
\end{subequations}
}
can be rewritten as
\begin{equation}
\begin{aligned}
	\left| \Rh(\tilde{\v}) \right| &\lesssim \sum \limits_{\element \in \mesh} \Bigg\{ \mu^{-\frac{1}{2}} h_K \| \rintu \|_{L^2(\element \cap \domain)}
	+ \mu^{\frac{1}{2}} \| \rintp \|_{L^2(\element\cap \domain)} + \mu^{-\frac{1}{2}} h_K^{\frac{1}{2}} \| \rneumannu \|_{L^2({ \element \cap \neumannboundary})} \Bigg. \\
	&\phantom{\leq  \sum \limits_{\element \in \mesh}}\Bigg. + 3 \mu^{\frac{1}{2}} h_K^{-\frac{1}{2}}  \| \rnitscheu \|_{L^2({\element \cap \dirichletboundary})} +  \mu^{\frac{1}{2}} \beta \meshsize_{\element}^{-\frac{1}{2}} \| \rnitscheu \|_{L^2(\element \cap \partial \domain_D)} 
	\Bigg. \\
	&\phantom{\leq  \sum \limits_{\element \in \mesh}} \Bigg. + \sum_{\face \in \mathcal{F}_{\rm skeleton}} \mu^{-\frac{1}{2}} h_K^{\frac{1}{2}} \| \rjumpu \|_{L^2({\partial \element \cap F})}  \Bigg. \\
	&\phantom{\leq  \sum \limits_{\element \in \mesh}} \Bigg. + \sum \limits_{\face \in \ghost^\element}  \mu^{\frac{1}{2}} \gamma_g \meshsize_{\face}^{\order-\frac{1}{2}}  \| \rghostu \|_{L^2(\partial \element \cap \face)}
\Bigg.  \\
	&\phantom{\leq  \sum \limits_{\element \in \mesh}}+\sum \limits_{\face \in \skeleton^\element}  \mu^{-\frac{1}{2}} \gamma_s \meshsize_{\face}^{\order+\frac{1}{2}}  \| \rskeleton \|_{L^2(\partial \element \cap \face)}  
	 \Bigg\} \ltrivert {v} \rtrivert_{\widetilde{K}\cup \widetilde{K}'}.
\end{aligned}
\end{equation}
Note that the factor 3 in front of the Nitsche residual results from the fact that both terms $2 \mu \| \left( \nabla^s \tilde{\vv} \right) \nn \|_{L^2({ \element \cap \dirichletboundary})}$ and $\| \tilde{\q} \|_{L^2({\element \cap \dirichletboundary})}$ are bound by the same norm. Following the same steps as for the heat conduction problem we then obtain the element error indicators as
\begin{equation}
\begin{aligned}
	\eta_K^2 &=  \Bigg. \mu^{-1} h_K^2 \| \rintu \|_{L^2(\element \cap \domain)}^2
	+ \mu \| \rintp \|_{L^2(\element\cap \domain)}^2 + \mu^{-1} h_K \| \rneumannu \|_{L^2({ \element \cap \neumannboundary})}^2 \Bigg. \\
	& \Bigg. + 9 \mu h_K^{-1}  \| \rnitscheu \|_{L^2({\element \cap \dirichletboundary})}^2 +  \mu \beta^2 \meshsize_{\element}^{-1} \| \rnitscheu \|_{L^2(\element \cap \partial \domain_D)}^2 \Bigg. \\
	& \Bigg. + \sum_{\face \in \mathcal{F}_{\rm skeleton}} \mu^{-1} h_K \| \rjumpu \|_{L^2({\partial \element \cap F})}^2 \Bigg. \\
	& \Bigg. + \sum \limits_{\face \in \ghost^\element}  \mu \gamma_g^2 \meshsize_{\face}^{2\order-1}  \| \rghostu \|_{L^2(\partial \element \cap \face)}^2 \Bigg. \\
	& \Bigg. \sum \limits_{\face \in \skeleton^\element}  \mu^{-1} \gamma_s^2 \meshsize_{\face}^{2\order+1}  \| \rskeleton \|_{L^2(\partial \element \cap \face)}^2. \Bigg.
\end{aligned}
\end{equation}
Compared to the error indicators for the heat conduction problem, we here get one additional term to represent the error in the balance of mass, \emph{i.e.}, $\| \rintp \|_{L^2(\element\cap \domain)}$, and one term related to the skeleton-stabilization, \emph{i.e.}, $\| \rskeleton \|_{L^2(\partial \element \cap \face)}$. Moreover, note that the mass and momentum balance terms are scaled with $\mu^{-\frac{1}{2}}$ and $\mu^{\frac{1}{2}}$, respectively, in order to be dimensionally-consistent with the energy norm \eqref{eq:energynormstokes}.

\subsection{Adaptive solution procedure}\label{section:eeaprocedure}
We employ the residual-based error estimator introduced above in an iterative mesh refinement procedure. In each iteration, for the given mesh we solve the Galerkin problem \eqref{equation:galerkin} and subsequently compute the element-wise error indicators \eqref{equation:estimator} (and the corresponding estimator). Based on the indicators, certain elements are then refined, after which the procedure is repeated on the refined mesh. These iterations are continued until a stopping criterion is satisfied.

We consider D\"orfler marking \cite{dorfler1996} to select the elements to be refined. In this marking strategy, the marked set, $\mathcal{M}$, is defined as a minimal set of elements such that
\begin{align}
 \sqrt{ \sum_{K \in \mathcal{M}} \eta_{K}^2} \geq \lambda  \sqrt{\sum_{K\in \mesh} \eta_{K}^2} = \lambda  \mathcal{E},
\end{align}
with $\lambda$ a selected fraction of the error estimator. For the considered (truncated) hierarchical spline meshes, refining elements does not necessarily result in a refinement of the approximation space \cite{kuru2014, brummelen2020}. To ensure that the approximation space is refined, an additional step is required in which a refinement mask $\widetilde{\mathcal{M}} \supset \mathcal{M}$ is defined. To determine the refinement mask, for each element $K$ in the marked set $\mathcal{M}$ we determine the support extension 
\begin{align}
 \widetilde{\mathcal{K}} =  \bigcup \left\{ {\rm supp}(N) \mid {\rm supp}(N) \cap K \neq \emptyset,~N \in \mathcal{H}(\mathcal{T}) \right\},
\end{align}
and then refine the elements in each support extension which are not smaller than the element $K$, \emph{i.e.},
\begin{align}
 \widetilde{\mathcal{M}} = \bigcup \limits_{K \in \mathcal{M}} \left\{ K' \in \widetilde{\mathcal{K}}   \mid K'\in \cup_{\ell=0}^{\ell_K} \mathcal{T}^\ell,~K \in \mathcal{T}^{\ell_K} \right\}.
\end{align}

\begin{figure}[htp]
	\centering
	\begin{subfigure}[b]{0.28\textwidth}
		\centering
		\includegraphics[width=\textwidth]{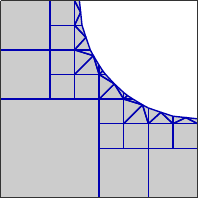}
		\caption{Original element}
	\end{subfigure}	\hfill%
	\begin{subfigure}[b]{0.28\textwidth}
		\centering
		\includegraphics[width=\textwidth]{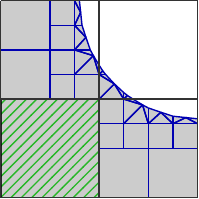}
		\caption{First refinement}
	\end{subfigure}	\hfill%
	\begin{subfigure}[b]{0.28\textwidth}
		\centering
		\includegraphics[width=\textwidth]{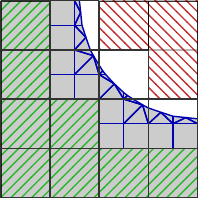}
		\caption{Second refinement}
	\end{subfigure}%

\vspace{5mm}

		\centering
		\includegraphics[width=0.6\textwidth]{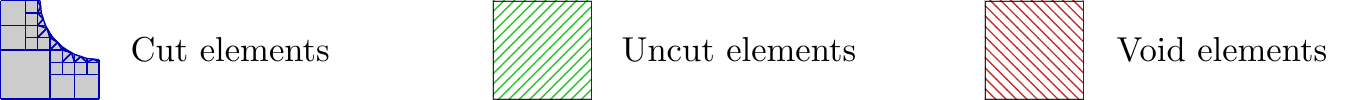}
\vspace{5mm}		
	\caption{Illustration of the refinement procedure for cut elements. The original element is subdivided in integration subcells (blue borders) using the recursive bisectioning procedure detailed in Ref.~\cite{divi2020}. At the lowest level of bisectioning, a triangulation procedure is employed. After one refinement of the original element, the original element is split into 4 elements, of which one is now an uncut element and the other three are cut elements. The bisectioning depth for the determination of the integration subcells is reduced by one level compared to the original element, so that the subcells remain identical under the element refinement operation. After one further refinement step, each of the four elements in the first refinement is further refined, resulting now also in elements that are void and are hence discarded from the background mesh.}
	\label{fig:geometrydef}
\end{figure}

During the element refinement procedure the geometry approximation is not altered, as illustrated in Figure~\ref{fig:geometrydef}. In our implementation, the bisectioning depth used to determine the integration subcells is lowered under refinement, resulting in the preservation of the integration subcells under refinement. This ensures that the boundary of the segmented geometry is invariant under mesh refinement. A consequence of this choice is that an element can only be refined up to the level of the integration subcells. Elements requiring refinement beyond the level of the integration subcells are discarded from the refinement list, and the adaptive refinement procedure is stopped if there are no more elements that can be refined.

%% file: sections/results.tex
\section{Benchmark simulations}\label{sec:results}
In this section we assess the developed residual-based adaptive refinement technique on a range of numerical experiments. For both the heat conduction problem (Section~\ref{sec:laplaceexamples}) and the viscous flow problem (Section~\ref{sec:stokesexamples}), both singular and non-singular test cases are considered. For all simulations exact reference solutions are available, allowing for a rigorous study of the stability and accuracy of the developed adaptive immersed isogeometric analysis framework. For all simulations the octree subdivision depth is set equal to the desired maximum number of refinements (see Section~\ref{section:eeaprocedure}) and the refinement threshold is set to $\lambda=0.8$. Throughout this section, the problems are considered to be in dimensionless form.


\subsection{Steady heat conduction} \label{sec:laplaceexamples}
We consider the two-dimensional heat conduction problem on a unit square and on a star-shaped domain with a smooth exact solution, and on a domain with a re-entrant corner, for which the exact solution has a reduced regularity (Section~\ref{sec:laplacelshape}). The problems are discretized with linear ($\order=1$) and quadratic ($\order=2$) (TH)B-splines using both uniform and adaptive refinement. All examples consider a non-conforming ambient mesh positioned at an angle of $20$ degrees (see Figure~\ref{fig:square} and Figure~\ref{fig:Lshape}), unless specified otherwise. The empirically selected Nitsche and ghost penalty parameters are set to $\beta=50$ and $\gamma_g=10^{-(\order+2)}$, respectively.

\subsubsection{Unit square}\label{sec:laplaceunitsquare}
\begin{figure}
	\centering
	\begin{subfigure}[b]{0.45\textwidth}
		\centering
		\includegraphics[width=\textwidth]{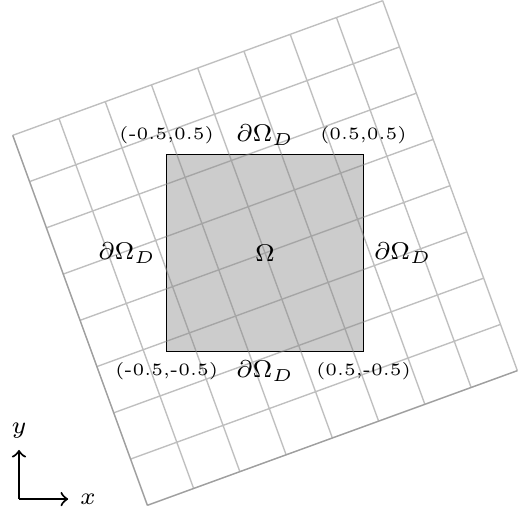}
		\caption{}
		\label{fig:square}
	\end{subfigure}
	\hspace{0.08\textwidth}
	\begin{subfigure}[b]{0.45\textwidth}
		\centering
		\includegraphics[width=\textwidth]{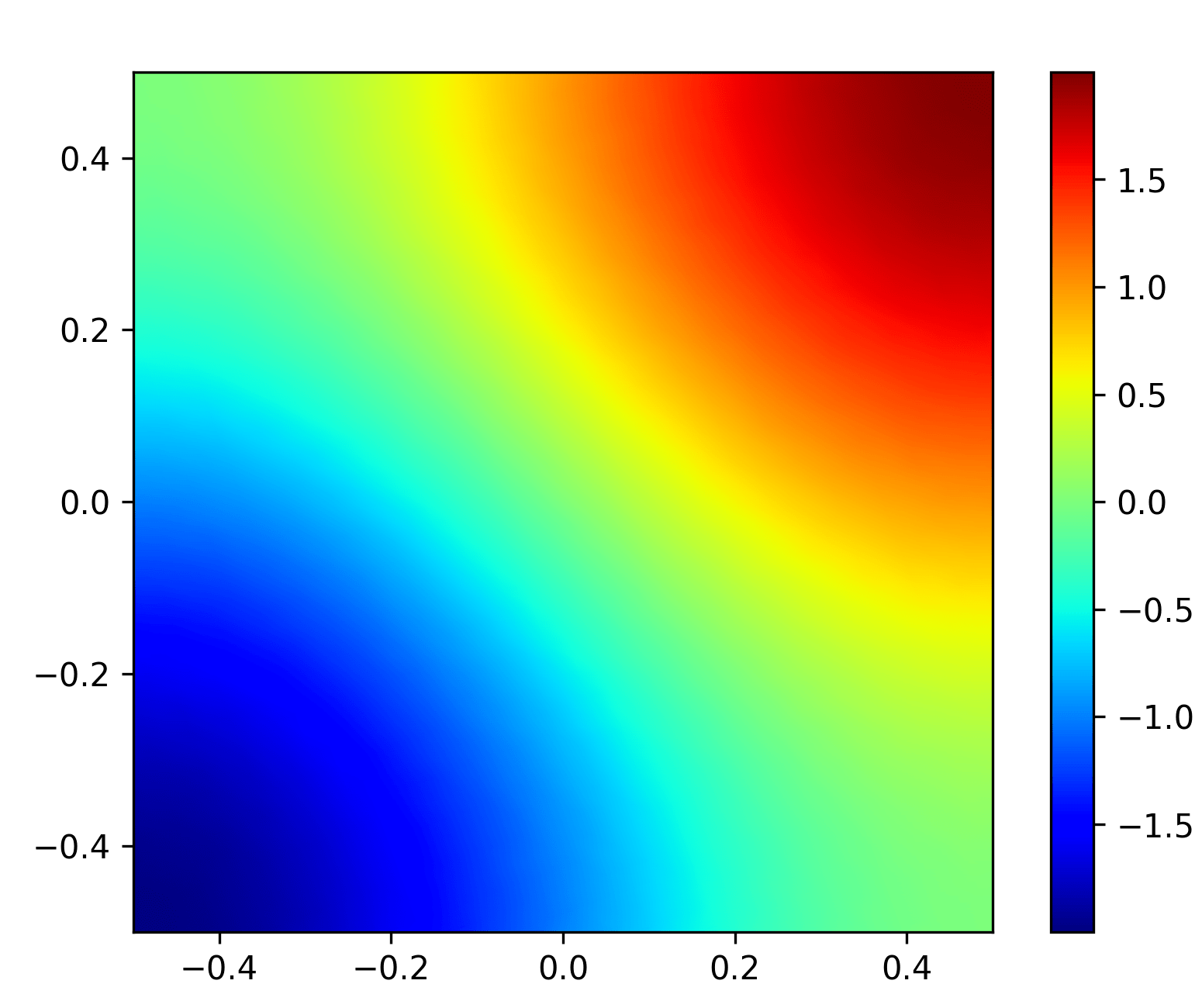}
		\caption{}
		\label{fig:square_exact}
	\end{subfigure}
	\hspace{0.08\textwidth}
	\caption{(\subref{fig:square}) Problem setup, and (\subref{fig:square_exact}) the exact solution $u(x,y)$, Eq.~\eqref{equation:exactlaplace}, for the Laplace problem on the unit square domain.}
	\label{fig:square_laplace}
\end{figure}

Let $\Omega = [-\frac{1}{2},\frac{1}{2}]^2$ be a unit square with Dirichlet boundary $\partial \Omega_{D}$ (see Figure \ref{fig:square}). We define the exact solution of the problem \eqref{equation:stronglaplace} as
\begin{equation}\label{equation:exactlaplace}
	u(x,y) = \sin(\pi x) + \sin(\pi y),
\end{equation}
which is shown in Figure \ref{fig:square_exact}. The heat source $f$ corresponding to this exact solution is equal to zero, and the Dirichlet data is set to $g = u |_{\partial \Omega_{D}}$, matching the exact solution.
 
\begin{figure}
	\centering
	\begin{subfigure}[b]{0.45\textwidth}
		\centering
		\includegraphics[width=\textwidth]{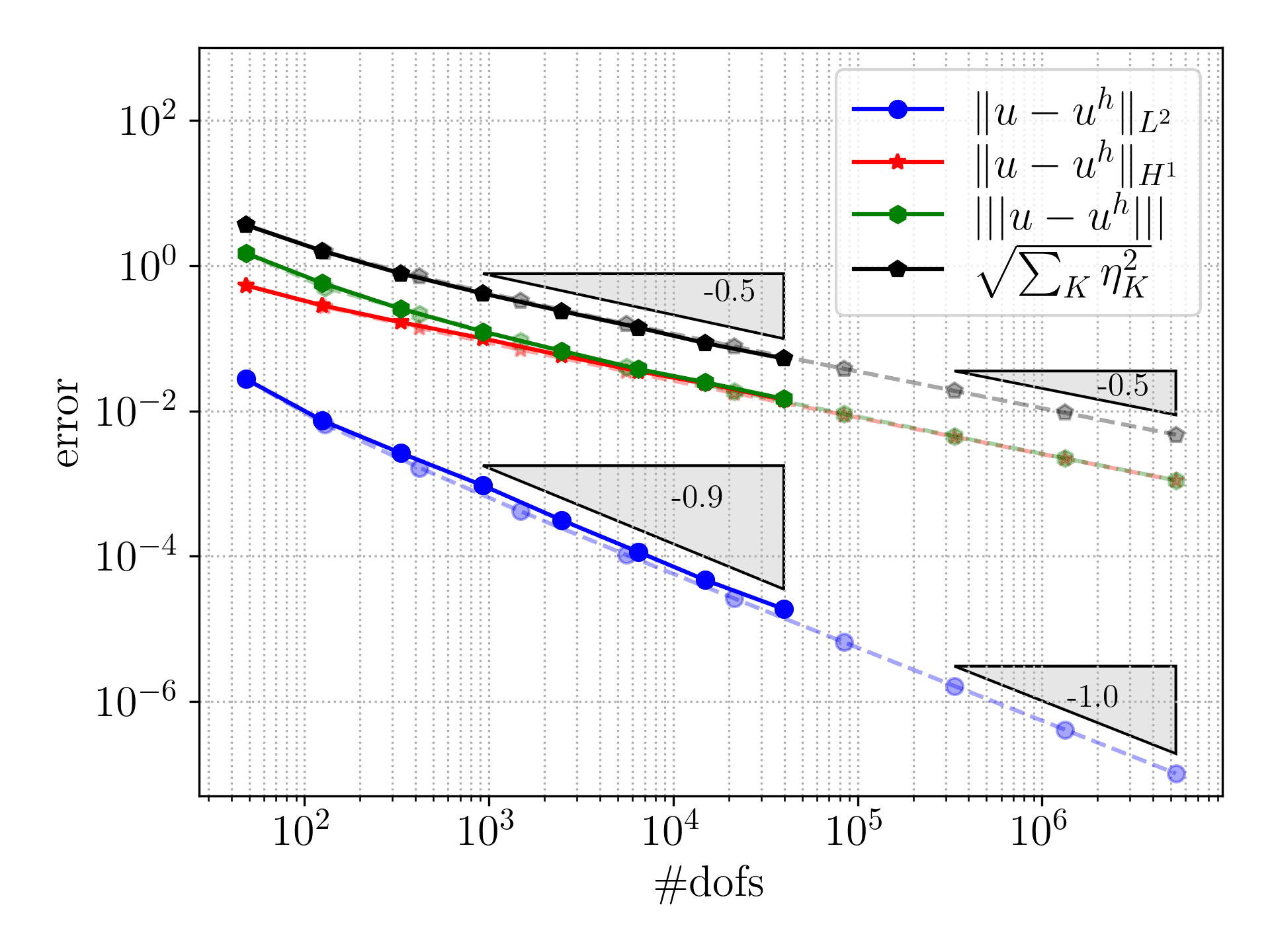}
		\caption{$\order = 1$}
		\label{fig:square_p1}
	\end{subfigure}
	\hspace{0.08\textwidth}
	\begin{subfigure}[b]{0.45\textwidth}
		\centering
		\includegraphics[width=\textwidth]{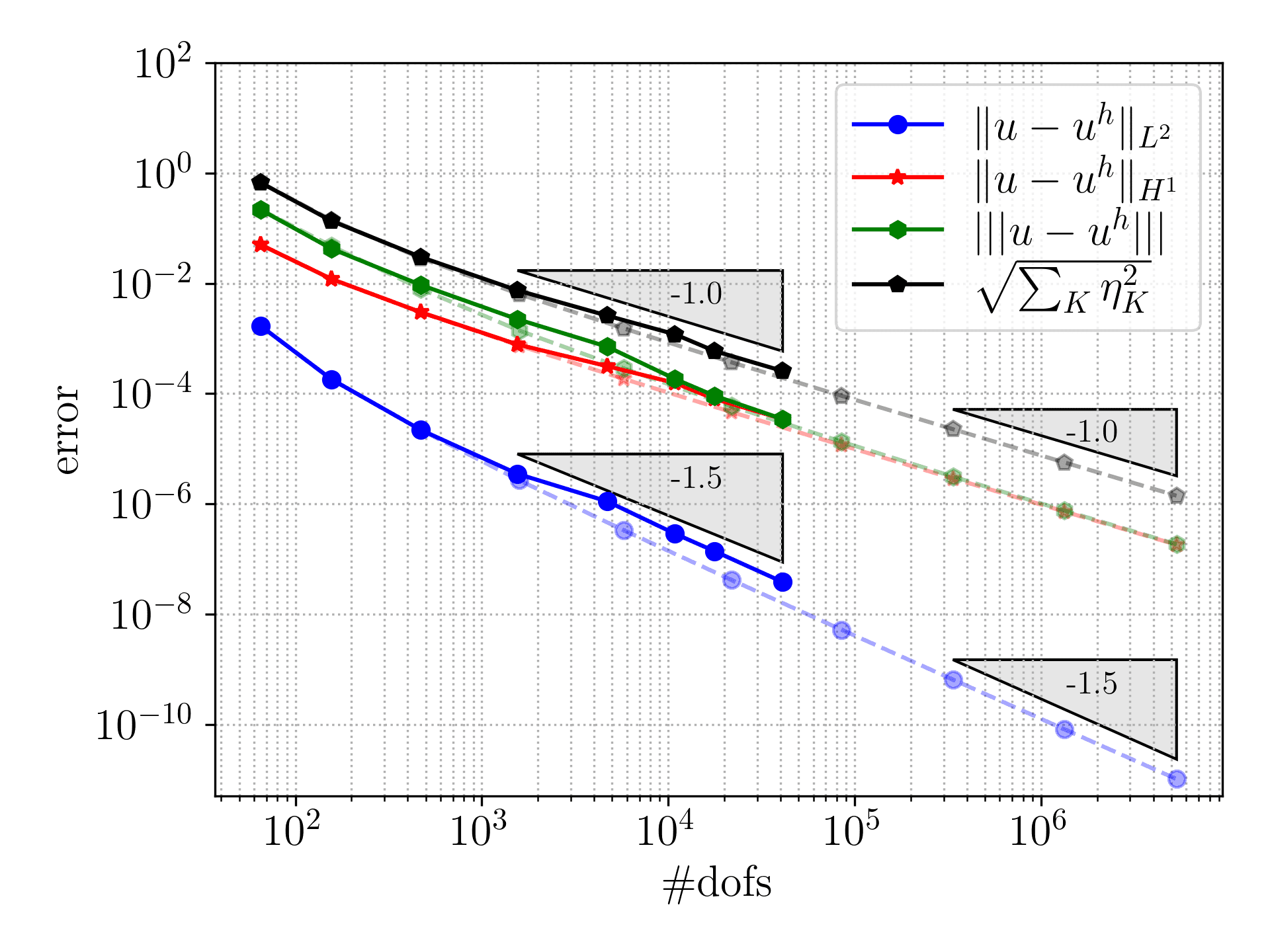}
		\caption{$\order = 2$}
		\label{fig:square_p2}
	\end{subfigure}
	\hspace{0.08\textwidth}
	\caption{Error convergence results for the Laplace problem on the unit square domain under residual-based adaptive refinement (solid) and uniform refinement (dashed) for linear ($\order =1$) and quadratic ($\order =2$) basis functions.}
	\label{fig:square_conv}
\end{figure}

Figure~\ref{fig:square_conv} shows error-analysis results using both uniform and adaptive refinements for the linear case (Figure~\ref{fig:square_p1}) and for the quadratic case (Figure~\ref{fig:square_p2}). Both refinement procedures start from an initial mesh consisting of $8 \times 8$ elements covering the ambient domain $[-1,1]^2$. Optimal convergence rates are obtained for both the error in the $L^2$-norm (\emph{i.e.}, $\mathcal{O}(n^{-\frac{1}{2}(k+1)})$) and in the $H^1$-norm (\emph{i.e.}, $\mathcal{O}(n^{-\frac{1}{2}k})$), with $n$ denoting the number of degrees of freedom. Moreover, as the number of refinement steps increases, the energy norm and $H^1$-norm of the error coincide, indicating that the error is dominated by the $H^1$-semi-norm contribution in Eq.~\eqref{eq:energynorm}. The estimator \eqref{equation:estimator} is observed to converge at the same rate as the energy norm, bounding the energy norm from above, consistent with~Eq.~\eqref{equation:boundonenergy}. Because of the smooth solution \eqref{equation:exactlaplace}, the refinement pattern following from the adaptive refinement procedure closely resembles the uniform refinements, as observed from the close correspondence between the error results for the uniform and adaptive simulations in Figure~\ref{fig:square_conv}.

\subsubsection{Star-shaped domain}\label{sec:laplacestar}
\begin{figure}
	\centering
	\begin{subfigure}[b]{0.45\textwidth}
		\centering
		\includegraphics[width=0.85\textwidth]{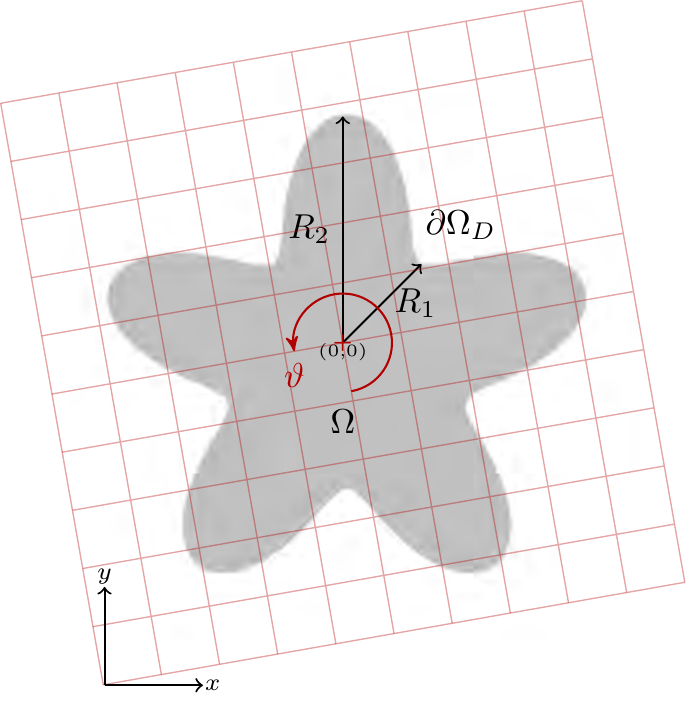}
		\caption{Star shaped domain}
		\label{fig:star_domain}
	\end{subfigure}%
	\hspace{0.08\textwidth}
	\begin{subfigure}[b]{0.45\textwidth}
		\centering
		\includegraphics[width=0.85\textwidth]{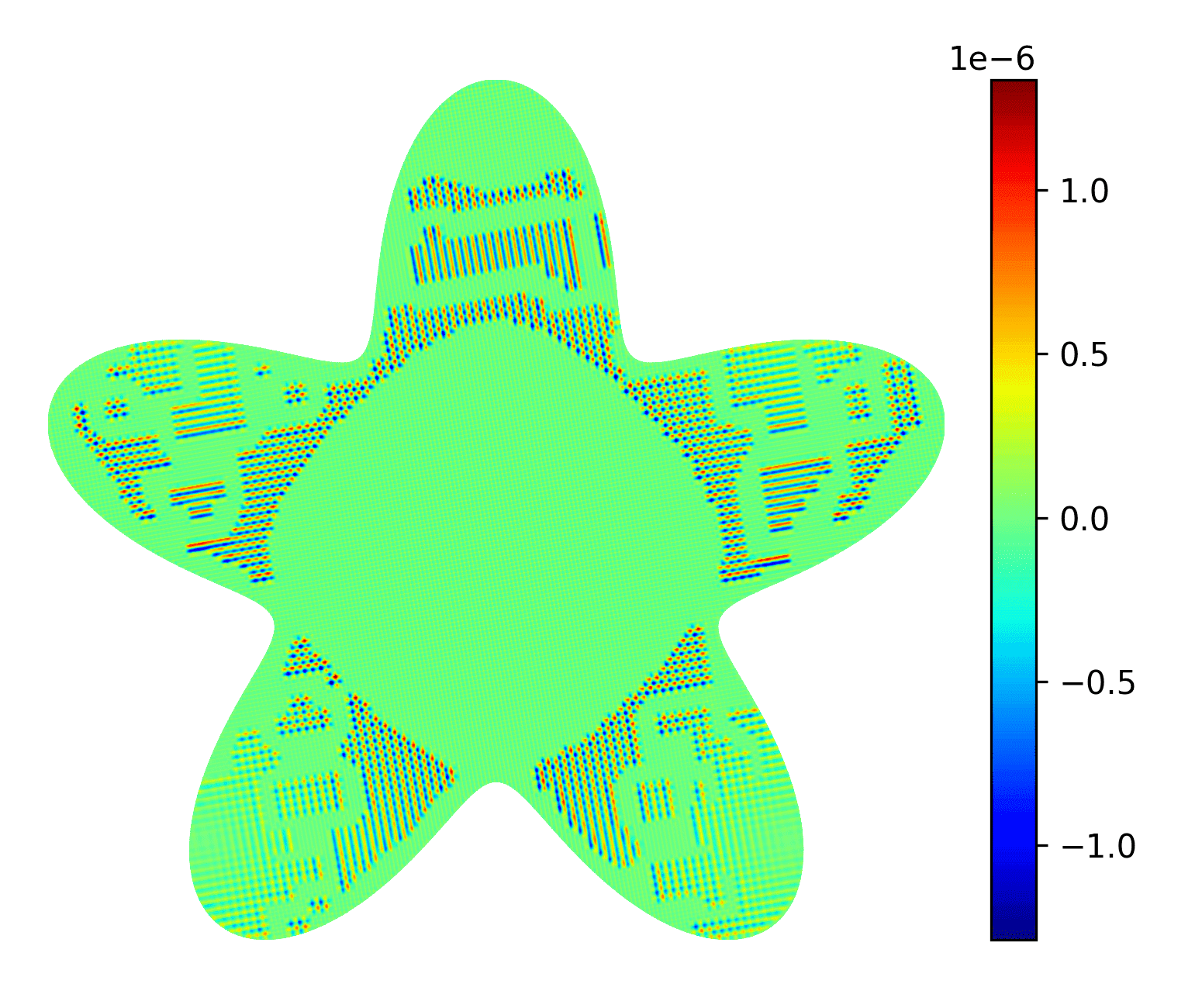}
		\caption{$\vartheta = 10 $}
		\label{fig:star_theta10}
	\end{subfigure}\\[12pt]
	\begin{subfigure}[b]{0.45\textwidth}
		\centering
		\includegraphics[width=0.85\textwidth]{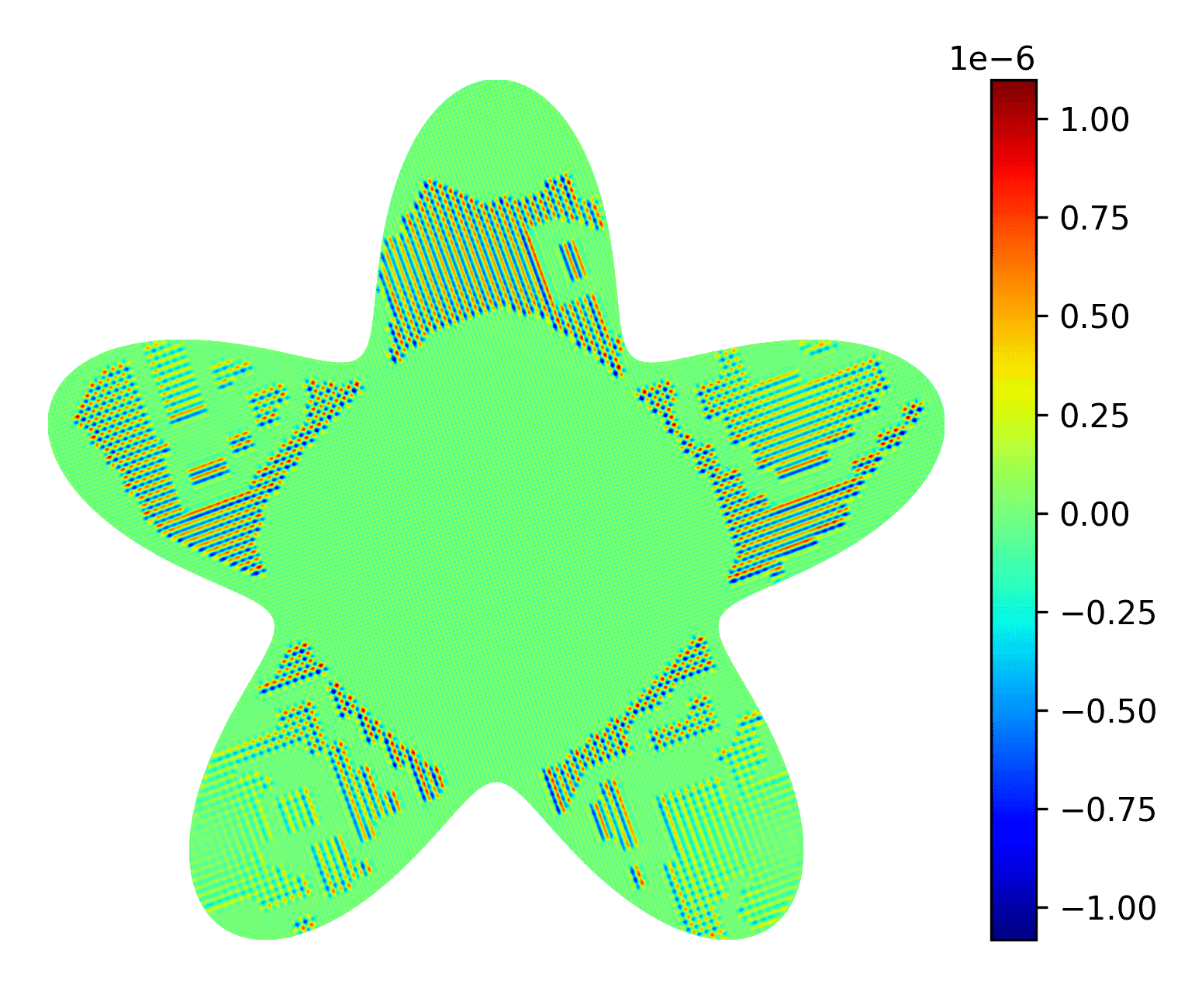}
		\caption{$\vartheta = 20 $}
		\label{fig:star_theta20}
	\end{subfigure}%
	\hspace{0.08\textwidth}
	\begin{subfigure}[b]{0.45\textwidth}
		\centering
		\includegraphics[width=0.85\textwidth]{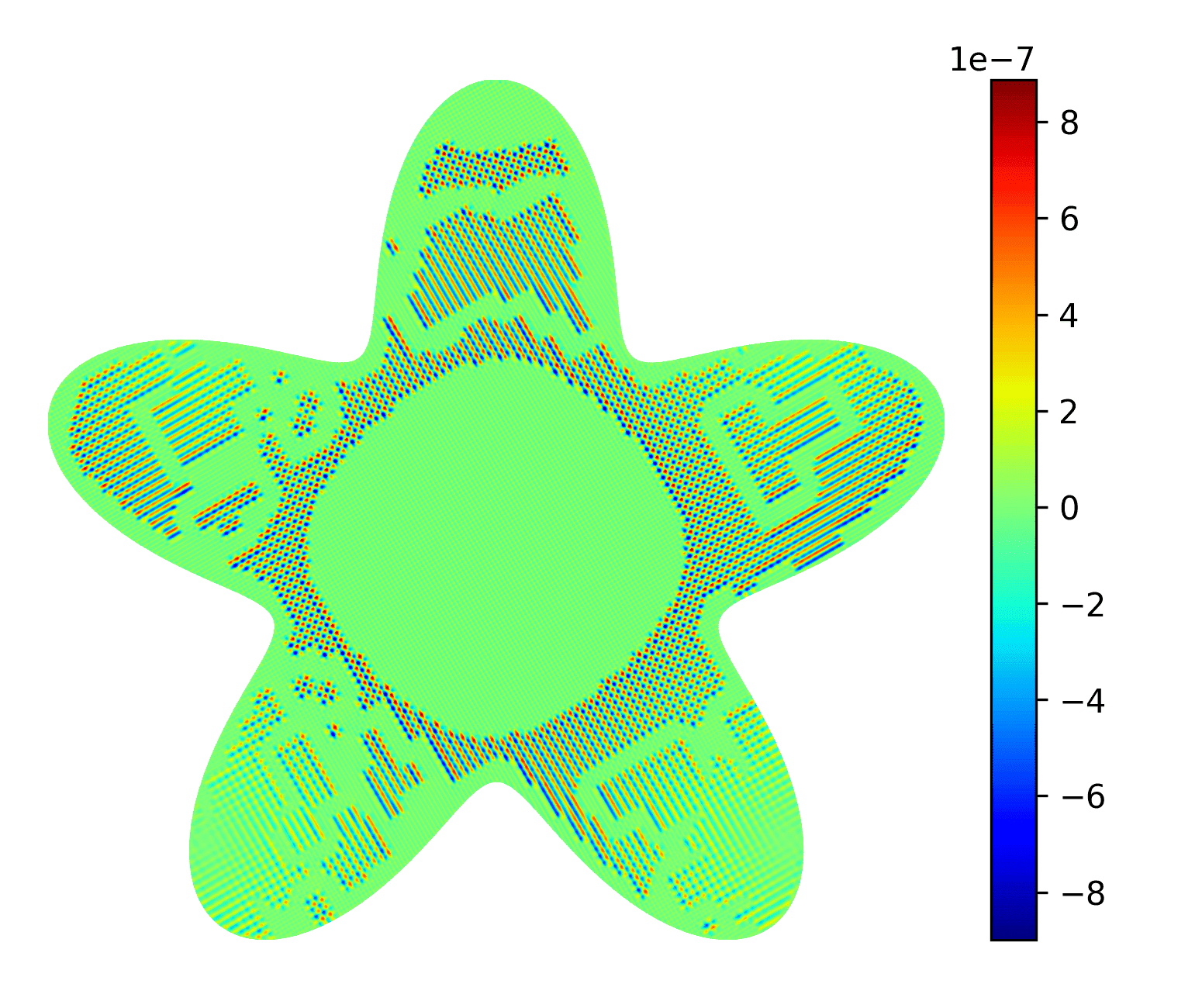}
		\caption{$\vartheta = 30$}
		\label{fig:star_theta30}
	\end{subfigure}\\[12pt]
		\begin{subfigure}[b]{0.45\textwidth}
		\centering
		\includegraphics[width=0.85\textwidth]{sensitivity_analysis_theta20.png}
		\caption{$\vartheta = 40 $}
		\label{fig:star_theta40}
	\end{subfigure}%
	\hspace{0.08\textwidth}
	\begin{subfigure}[b]{0.45\textwidth}
		\centering
		\includegraphics[width=0.85\textwidth]{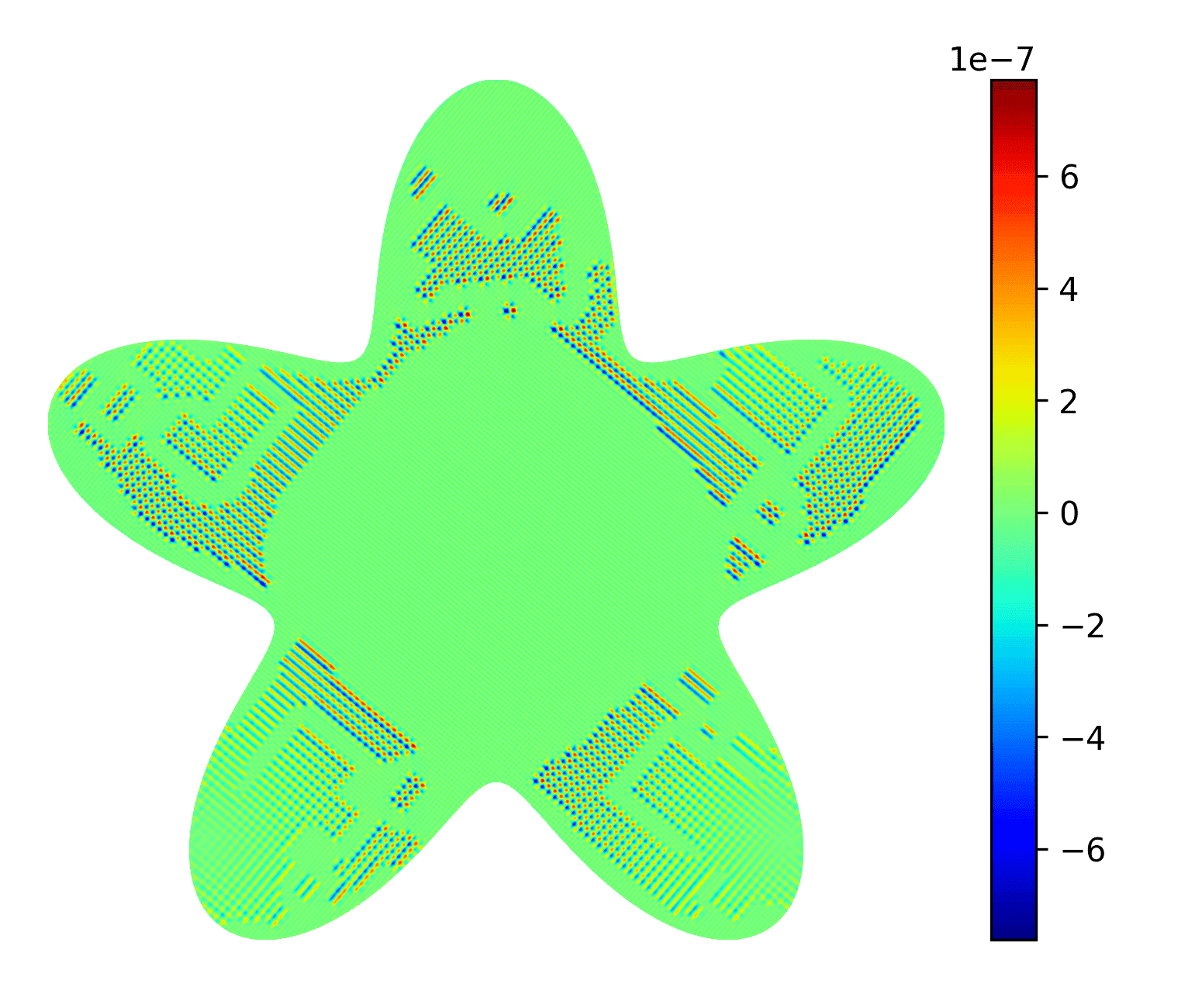}
		\caption{$\vartheta = 50$}
		\label{fig:star_theta50}
	\end{subfigure}
	\caption{(\subref{fig:star_domain}) Problem setup, and (\subref{fig:star_theta10})-(\subref{fig:star_theta50})  contour plots of the error, $u-u^h$, for the Laplace problem on the star shaped domain at the end of $6$ adaptive refinement steps for different angles of mesh rotation $\vartheta$.}
	\label{fig:star_conv}
\end{figure}

To study the sensitivity of the adaptive simulation framework to the cut-cell configurations, we consider the star-shaped domain shown in Figure~\ref{fig:star_domain} for various orientation angles $\vartheta$. The star-shaped domain is constructed using the level set function
\begin{equation*}
\psi (x,y) = R_{1} + R_{2} \sin(n_{\rm fold} \, {\rm arctan2}( y,x ))- \sqrt{ x^{2} + y^{2}},
\end{equation*}
with $R_{1}= 0.6$, $R_{2} = 0.2$ and $n_{\rm fold} = 5$ \cite{deprenter2020}. On the boundary of the domain, the Dirichlet data is set equal to the same exact solution \eqref{equation:exactlaplace} as in the previous example. For all orientations, an initial mesh of $10 \times 10$ elements covering the ambient domain $[-1,1]^2$ is considered, after which local refinements using second-order THB-splines are performed until the smallest elements have been refined six times.

\begin{figure}
	\centering
		\begin{subfigure}[b]{0.45\textwidth}
		\centering
		\includegraphics[width=\textwidth]{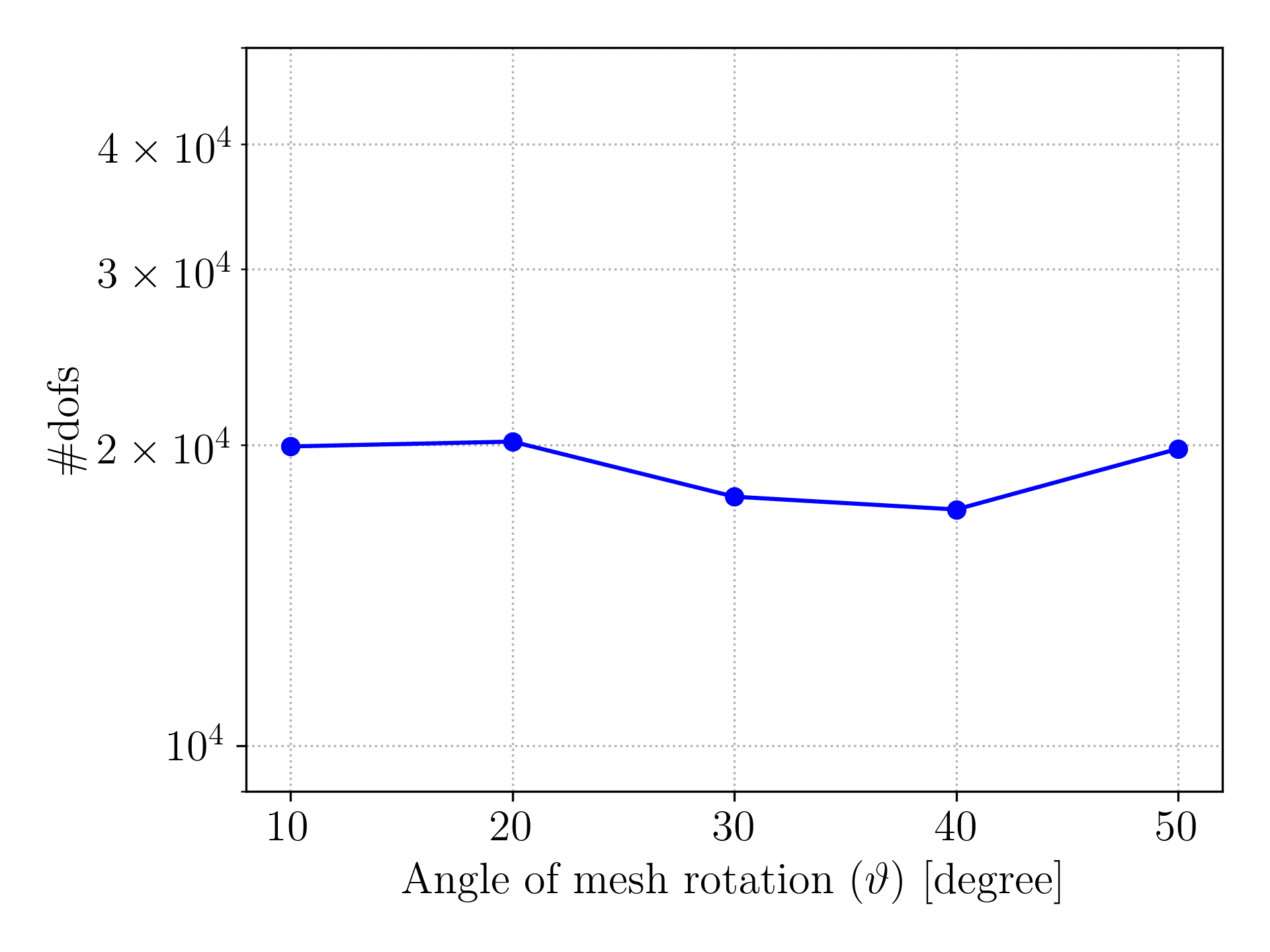}
		\caption{}\label{fig:rotationdofs}
		\end{subfigure}
	\hfill
		\begin{subfigure}[b]{0.45\textwidth}
		\centering
	    \includegraphics[width=\textwidth]{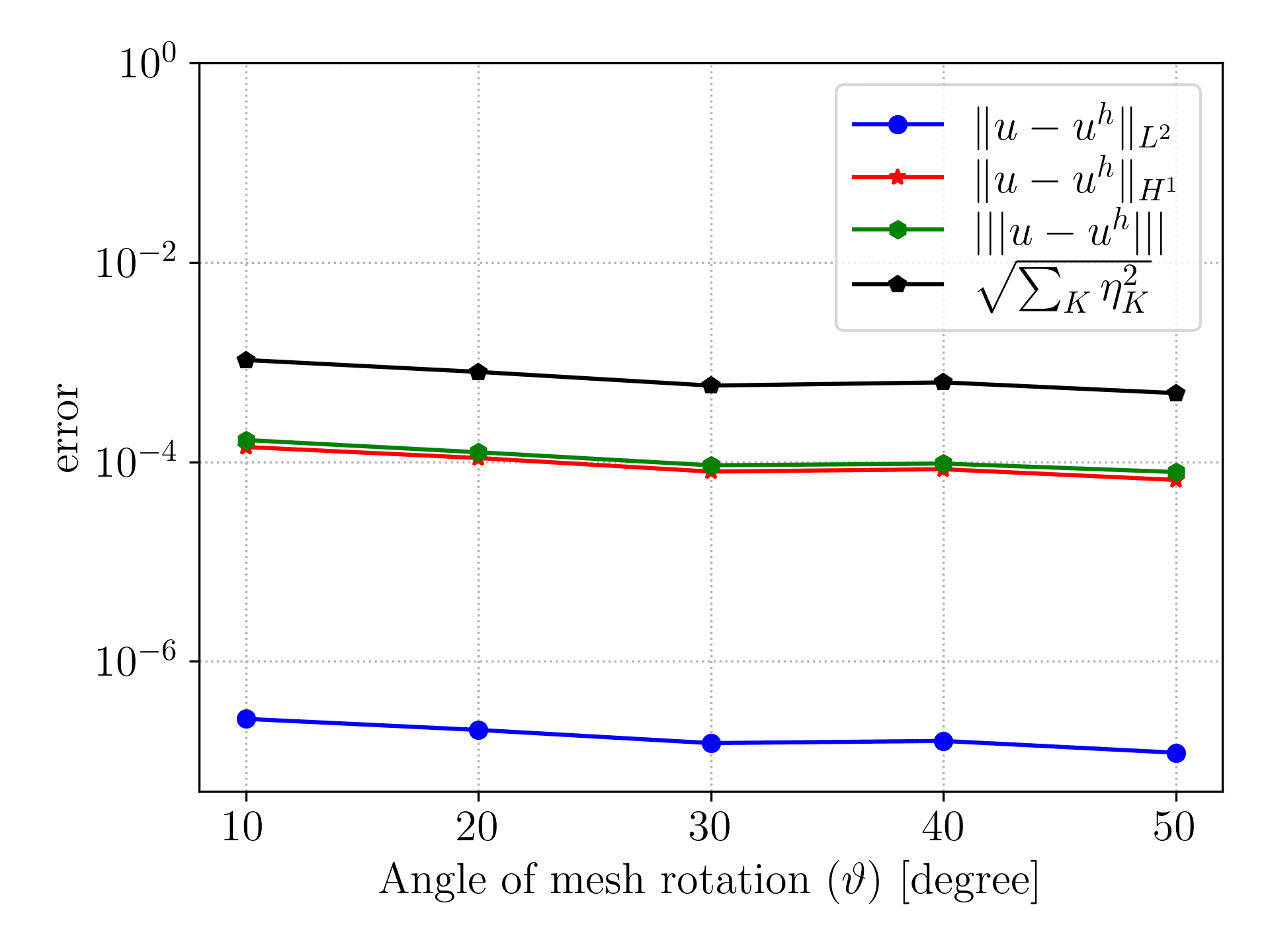}
	  \caption{}\label{fig:rotationerrors}
		\end{subfigure}
	\caption{(\subref{fig:rotationdofs}) Degrees of freedom, and (\subref{fig:rotationerrors}) error norms for the Laplace problem on the star shaped domain after 6 adaptive refinement steps for different angles of mesh rotation $\vartheta$.}
	\label{fig:sensitivity_analysis}
\end{figure}

Figures~\ref{fig:star_theta10}--\ref{fig:star_theta50} show the error  $u-u^h$ after completion of the refinement procedure. These figures convey that both the error and the refinement pattern are similar for all orientations. This is corroborated by the results in Figure~\ref{fig:sensitivity_analysis}, which indicates that both the number of degrees of freedom and the errors (in various norms) are insensitive to the orientation angle.

\subsubsection{Re-entrant corner}\label{sec:laplacelshape}
\begin{figure}
	\centering
	\begin{subfigure}[b]{0.45\textwidth}
		\centering
		\includegraphics[width=\textwidth]{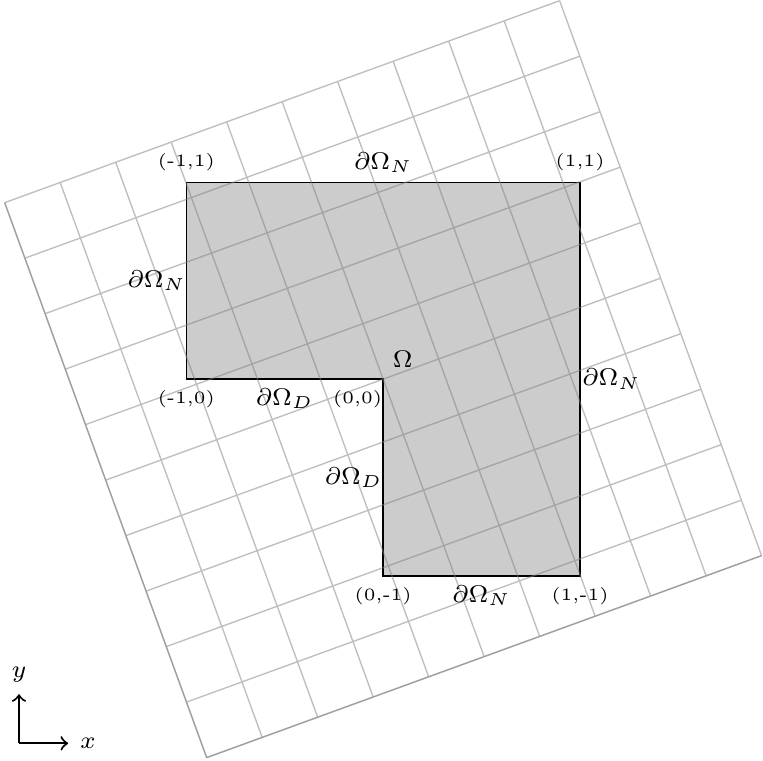}
		\caption{}
		\label{fig:Lshape}
	\end{subfigure}\hfill%
	\begin{subfigure}[b]{0.45\textwidth}
		\centering
		\includegraphics[width=\textwidth]{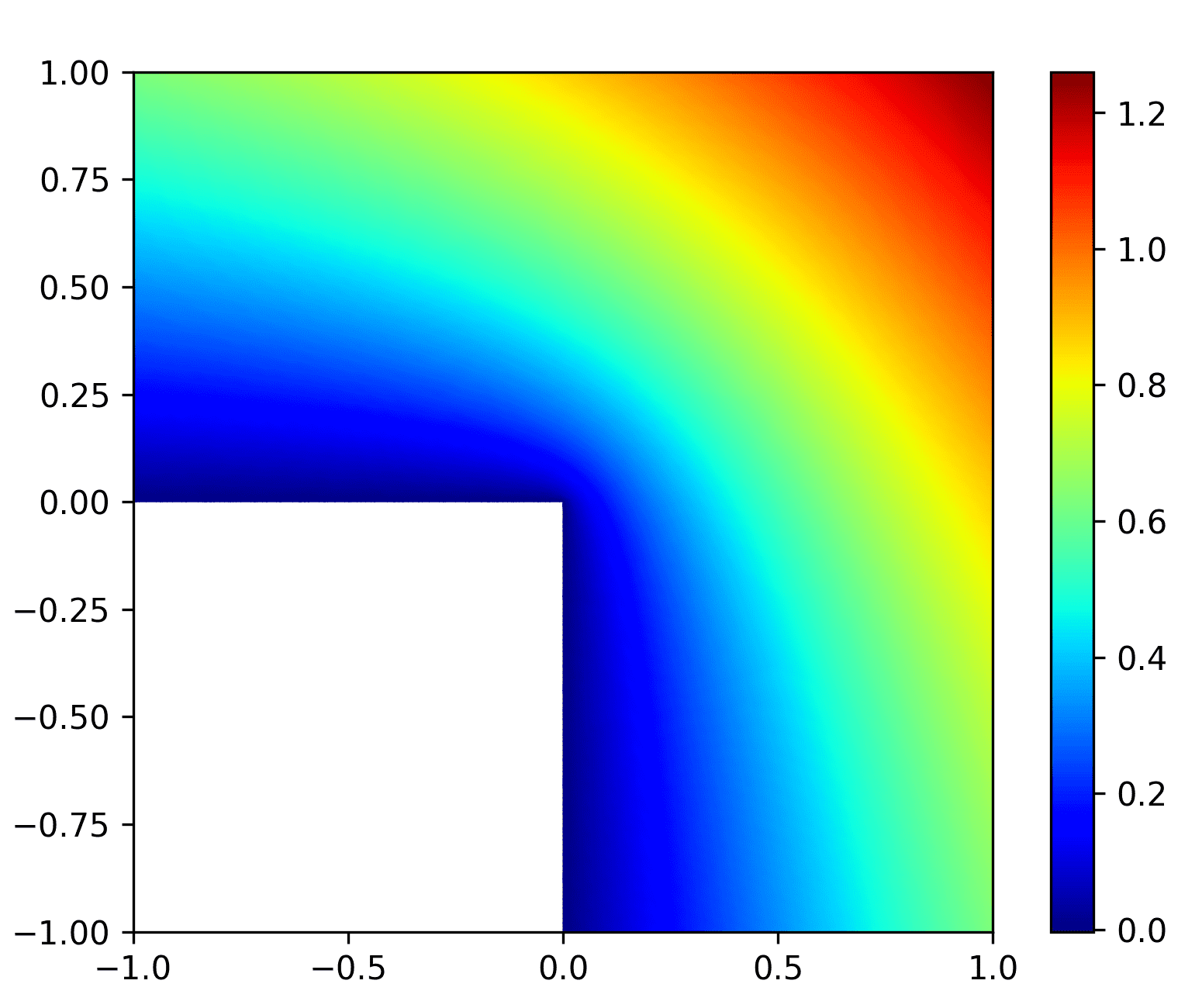}
		\caption{}
		\label{fig:Lshape_exact}
	\end{subfigure}
	\hspace{0.08\textwidth}
	\caption{(\subref{fig:Lshape}) Problem setup, and (\subref{fig:Lshape_exact}) the exact solution $u(x,y)$, Eq.~\eqref{equation:lshape_exactlaplace}, for the Laplace problem on the re-entrant corner domain.}
	\label{fig:Laplace}
\end{figure}

To study the behavior of the adaptive simulation strategy for problems with (weakly) singular solutions, we consider a domain with a re-entrant corner, as shown in Figure~\ref{fig:Lshape}. The data on the Dirichlet and Neumann boundaries, $u |_{\partial \Omega_{D}}=g=0$ and $\partial_n u |_{\partial \Omega_{N}} = q$, is set to match the exact solution \cite{kuru2014,dangella2016}
\begin{equation}\label{equation:lshape_exactlaplace}
u(x,y) = (x^2 + y^2)^\frac{1}{3} \cos \left( \frac{2}{3} {\rm arctan2}(x - y, x + y) \right).
\end{equation}

\begin{figure}
	\centering
	\begin{subfigure}[b]{0.45\textwidth}
		\centering
		\includegraphics[width=\textwidth]{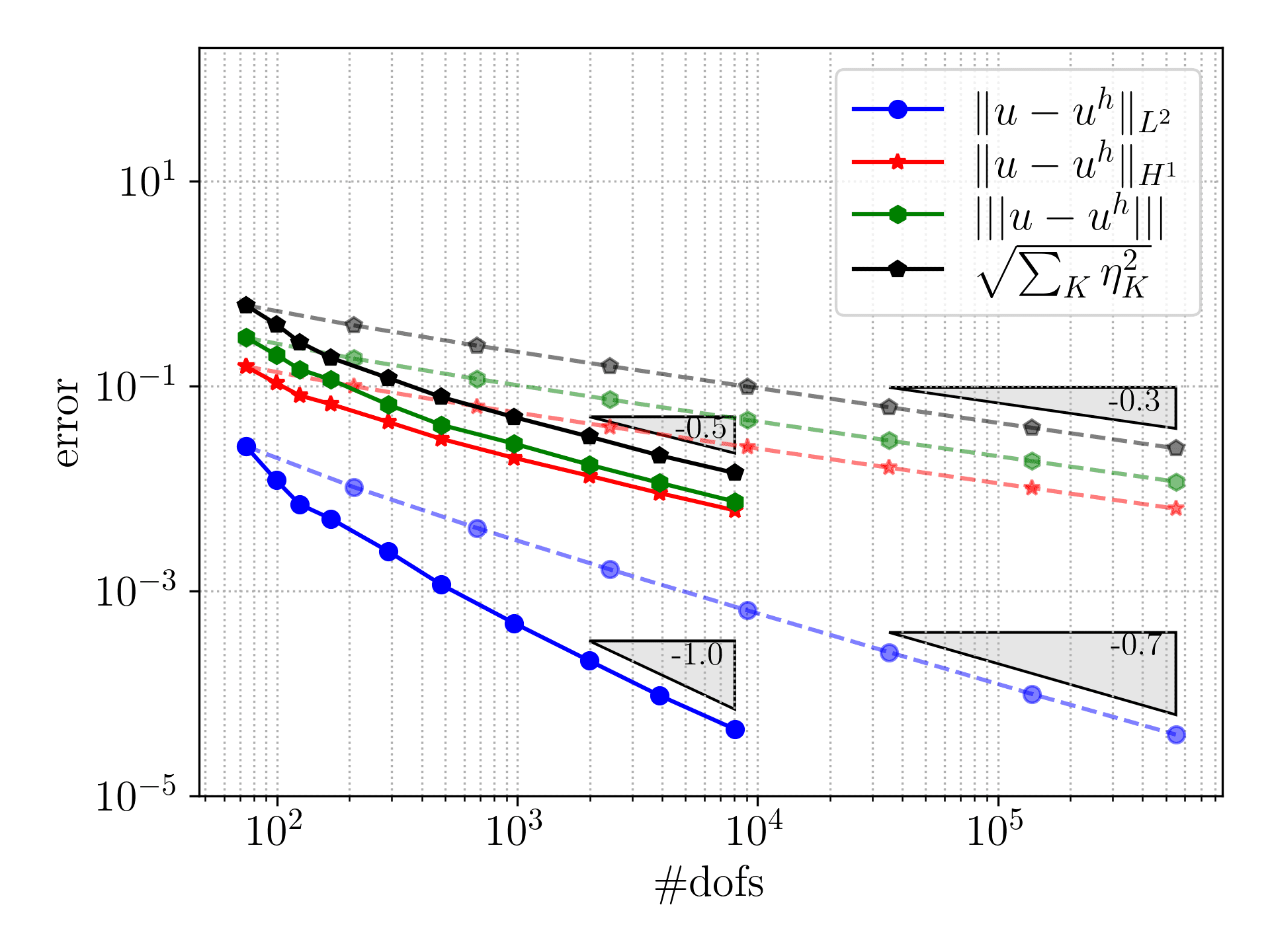}
		\caption{$\order = 1$}
		\label{fig:Lshape_p1}
	\end{subfigure}\hspace{0.08\textwidth}%
	\begin{subfigure}[b]{0.45\textwidth}
		\centering
		\includegraphics[width=\textwidth]{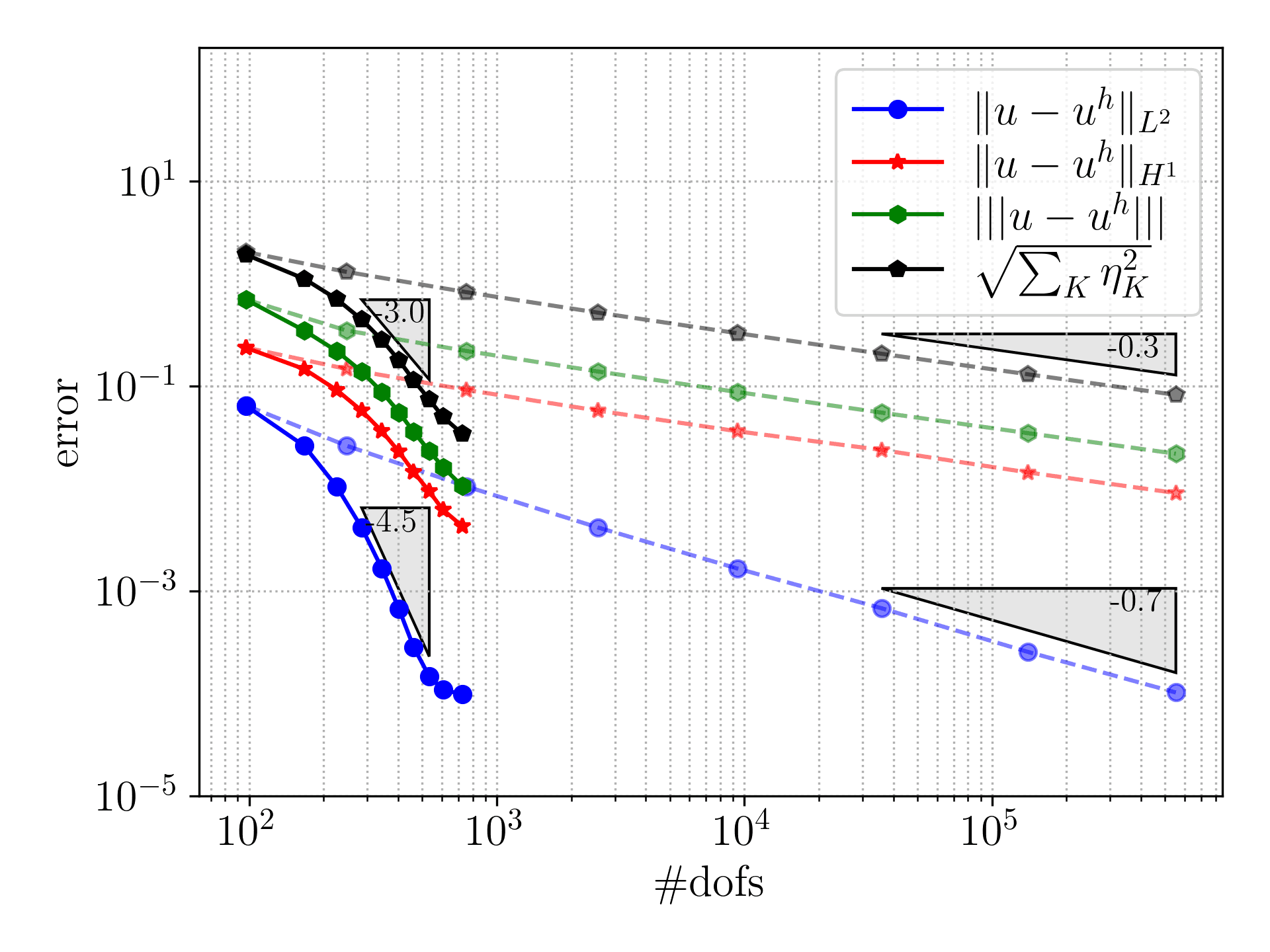}
		\caption{$\order = 2$}
		\label{fig:Lshape_p2}
	\end{subfigure}
	\hspace{0.08\textwidth}
		\caption{Error convergence results for the Laplace problem on the re-entrant corner domain under residual-based adaptive refinement (solid) and uniform refinement (dashed) for linear ($\order =1$) and quadratic ($\order =2$) basis functions.}
	\label{fig:Lshape_conv}
\end{figure}

The convergence behavior of the $L^{2}$-error, $H^{1}$-error, energy norm error \eqref{eq:energynorm} and the residual-based estimator \eqref{equation:laplaceindicators} is studied for uniform refinement and residual-based adaptive refinement. Both refinement procedures start from an initial mesh of $10 \times 10$ elements formed on the ambient domain $[-\frac{3}{2},\frac{3}{2}]^2$. The convergence results for first and second order B-splines are shown  in Figure~\ref{fig:Lshape_p1} and Figure~\ref{fig:Lshape_p2}, respectively.

Under uniform refinement, the convergence rates are impeded by the weak singularity at the re-entrant corner. For the $L^2$-error and $H^1$-error, suboptimal rates of $\mathcal{O}(n^{-\frac{2}{3}})$ and $\mathcal{O}(n^{-\frac{1}{3}})$ are observed, which is in agreement with the expected rates \cite{babuvska1996}. These rates are independent of the order of the approximation, as the regularity of the exact solution limits the rate already for the linear case. As for the cases considered above, the energy error and estimator follow the convergence of the $H^1$-error.

Using the adaptive refinement strategy with linear basis functions, the optimal rates of $\mathcal{O}(n^{-1})$ and $\mathcal{O}(n^{-\frac{1}{2}})$ are recovered for the $L^2$-error and $H^1$-error, respectively. For the quadratic case, rates that are substantially higher than the theoretical rates are observed. We attribute this to pre-asymptotic behavior, in which the refinement pattern as shown in Figure~\ref{fig:Lshape_conv} is strongly focused on the re-entrant corner singularity. After the first two steps, the errors become dominated by the singularity at the re-entrant corner, which results in the further refinement of the few elements in the vicinity of the corner. These refinements do reduce the error, while they only introduce a limited number of additional degrees of freedom. The observed flattening in the rate of the $L^2$-error in the quadratic case is caused by the refinement reaching the maximum level in the elements in the corner, which causes the marking strategy to tag elements that do not carry the largest error contributions.

\begin{figure}
	\centering
	\begin{subfigure}[b]{0.45\textwidth}
		\centering
		\includegraphics[width=0.8\textwidth]{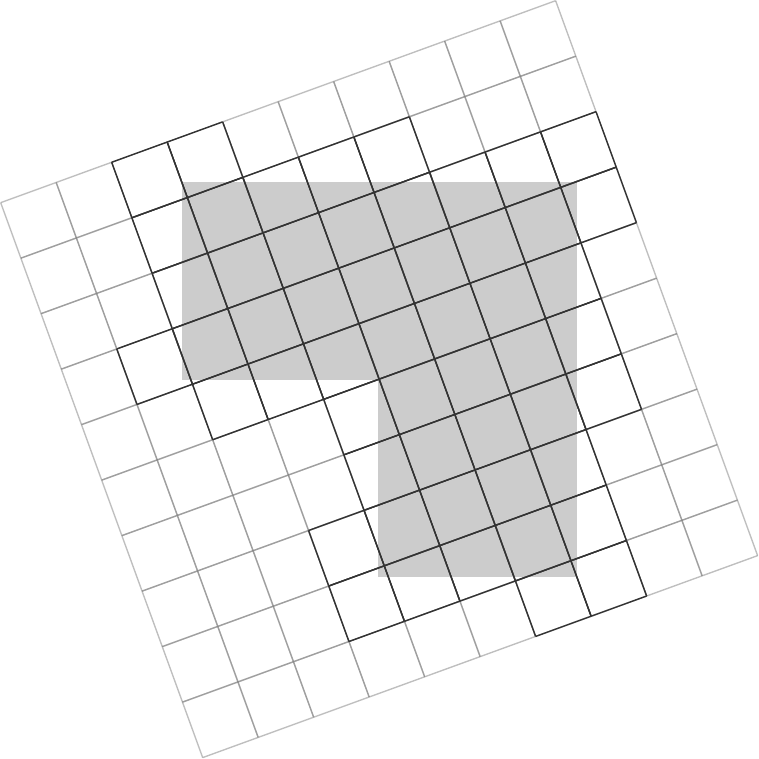}
		\caption{Initial mesh}
	\end{subfigure}	\hspace{0.08\textwidth}%
	\begin{subfigure}[b]{0.45\textwidth}
		\centering
		\includegraphics[width=0.8\textwidth]{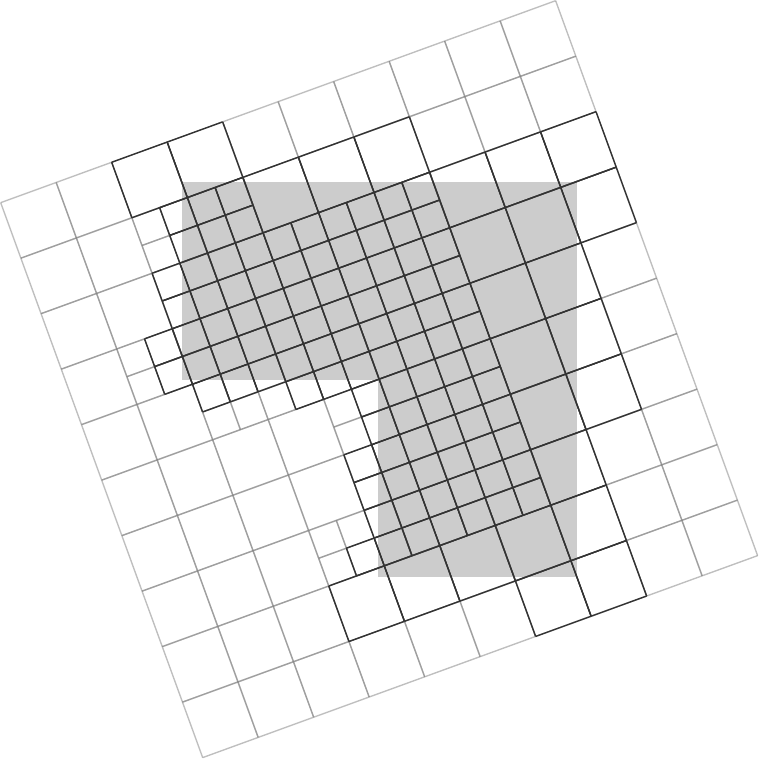}
		\caption{Step $1$}
	\end{subfigure}\\[12pt]
	\begin{subfigure}[b]{0.45\textwidth}
		\centering
		\includegraphics[width=0.8\textwidth]{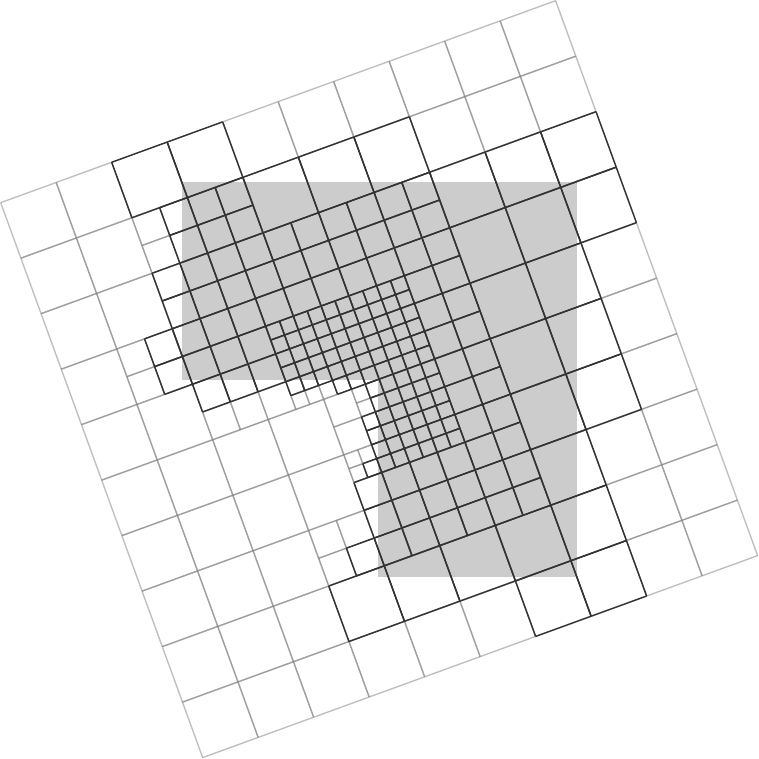}
		\caption{Step $2$}
	\end{subfigure}	\hspace{0.08\textwidth}%
	\begin{subfigure}[b]{0.45\textwidth}
		\centering
		\includegraphics[width=0.8\textwidth]{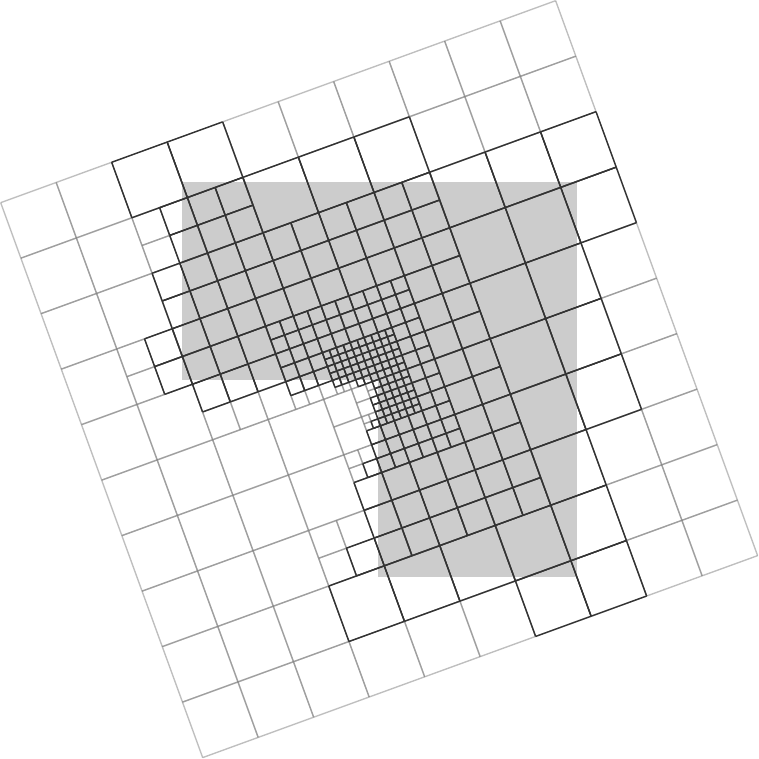}
		\caption{Step $3$}
	\end{subfigure}\\[12pt]
	\begin{subfigure}[b]{0.45\textwidth}
		\centering
		\includegraphics[width=0.8\textwidth]{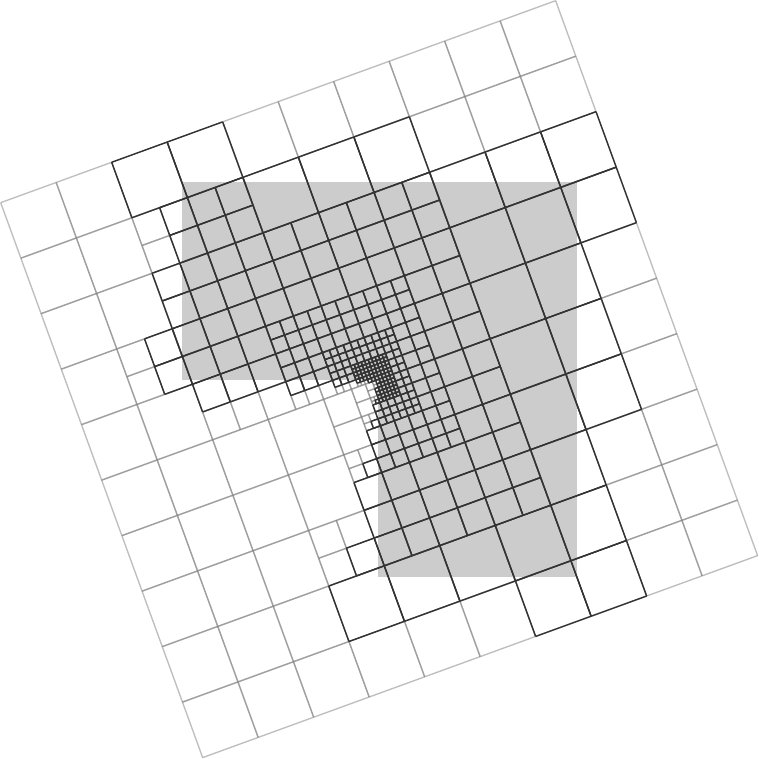}
		\caption{Step $4$}
	\end{subfigure}	\hspace{0.08\textwidth}%
	\begin{subfigure}[b]{0.45\textwidth}
		\centering
		\includegraphics[width=0.8\textwidth]{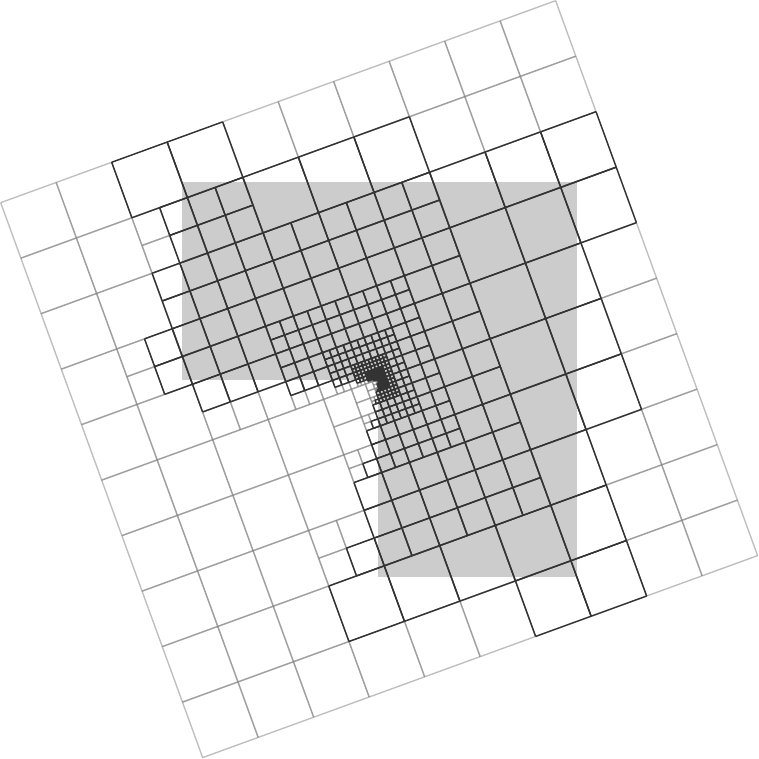}
		\caption{Step $5$}
	\end{subfigure}
	\caption{Evolution of the mesh using the adaptive refinement procedure for the Laplace problem on the re-entrant corner domain  using $\order = 2$.}
	\label{fig:Lshape_adap_mesh}
\end{figure}


\subsection{Steady viscous flow} \label{sec:stokesexamples}
We regard the two-dimensional Stokes flow problem on a quarter annulus ring domain with a smooth solution and on the above-introduced re-entrant corner domain with a singular solution. We consider equal-order discretizations for the velocity and pressure fields using optimal regularity (TH)B-splines of degree $\order=1$ and $\order=2$. For the Nitsche and ghost-penalty parameter the same settings are used as for the Laplace problem considered above, \emph{i.e.}, $\beta=50$ and $\gamma_g=10^{-(\order+2)}$. In addition, a skeleton-penalty parameter of $\gamma_s= 10^{-(\order+1)}$ is used for all simulations.

\subsubsection{Quarter annulus ring}\label{sec:annulus}
\begin{figure}
	\centering
	\begin{subfigure}[b]{0.45\textwidth}
		\centering
		\includegraphics[width=0.8\textwidth]{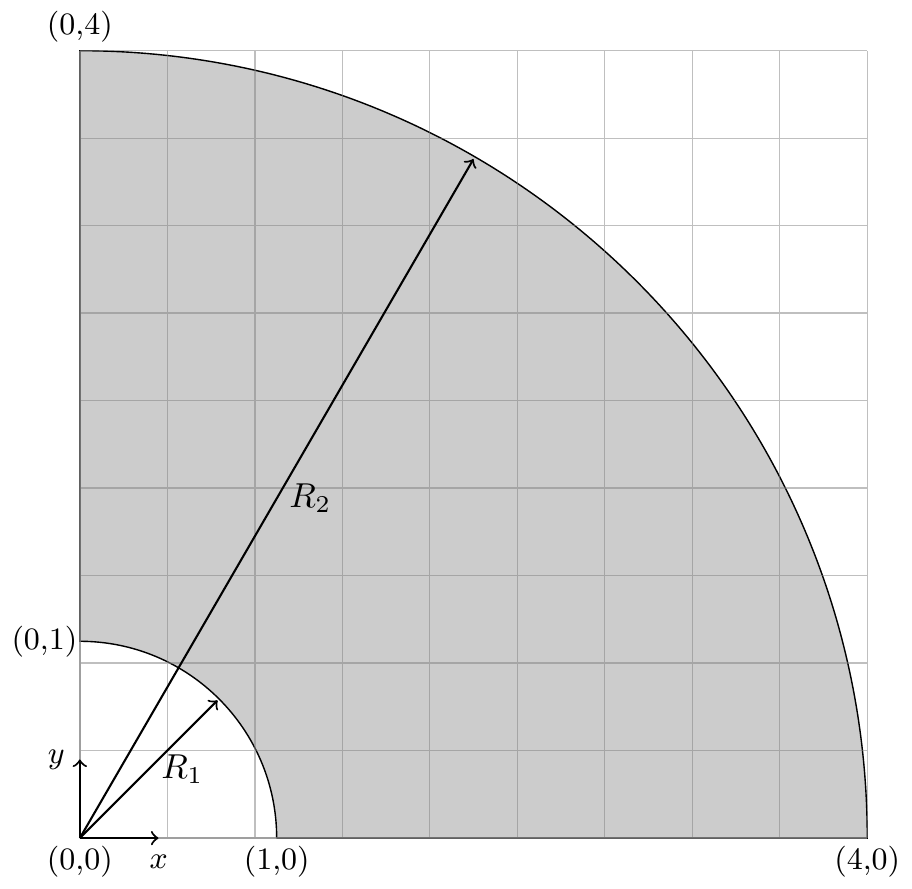}
		\caption{}
		\label{fig:annulus_domain}
	\end{subfigure}\hspace{0.08\textwidth}%
	\begin{subfigure}[b]{0.45\textwidth}
		\centering
		\includegraphics[width=\textwidth]{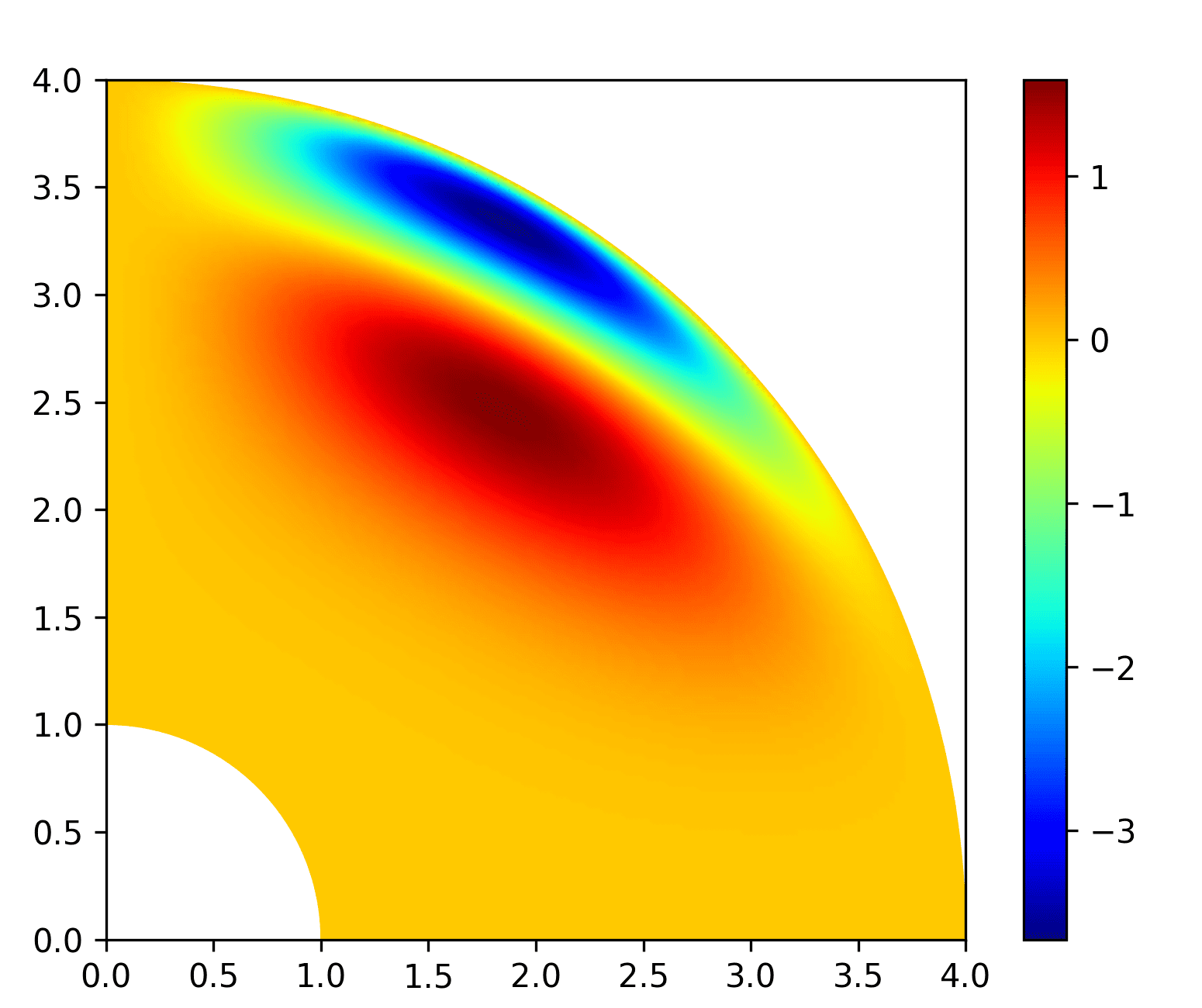}
		\caption{$u_1 (x,y)$}
		\label{fig:annulus_exact_u1}
	\end{subfigure}\\[12pt]
	\begin{subfigure}[b]{0.45\textwidth}
		\centering
		\includegraphics[width=0.96\textwidth]{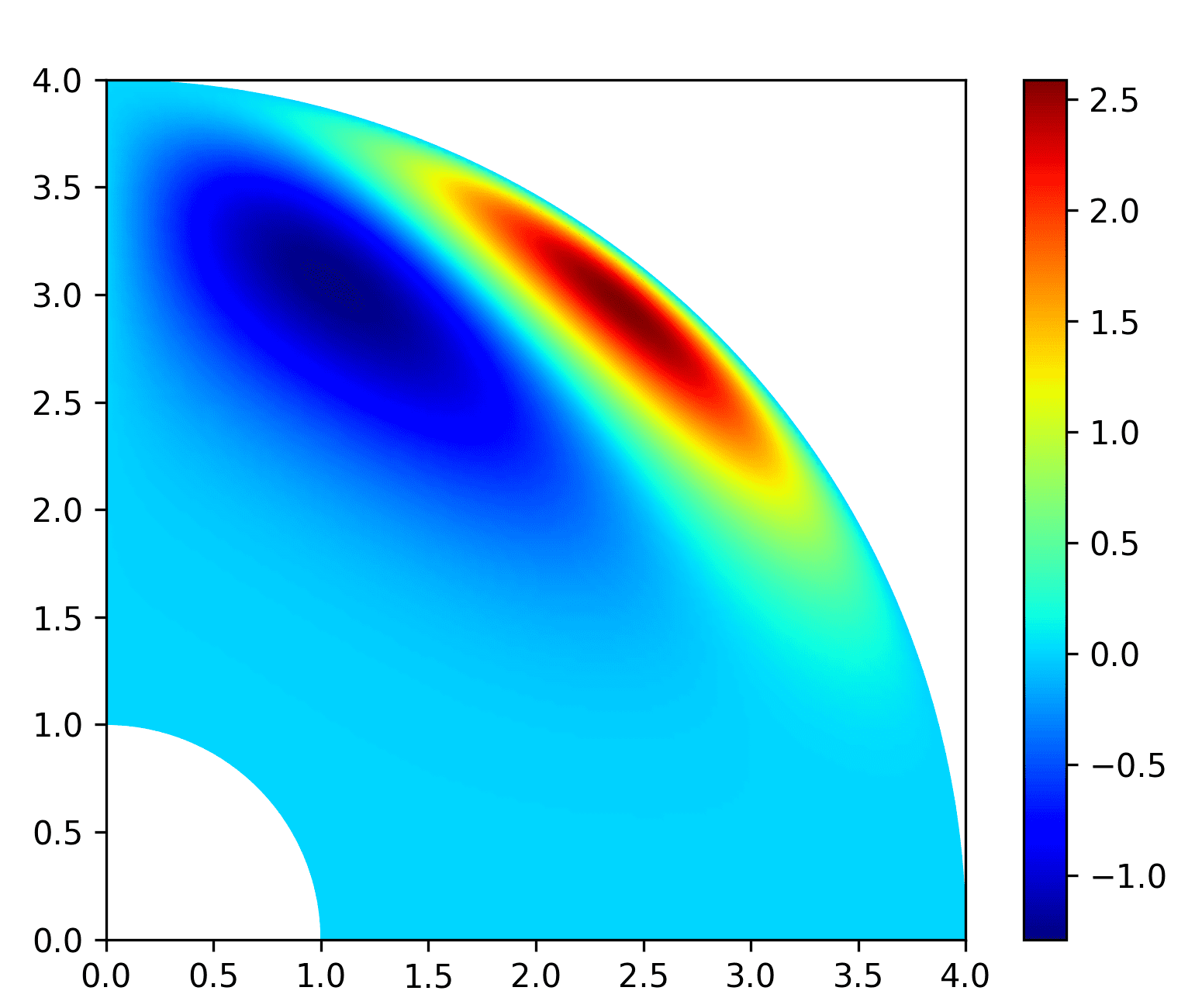}
		\caption{$u_2 (x,y)$}
		\label{fig:annulus_exact_u2}
	\end{subfigure}\hspace{0.08\textwidth}%
	\begin{subfigure}[b]{0.45\textwidth}
		\centering
		\includegraphics[width=\textwidth]{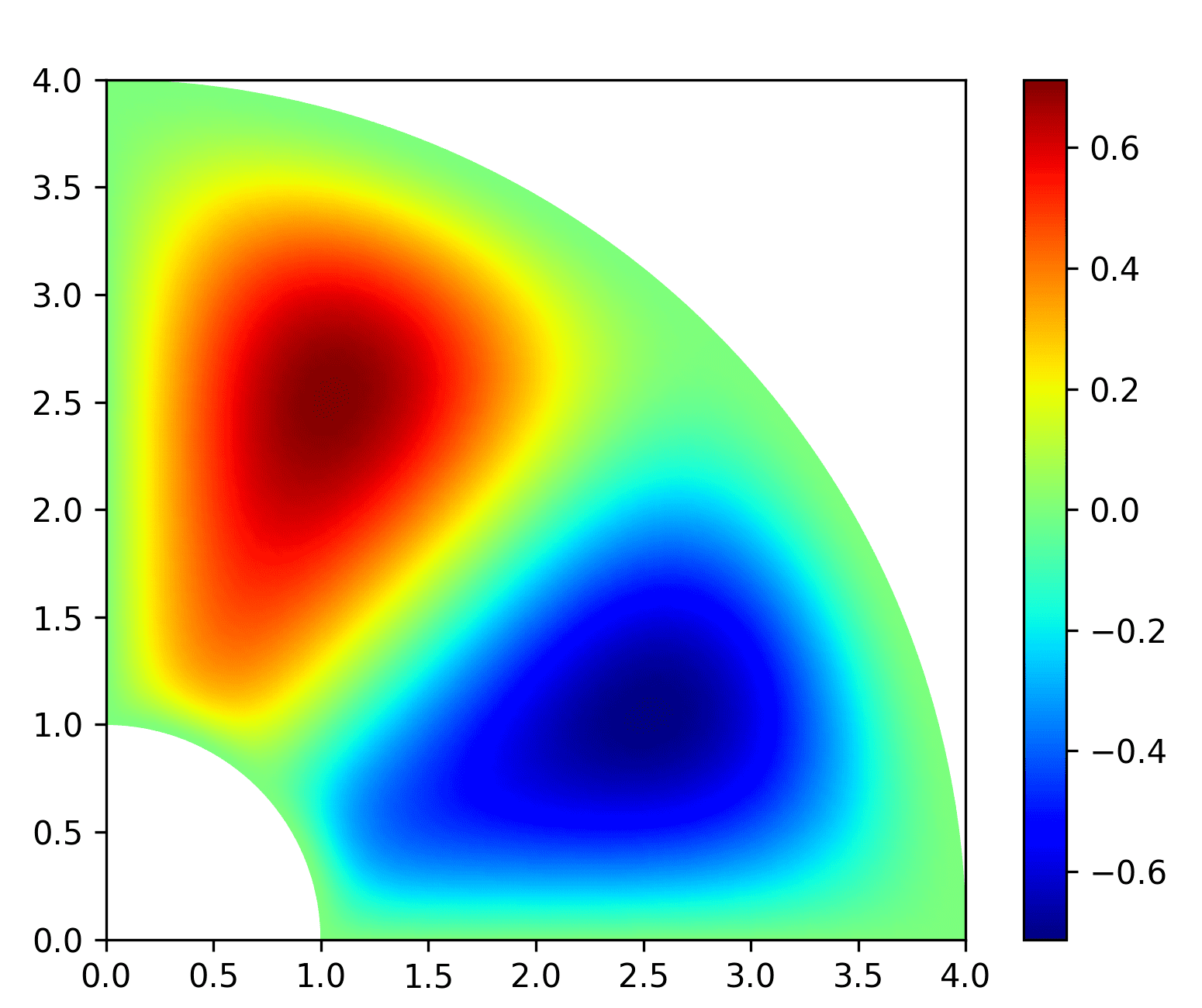}
		\caption{$p (x,y)$}
		\label{fig:annulus_exact_p}
	\end{subfigure}
		\caption{(\subref{fig:annulus_domain}) Problem setup, and (\subref{fig:annulus_exact_u1})-(\subref{fig:annulus_exact_p}) the exact solution components $u_1(x,y)$, $u_2(x,y)$, and $p(x,y)$, defined in Eq.~\eqref{eq:annulusexact}, for the Stokes problem on the quarter annulus ring domain.}
	\label{fig:annulus}
\end{figure}

We consider an annulus ring domain $\Omega = \{ (x,y) \in \mathbb{R}_{>0}^{2} : R_{1}^{2} < x^2 + y^2 < R_{2}^2\}$ with inner radius $R_{1} = 1$, outer radius $R_{2} = 4$, Dirichlet boundary $\partial \Omega_{D}$ and Neumann boundary $\partial \Omega_{N}$, as shown in  Figure~\ref{fig:annulus_domain}. The Dirichlet data $\boldsymbol{g}$ and Neumann data $\boldsymbol{t}$ are prescribed in accordance with the divergence-free manufactured solution \cite{hoang2017}
\begin{equation}
\begin{aligned}
	u_1 (x, y) &= 10^{-6} x^2 y^4 (x^2 + y^2 -1) (x^2 + y^2 - 16) \left(5 x^4 + 18 x^2 y^2 - 85 x^2 + 13 y^4 - 153 y^2 + 80 \right), \\
	u_2 (x, y) &= 10^{-6} x y^5 (x^2 + y^2 - 1) (x^2 + y^2 - 16) \left(102 x^2 + 34 y^2 - 10 x^4 - 12 x^2 y^2 - 2y^4 - 32 \right), \\
	p (x,y)&= 10^{-7} x y (y^2 - x^2) (x^2 + y^2 - 16)^2 (x^2 + y^2 - 1)^2 \exp( 14 (x^2 + y^2)^{- 1/2}).
\end{aligned}
\label{eq:annulusexact}
\end{equation}
The body force $\boldsymbol{f}$ in the Stokes problem \eqref{equation:stokesequations} is determined based on this manufactured solution, with the viscosity set to $\mu = 1$.

%
%

\begin{figure}
	\centering
	\begin{subfigure}[b]{0.45\textwidth}
		\centering
		\includegraphics[width=\textwidth]{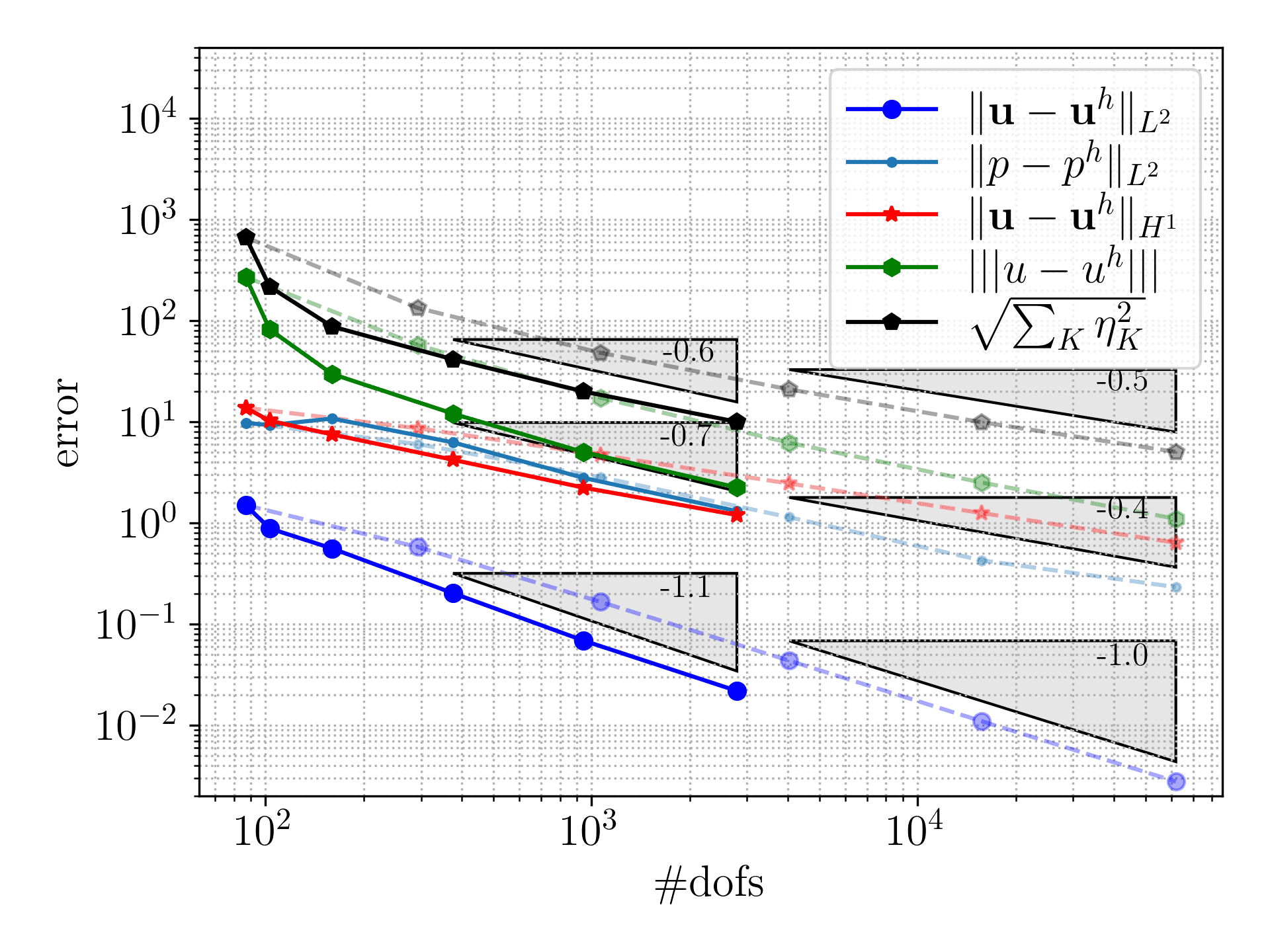}
		\caption{$\order = 1$}
		\label{fig:annulus_p1}
	\end{subfigure}\hspace{0.08\textwidth}%
	\begin{subfigure}[b]{0.45\textwidth}
		\centering
		\includegraphics[width=\textwidth]{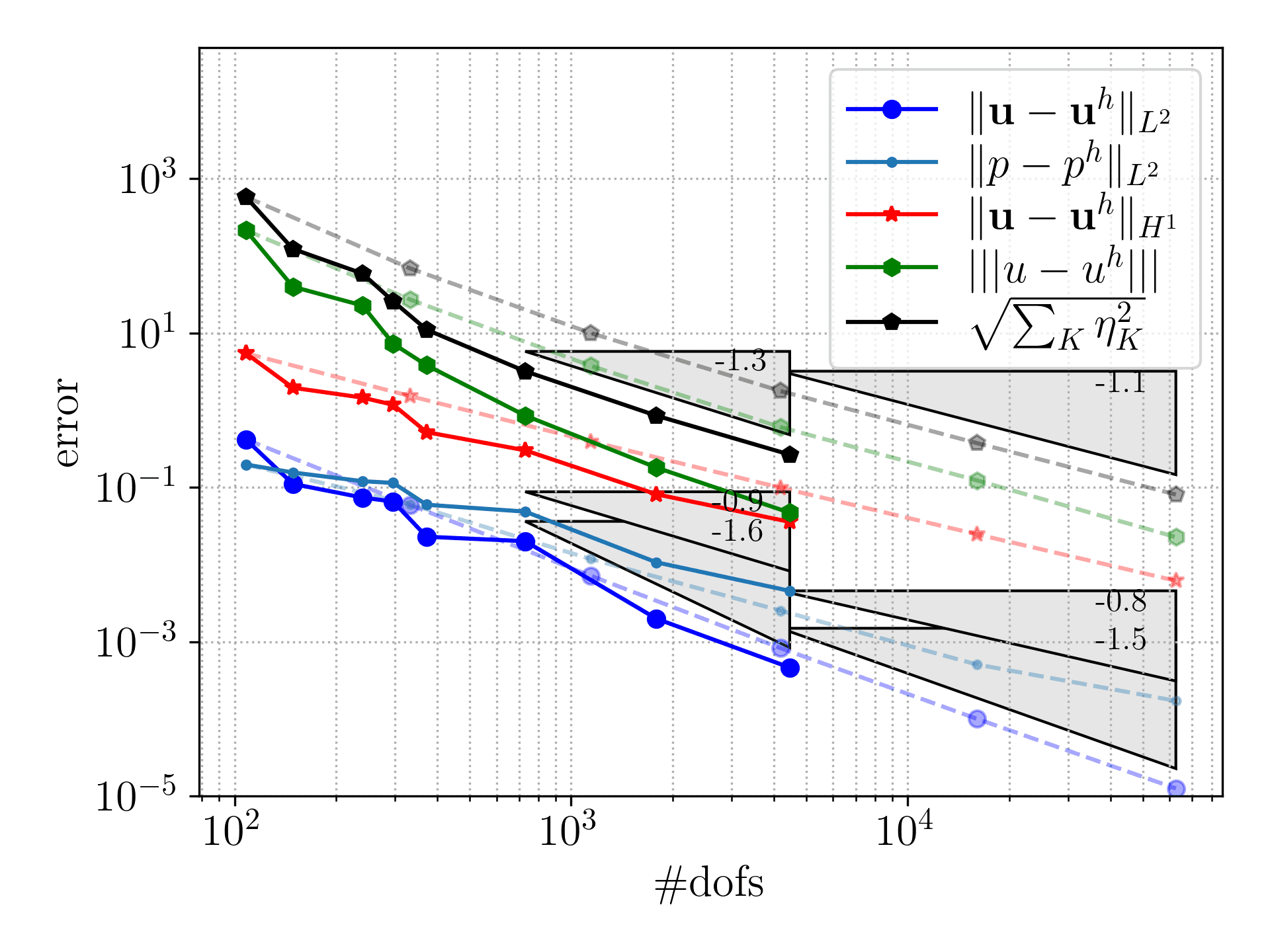}
		\caption{$\order = 2$}
		\label{fig:annulus_p2}
	\end{subfigure}
	\hspace{0.08\textwidth}
		\caption{Error convergence results for the Stokes  problem on the quarter annulus ring domain under residual-based adaptive refinement (solid) and uniform refinement (dashed) for linear ($\order =1$) and quadratic ($\order =2$) basis functions.}
	\label{fig:annulus_conv}
\end{figure}

Figure~\ref{fig:annulus_conv} displays the convergence results for the annulus ring problem. Both the uniform refinement results and the adaptive refinement results are obtained starting from a $9 \times 9$ uniform mesh on the ambient domain $[0,R_2]^2=[0,4]^2$. A good resemblance with the optimal rates of $\mathcal{O}(n^{-\frac{1}{2}k})$ in the velocity $H^1$-norm and pressure $L^2$-norm is observed, and, as expected, the rate of the velocity $L^2$-error is $\mathcal{O}(n^{-\frac{1}{2}(k+1)})$. The error in the energy norm \eqref{eq:energynormstokes} is observed to converge with the same rate as the $H^1$-norm velocity error and $L^2$-norm pressure error, which is in agreement with the definition of the energy norm. As expected, the error estimator bounds the error in the energy norm from above.

\begin{figure}
	\centering
	\begin{subfigure}[b]{0.45\textwidth}
		\centering
		\includegraphics[width=0.8\textwidth]{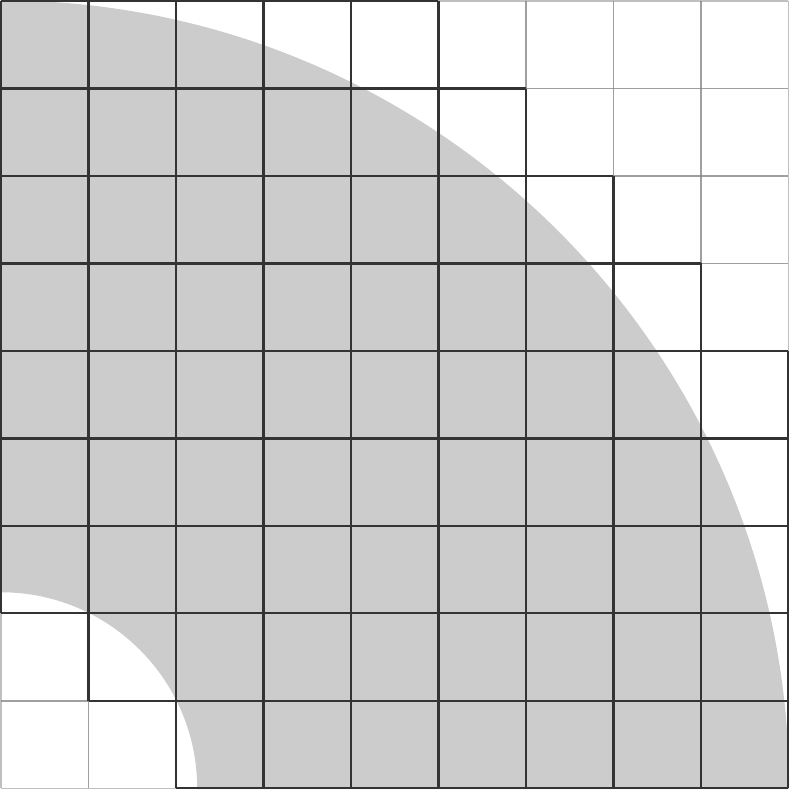}
		\caption{Initial mesh}
	\end{subfigure}%
	\hspace{0.08\textwidth}
	\begin{subfigure}[b]{0.45\textwidth}
		\centering
		\includegraphics[width=0.8\textwidth]{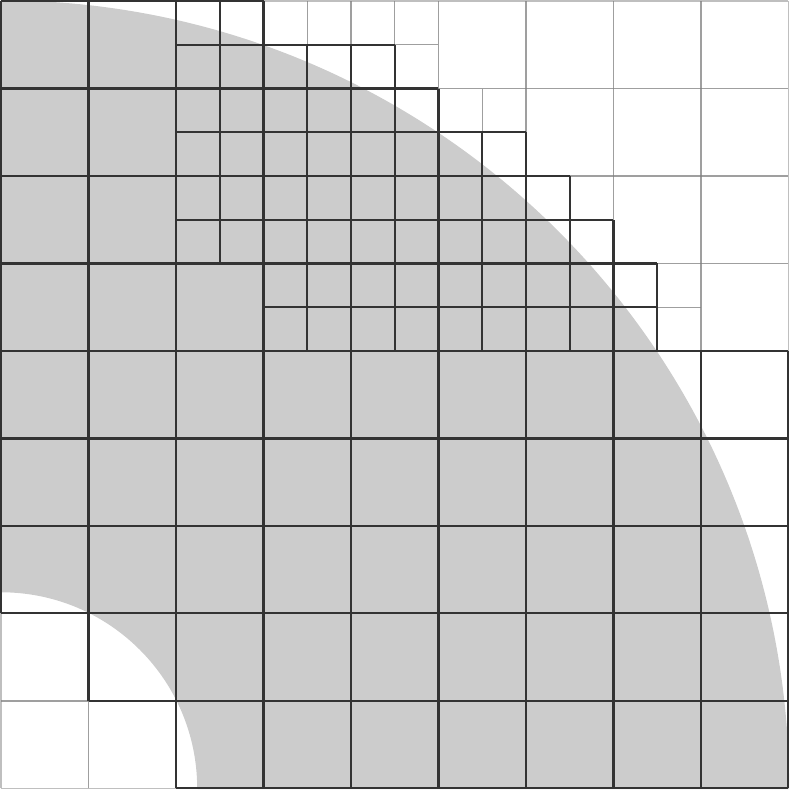}
		\caption{Step $1$}
	\end{subfigure}\\[12pt]
	\begin{subfigure}[b]{0.45\textwidth}
		\centering
		\includegraphics[width=0.8\textwidth]{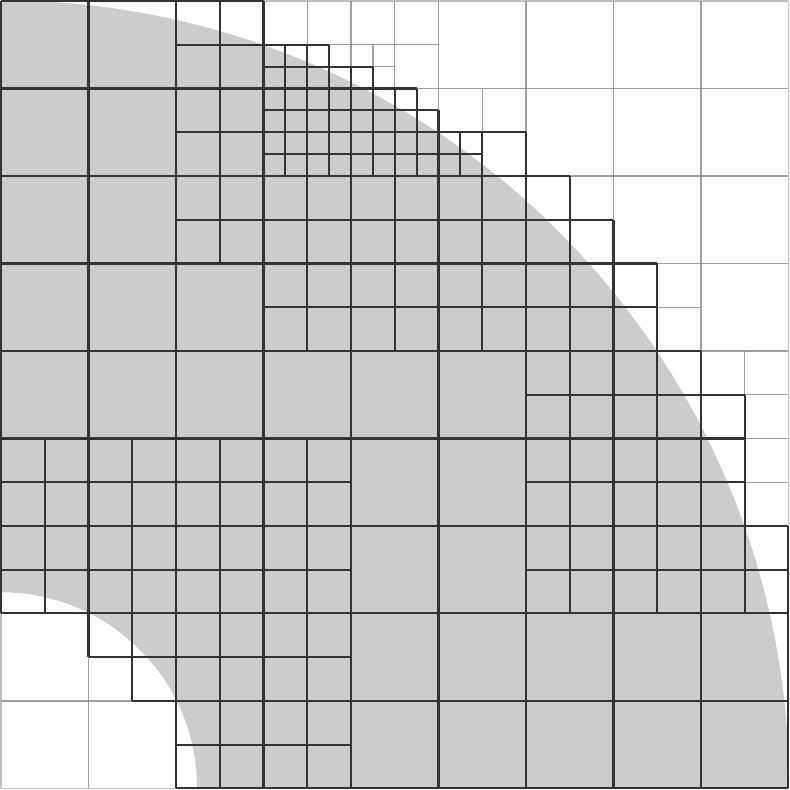}
		\caption{Step $2$}
	\end{subfigure}%
	\hspace{0.08\textwidth}
	\begin{subfigure}[b]{0.45\textwidth}
		\centering
		\includegraphics[width=0.8\textwidth]{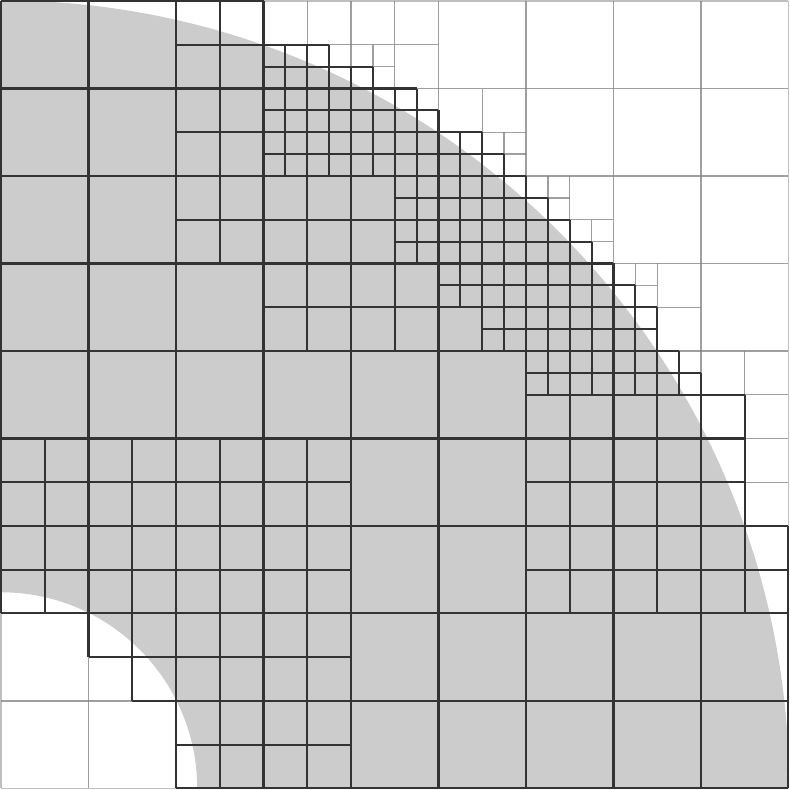}
		\caption{Step $3$}
	\end{subfigure}\\[12pt]
		\begin{subfigure}[b]{0.45\textwidth}
		\centering
		\includegraphics[width=0.8\textwidth]{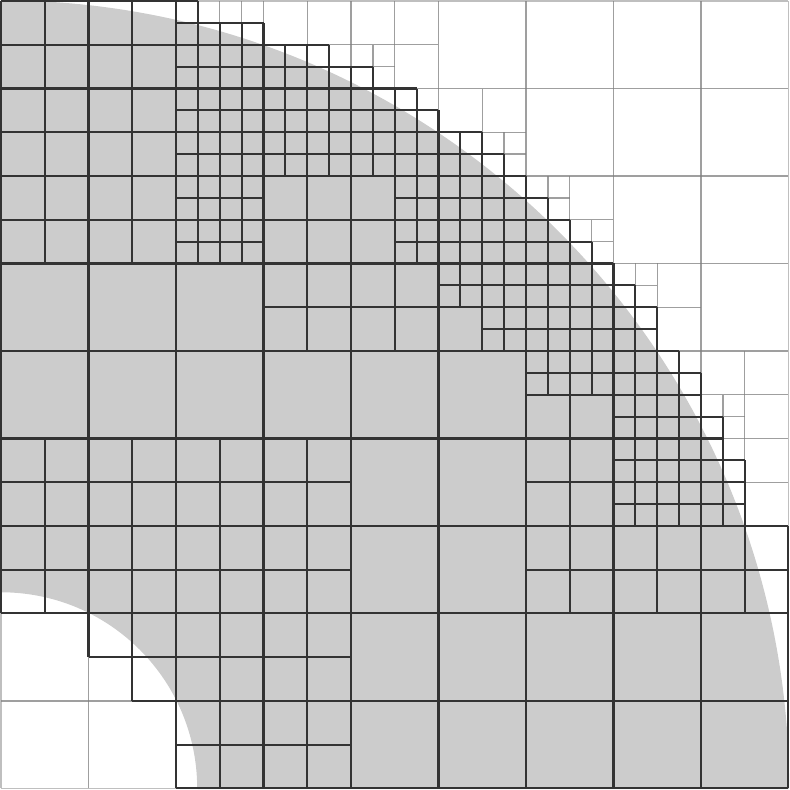}
		\caption{Step $4$}
	\end{subfigure}%
	\hspace{0.08\textwidth}
	\begin{subfigure}[b]{0.45\textwidth}
		\centering
		\includegraphics[width=0.8\textwidth]{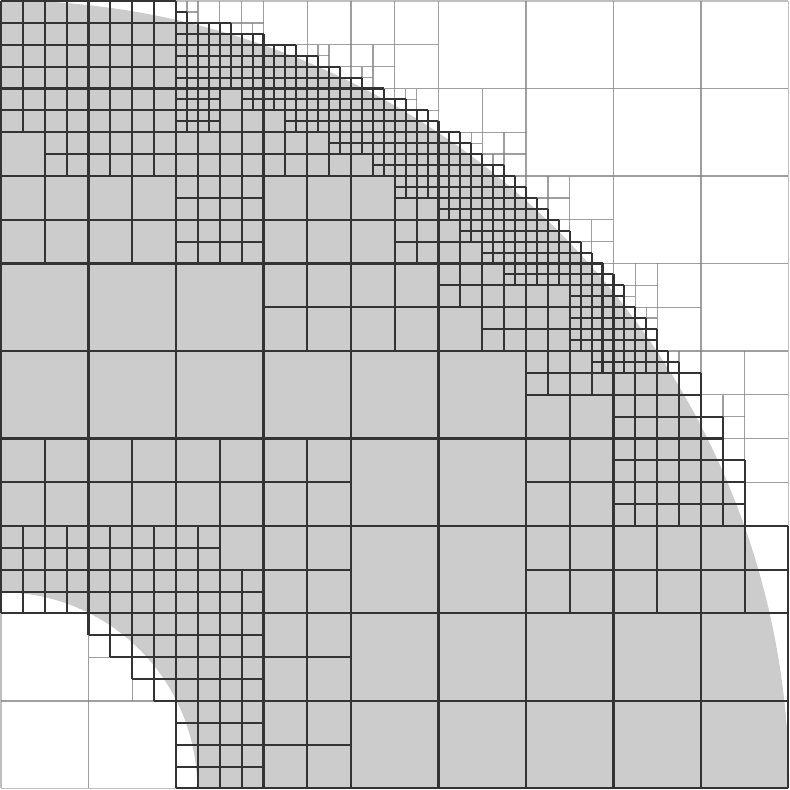}
		\caption{Step $5$}
	\end{subfigure}
	\caption{Evolution of the mesh using the adaptive refinement procedure for the Stokes problem on the quarter annulus ring domain  using $\order = 2$.}
	\label{fig:annulus_adap_mesh}
\end{figure}

Although optimal convergence rates are obtained using uniform refinements, the adaptive refinement procedure is observed to substantially improve the error for a fixed number of degrees of freedom. This behavior is explained by the observed refinement patterns, as shown in Figure~\ref{fig:annulus_adap_mesh}. Although the exact solution \eqref{eq:annulusexact} is smooth, in particular the steep gradients in the velocity solution lead to local refinements. This effectively reduces the error when compared to a uniform refinement with a similar number of degrees of freedom.

\subsubsection{Re-entrant corner}\label{sec:stokeslshape}
\begin{figure}
	\centering
	\begin{subfigure}[H]{0.45\textwidth}
		\centering
		\includegraphics[width=1\textwidth]{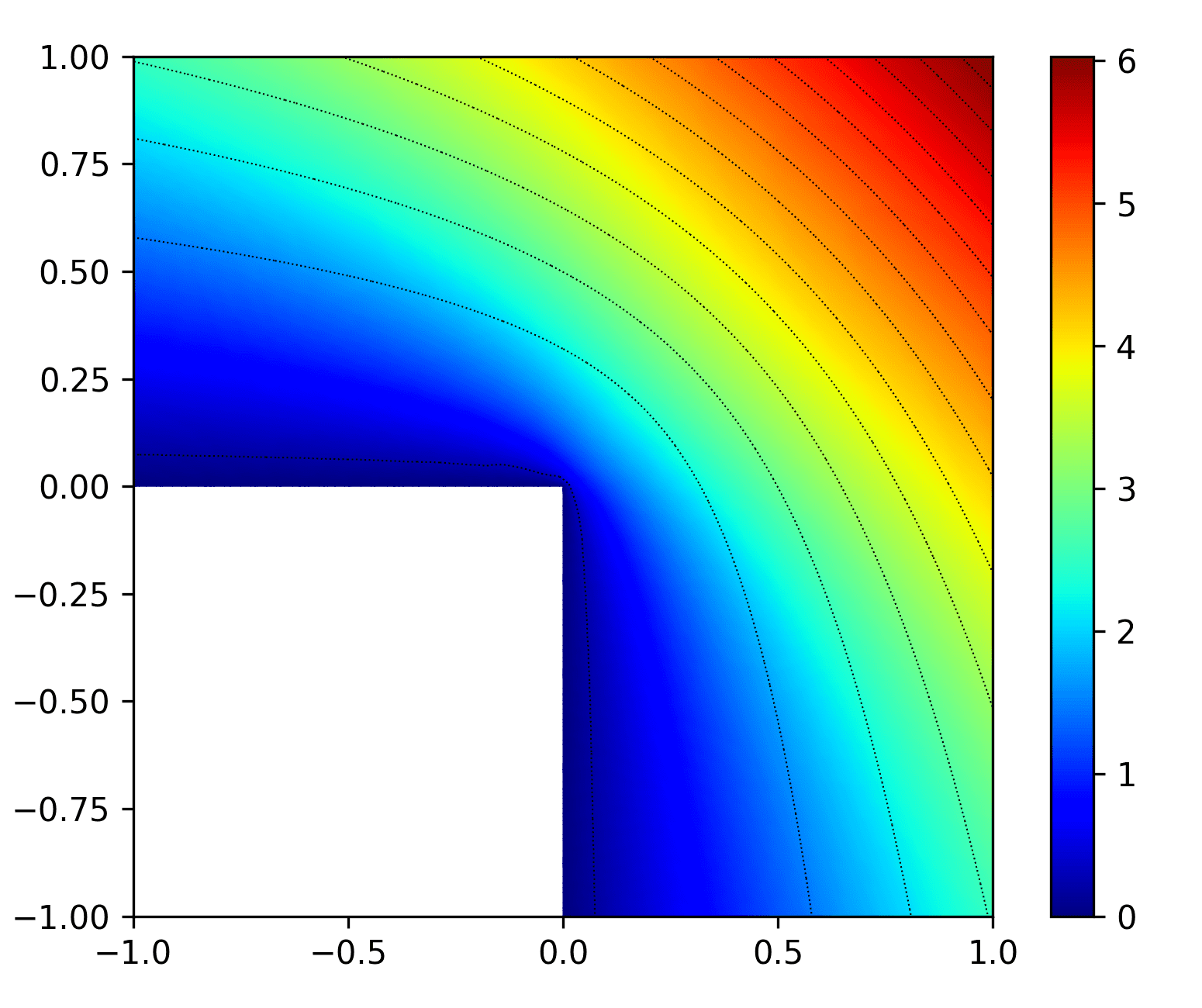}
		\caption{$| \uu |$}
		\label{fig:lshapevelocity}
	\end{subfigure}\hspace{0.08\textwidth}%
	\begin{subfigure}[H]{0.45\textwidth}
		\centering
		\includegraphics[width=1\textwidth]{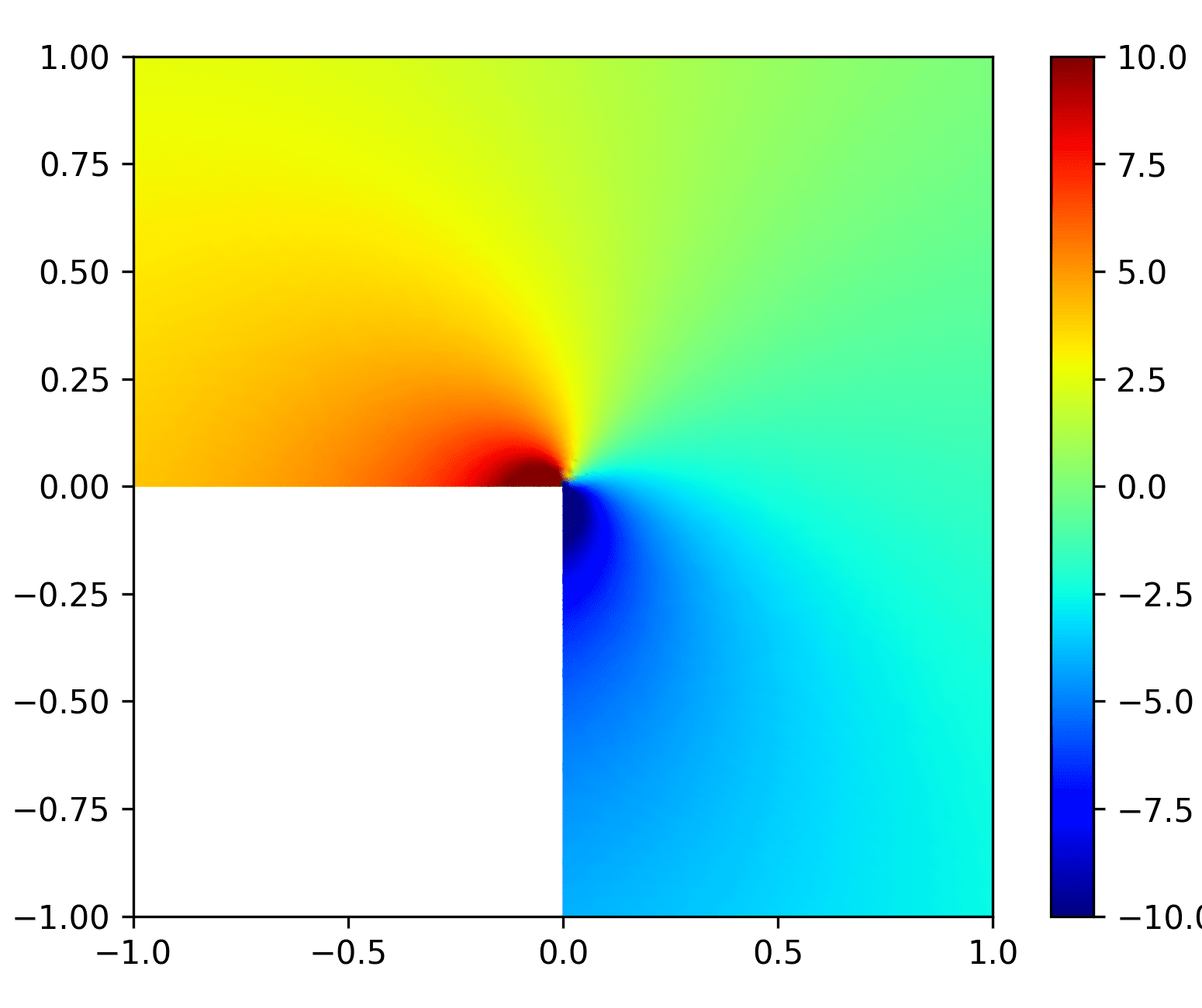}
		\caption{$p$}
		\label{fig:lshapepressure}
    \end{subfigure}
	\caption{(\subref{fig:lshapevelocity}) Velocity magnitude and streamlines, and (\subref{fig:lshapepressure}) pressure for the exact solution \eqref{equation:exactstokescorner} to the Stokes problem on the re-entrant corner domain. Because of the singular solution, the pressure color bar is truncated to the range $-10$ and $10$.}
	\label{fig:stokesonlshape}
\end{figure}

As a final benchmark problem we consider the Stokes problem \eqref{equation:stokesequations} on the re-entrant corner domain with mixed Dirichlet and Neumann boundaries introduced above, as shown in Figure \ref{fig:Lshape}. The weakly singular exact solution is taken from Ref.~\cite{verfurth1996} as
\begin{equation}\label{equation:exactstokescorner}
\begin{array}{l l}
u_1 = \phantom{-}R^\alpha \big{[}\sin(\theta) \frac{\partial \psi}{\partial \theta} -(1 + \alpha) \cos(\theta) \psi \big{]}\\
u_2 = -R^\alpha \big{[}\cos(\theta) \frac{\partial \psi}{\partial \theta} + (1 + \alpha) \sin(\theta) \psi \big{]}\\
\end{array} \qquad p= -\frac{R^{\alpha-1}}{1-\alpha} \Big{[}(1 + \alpha)^2 \frac{\partial \psi}{\partial \theta} + \frac{\partial^3 \psi}{\partial \theta^3} \Big{]},
\end{equation}
with constants $\alpha = 856399/1572864$ and $\omega = \frac{3}{2}\pi$, and with
\begin{equation}
\begin{split}
\psi(\theta) =& \frac{\cos(\alpha \omega)}{1 + \alpha}\sin((1 + \alpha) \theta)  - \frac{\cos(\alpha \omega)}{1 - \alpha} \sin((1 - \alpha) \theta) + \cos((1 - \alpha) \theta) - \cos((1 + \alpha) \theta).
\end{split}
\end{equation}
The exact pressure and velocity fields are illustrated in Figure~\ref{fig:stokesonlshape}. The corresponding Stokes problem \eqref{equation:stokesequations} is considered with the viscosity set to $\mu=1$, no body force, $\ff=\mathbf{0}$, a no slip condition on $\Gamma_D$, such that $\uu_D=\mathbf{0}$, and the Neumann data $\g$ on $\Gamma_N$ matching the exact solution.

\begin{figure}
	\centering
	\begin{subfigure}[b]{0.45\textwidth}
		\centering
		\includegraphics[width=\textwidth]{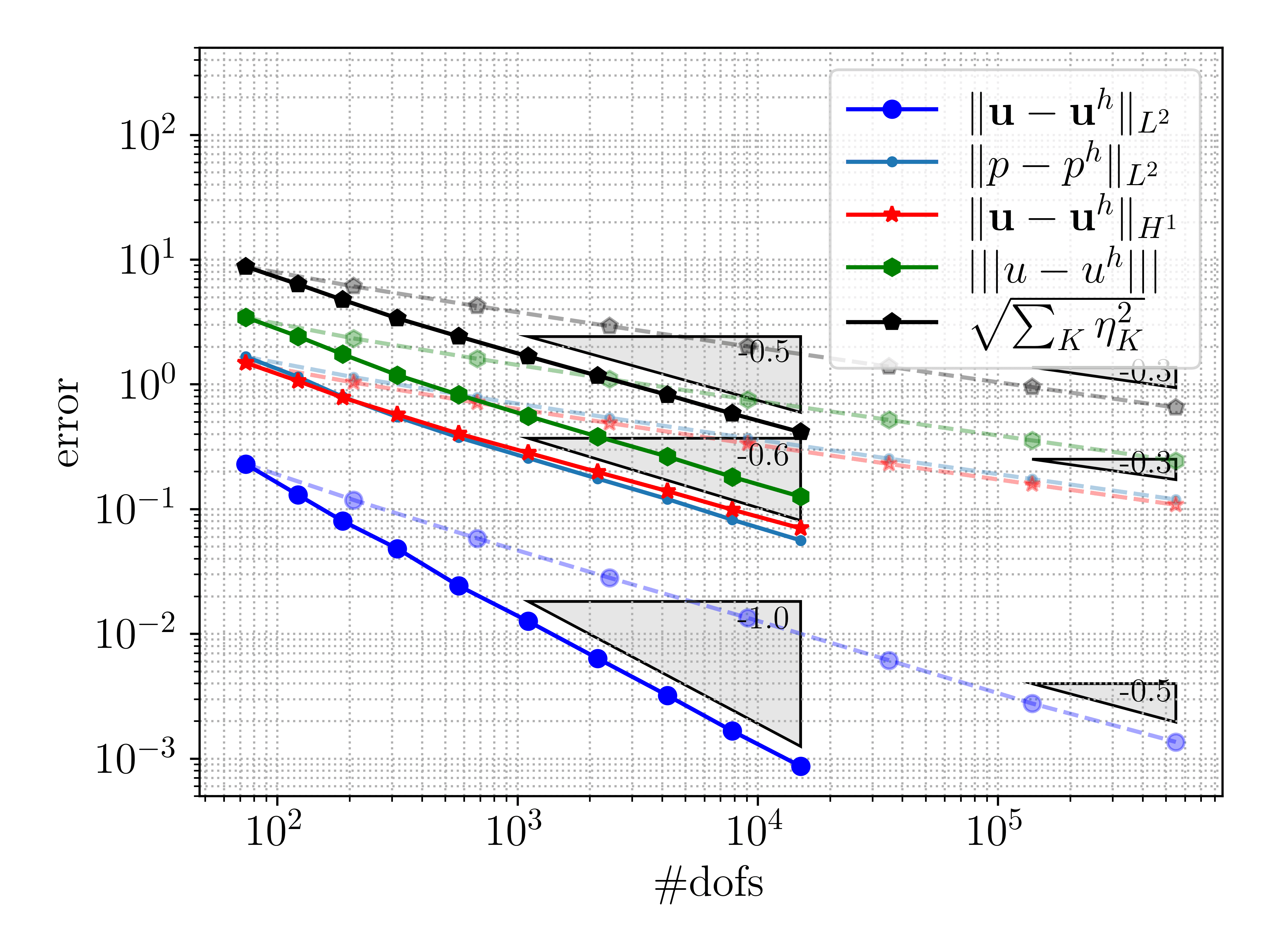}
		\caption{$\order = 1$}
		\label{fig:Lshape_stokes_p1}
	\end{subfigure}\hspace{0.08\textwidth}%
	\begin{subfigure}[b]{0.45\textwidth}
		\centering
		\includegraphics[width=\textwidth]{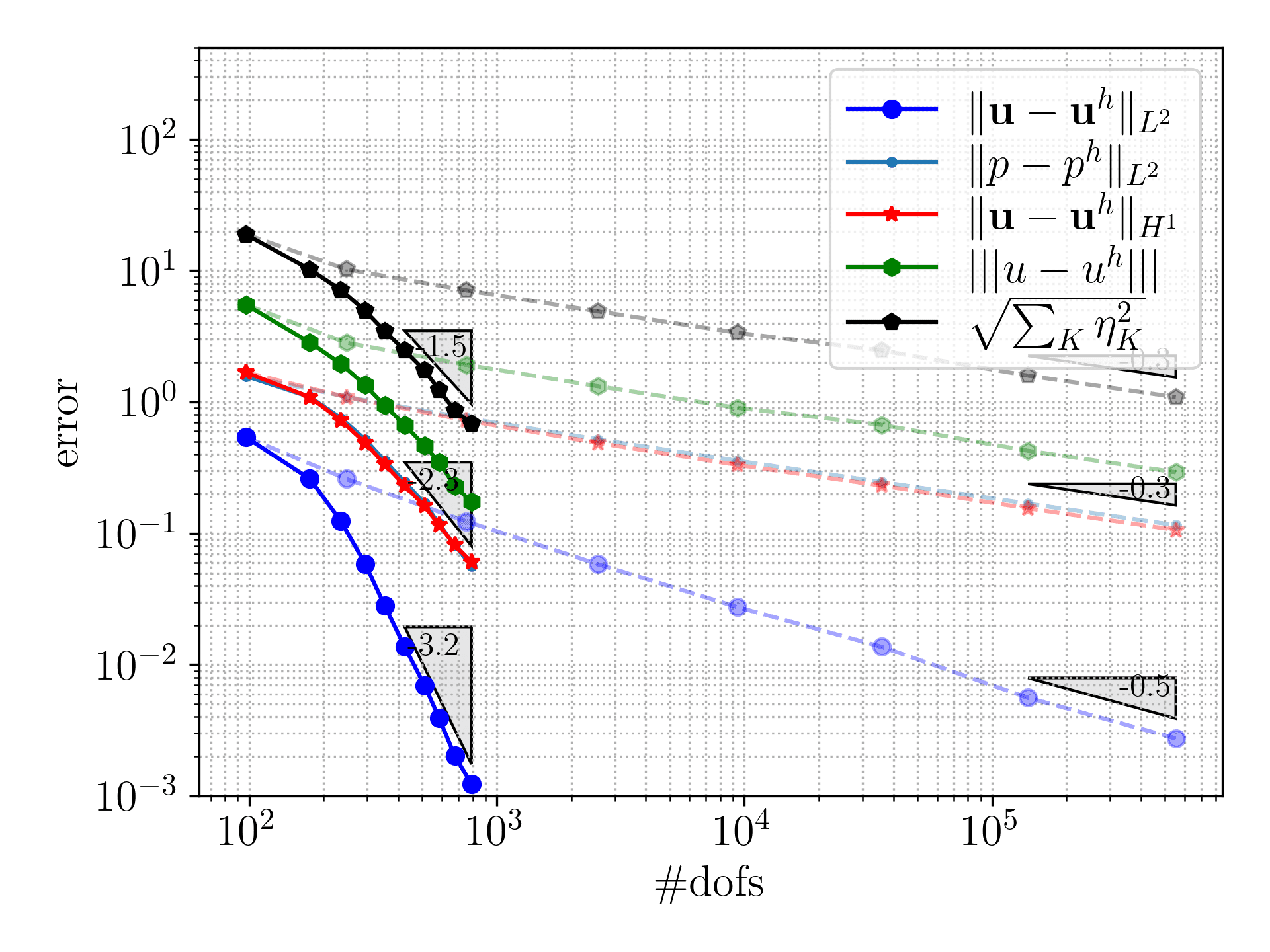}
		\caption{$\order = 2$}
		\label{fig:Lshape_stokes_p2}
	\end{subfigure}
	\hspace{0.08\textwidth}
		\caption{Error convergence results for the Stokes problem on the re-entrant corner domain under residual-based adaptive refinement (solid) and uniform refinement (dashed) for linear ($\order =1$) and quadratic ($\order =2$) basis functions.}
	\label{fig:Lshape_stokes_conv}
\end{figure}

\begin{figure}
	\centering
	\begin{subfigure}[b]{0.45\textwidth}
		\centering
		\includegraphics[width=0.8\textwidth]{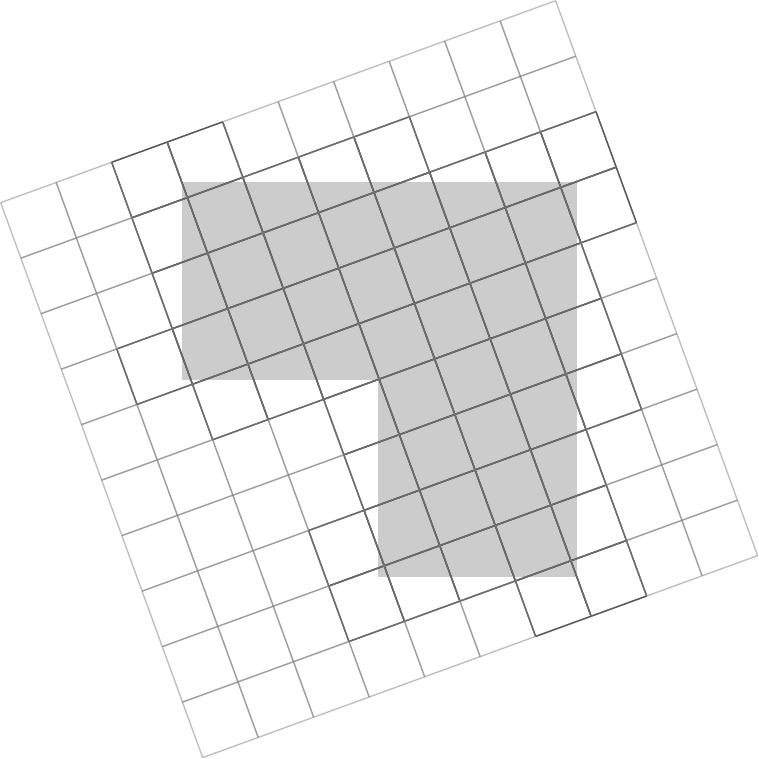}
		\caption{Initial mesh}
	\end{subfigure}	\hspace{0.08\textwidth}%
	\begin{subfigure}[b]{0.45\textwidth}
		\centering
		\includegraphics[width=0.8\textwidth]{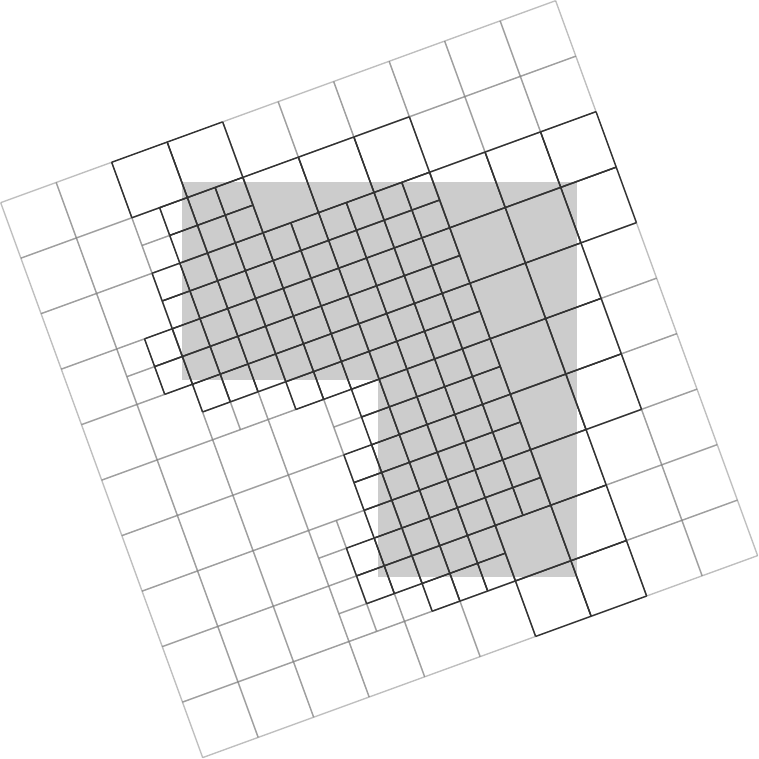}
		\caption{Step $1$}
	\end{subfigure}
	\begin{subfigure}[b]{0.45\textwidth}
		\centering
		\includegraphics[width=0.8\textwidth]{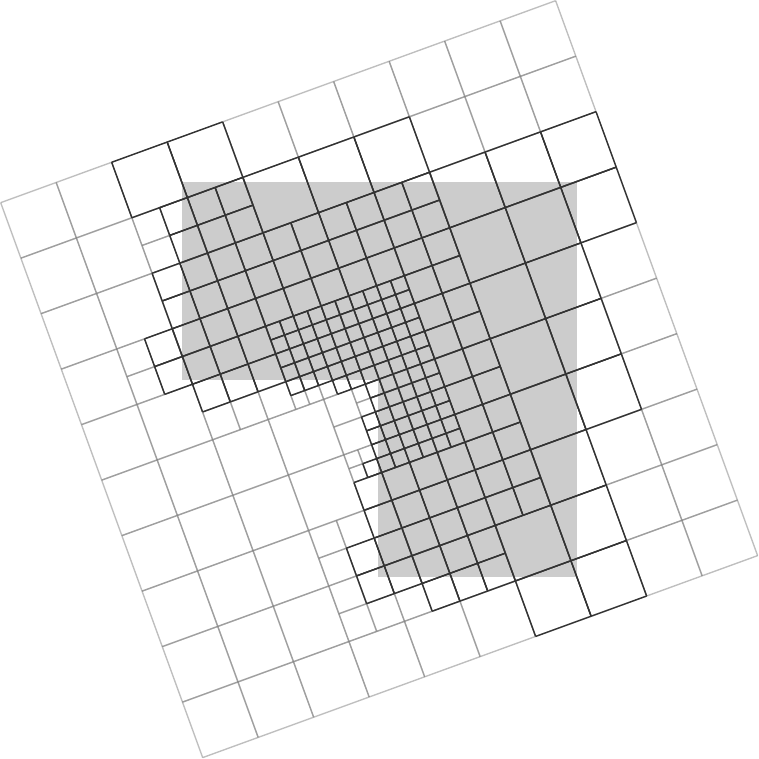}
		\caption{Step $2$}
	\end{subfigure}	\hspace{0.08\textwidth}%
	\begin{subfigure}[b]{0.45\textwidth}
		\centering
		\includegraphics[width=0.8\textwidth]{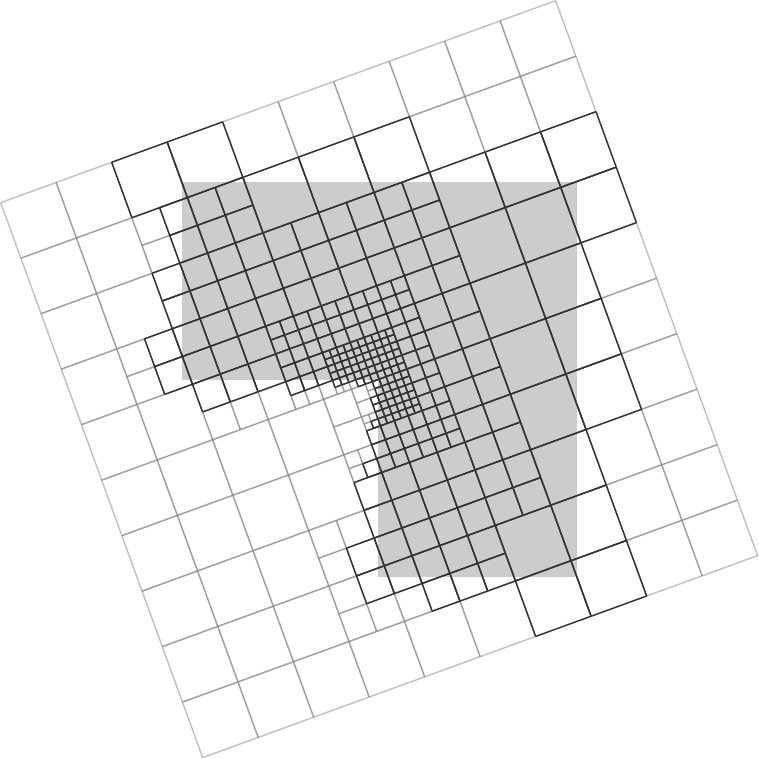}
		\caption{Step $3$}
	\end{subfigure}
	\begin{subfigure}[b]{0.45\textwidth}
		\centering
		\includegraphics[width=0.8\textwidth]{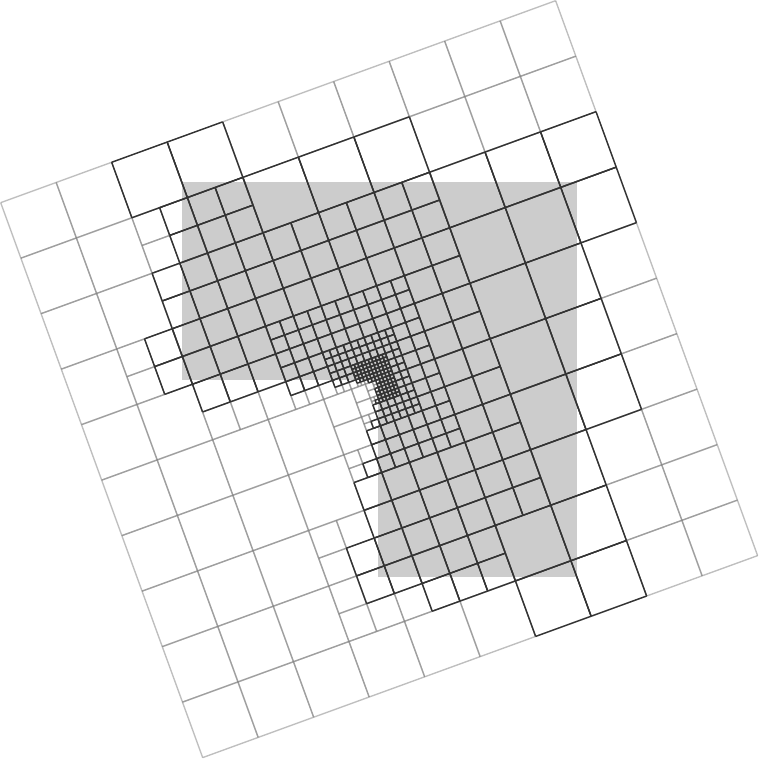}
		\caption{Step $4$}
	\end{subfigure}	\hspace{0.08\textwidth}%
	\begin{subfigure}[b]{0.45\textwidth}
		\centering
		\includegraphics[width=0.8\textwidth]{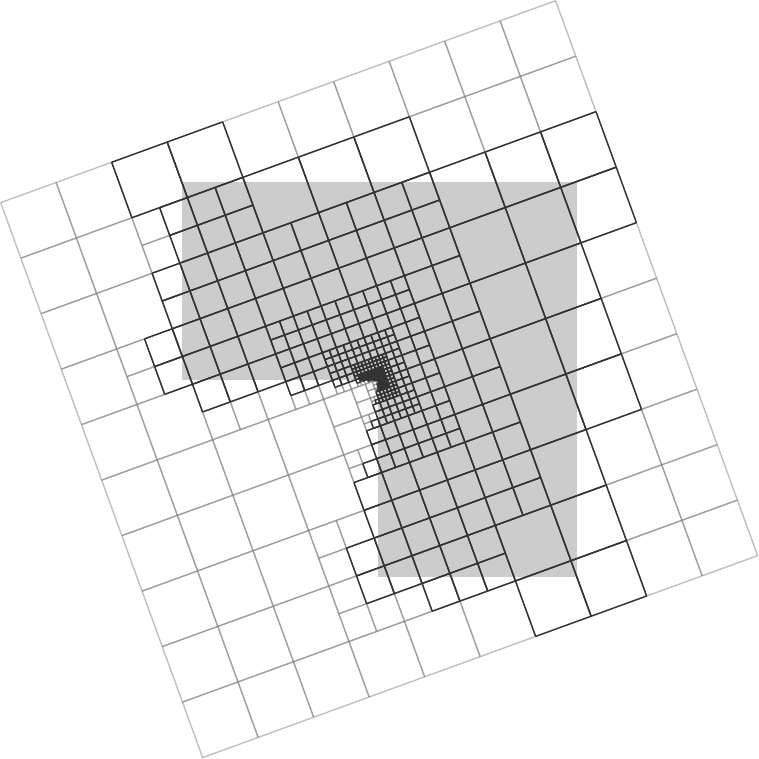}
		\caption{Step $5$}
	\end{subfigure}
	\caption{Evolution of the mesh using the adaptive refinement procedure for the Stokes problem on the re-entrant corner domain  using $\order = 2$.}
	\label{fig:Lshape_stokes_adap_mesh}
\end{figure}

Figure~\ref{fig:Lshape_stokes_conv} displays the error convergence results obtained using uniform and adaptive refinements, for both linear and quadratic (TH)B-splines. As for the Laplace case, the weak singularity in the exact solution \eqref{equation:exactstokescorner} limits the convergence rate when uniform refinements are considered. Using adaptive mesh refinement results in a recovery of the optimal rates in the case of linear basis functions, with even higher rates observed for the quadratic splines on account of the highly-focussed refinements resulting from the residual-based error estimator as observed in Figure~\ref{fig:Lshape_stokes_adap_mesh}.

\section{Scan-based simulations} \label{sec:scanbased}
In this section we apply the developed adaptive immersed isogeometric analysis framework in the context of scan-based analysis. We consider the viscous flow problem on a two-dimensional image domain and on a three-dimensional patient-specific problem based on a $\mu$CT-scan of a carotid artery, represented by grayscale voxels. The primary purpose of the two-dimensional setting is to test the scan-based analysis framework. For all simulations, the octree subdivision depth is set equal to $8$ in two dimensions and $3$ in three dimensions. The refinement threshold related to the D\"orfler marking is set to $\lambda=0.8$.

\begin{figure}
	\centering
	\begin{subfigure}[b]{0.3\textwidth}
		\centering
		\includegraphics[width=0.8\textwidth]{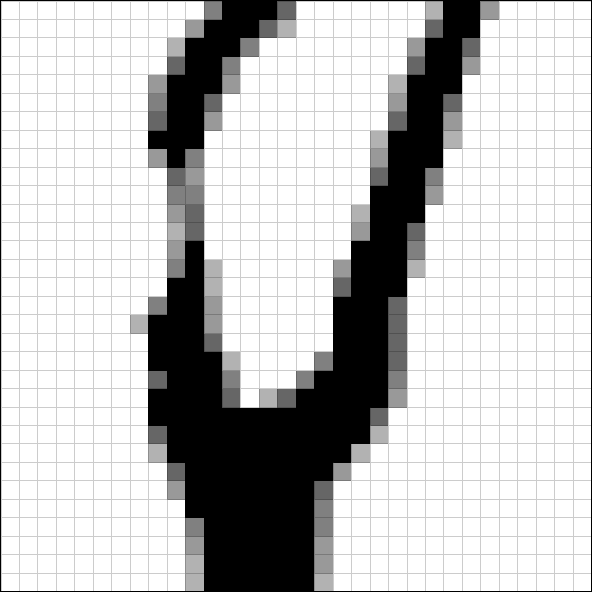}
		\caption{}
		\label{fig:voxeldata}
	\end{subfigure}\hfill%
	\begin{subfigure}[b]{0.3\textwidth}
		\centering
		\includegraphics[width=\textwidth]{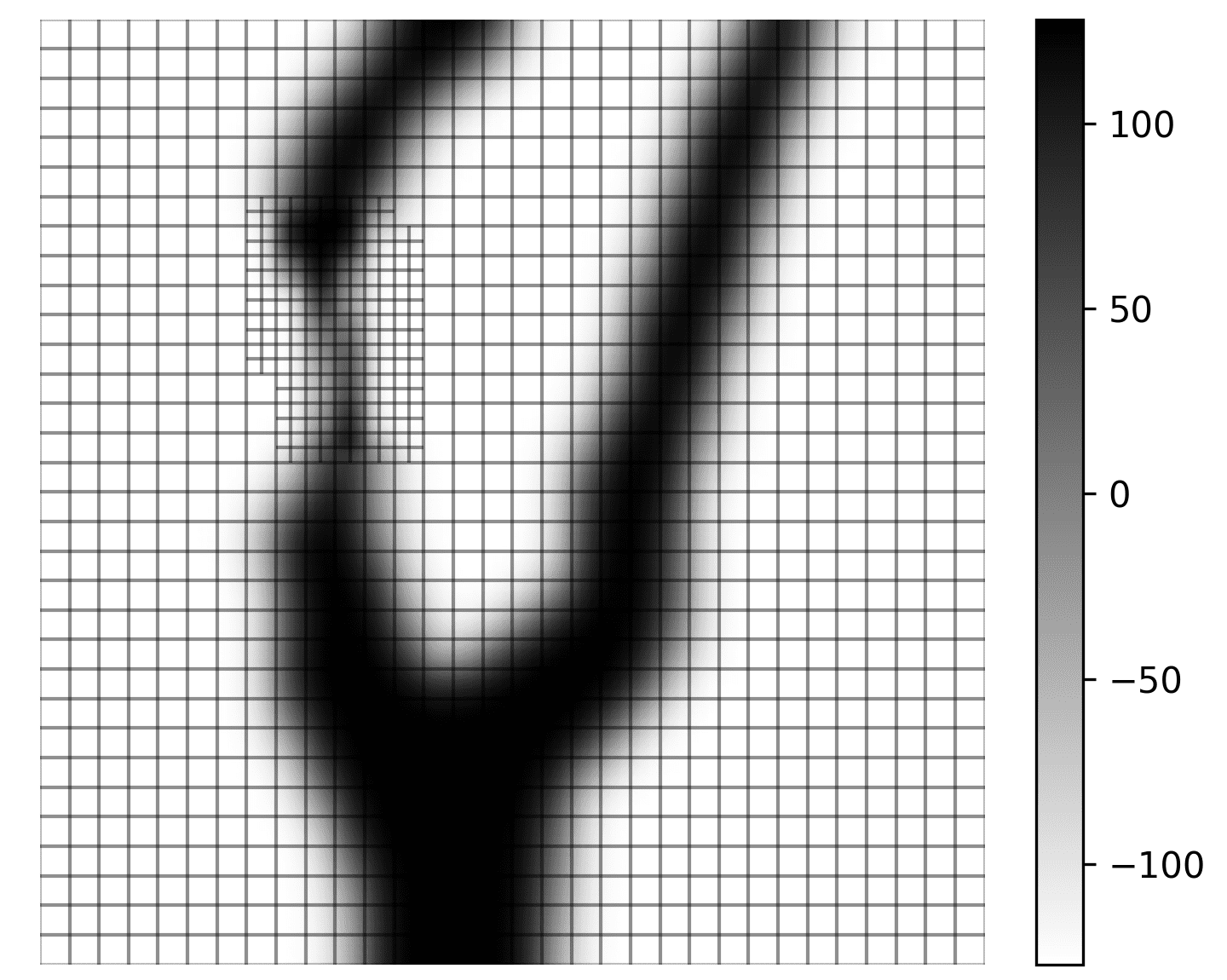}
		\caption{}
		\label{fig:levelset}
	\end{subfigure}\hfill%
	\begin{subfigure}[b]{0.3\textwidth}
		\centering
		\includegraphics[width=0.8\textwidth]{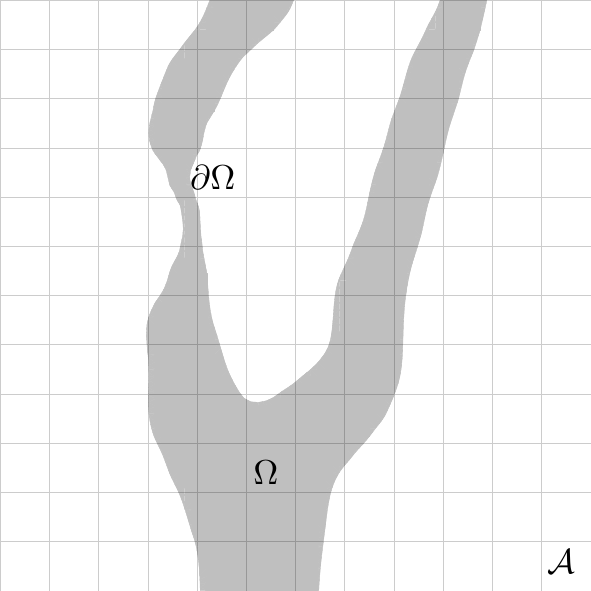}
		\caption{}
		\label{fig:trimmedtopo}
	\end{subfigure}
	\caption{Illustration of the scan-based analysis workflow. The original grayscale image in panel (\subref{fig:voxeldata}) is converted to a level set function, shown in panel (\subref{fig:levelset}), which is constructed using the topology-preserving segmentation algorithm of Ref.~\cite{divi2021}. The  trimmed geometry, shown in panel (\subref{fig:trimmedtopo}), is then extracted using the recursive bi-sectioning strategy with mid-point tessellation of Ref.~\cite{divi2020}.}
	\label{fig:scanstrategy}
\end{figure}

Our scan-based analysis workflow is illustrated in Figure~\ref{fig:scanstrategy}. The first step in this workflow is to smoothen the original grayscale voxel data using a convolution operation on a B-spline basis formed on the voxel grid \cite{verhoosel2015}. Since this smoothing operator behaves as a Gaussian filter, geometric features that are similar in size to the voxels can be lost \cite{divi2021}. To avoid this loss of features, the topology-preservation procedure proposed in Ref.~\cite{divi2021} is employed. This procedure locally refines the convolution basis to retain small geometric features in the smoothing procedure. Once the smooth level set representation has been obtained, the octree segmentation procedure with mid-point tessellation of Ref.~\cite{divi2020} is used to obtain the immersed geometry represented on an ambient domain mesh. It is important to note that this ambient domain mesh, on which the  solution to the flow problem is computed, can be chosen independently of the voxel size, and hence it is independent of the mesh on which the level set function is constructed.

\begin{figure}
	\centering
	\begin{subfigure}[b]{0.45\textwidth}
		\centering
		\includegraphics[width=\textwidth]{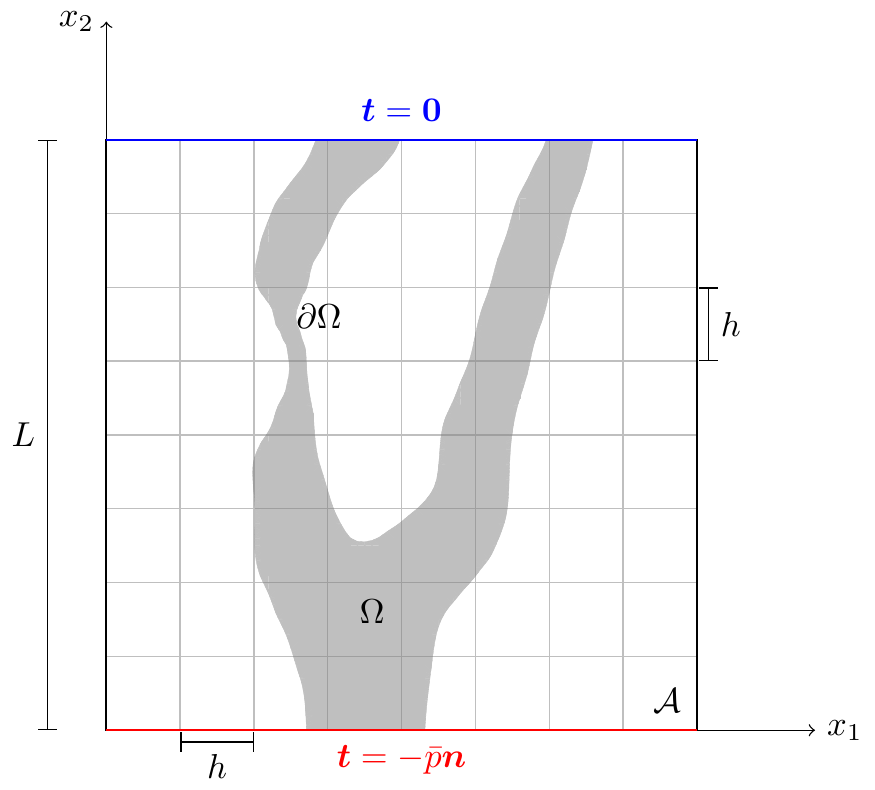}
		\caption{ }
		\label{fig:2dstokes_domain}
	\end{subfigure}\hspace{0.08\textwidth}%
	\begin{subfigure}[b]{0.45\textwidth}
		\centering
			\includegraphics[width=\textwidth]{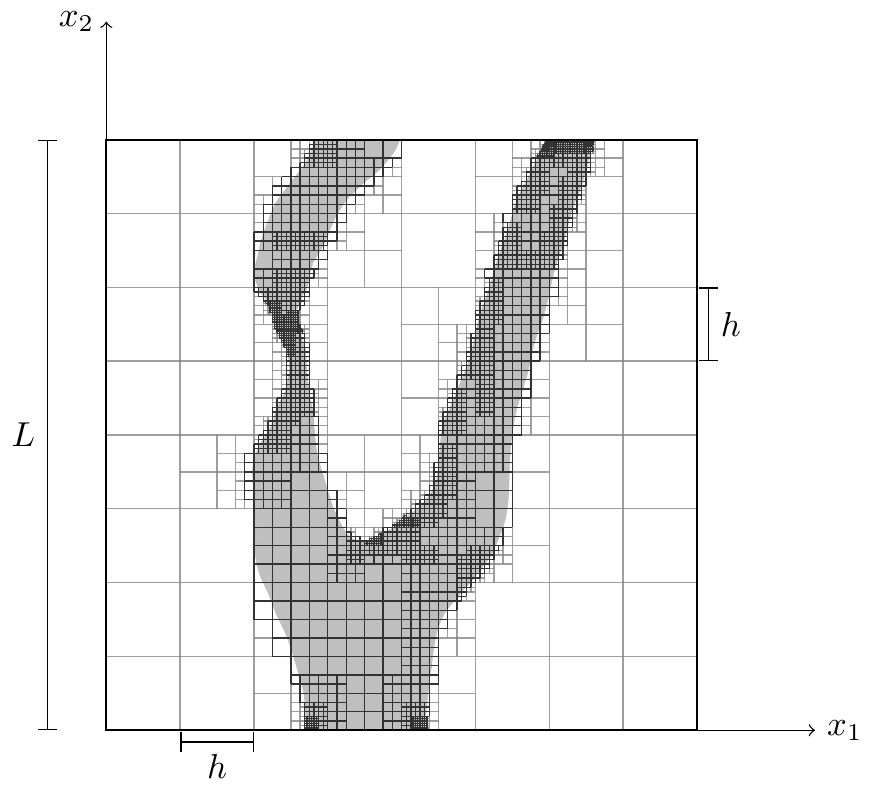}
		\caption{ }
		\label{fig:2dstokes_mesh}
	\end{subfigure}
	\caption{(\subref{fig:2dstokes_domain}) Illustration of the domain and boundary conditions for the scan-based viscous flow problem, and (\subref{fig:2dstokes_mesh}) a typical locally refined mesh resulting from the adaptive procedure.}
	\label{fig:2dstokes}
\end{figure}

The considered computational domain is illustrated in Figure~\ref{fig:2dstokes_domain}. Neumann conditions are imposed on the inflow and outflow boundaries, with the traction on the inflow boundary acting in the normal direction with a traction data, $\boldsymbol{t}=-\bar{p} \boldsymbol{n}$, where $\bar{p}$ is the pressure magnitude. Homogeneous Dirichlet conditions are imposed along the immersed boundaries in accordance with the no slip condition. It is to be noted that a Neumann condition at an inflow boundary generally leads to an ill-posed boundary value problem for the Navier-Stokes equations, but the Stokes problem is well-posed. 
In all simulations we consider second-order ($\order = 2$) (TH)B-splines  and set the stabilization parameters to $\beta = 100$, $\gamma_{g} = 10^{-(\order + 2)}$ and $\gamma_{s}=10^{-(\order+1)}$, which have been determined empirically.
 
\subsection{Two-dimensional prototypical geometry}\label{sec:2dstokes}
To test the developed methodology in the scan-based setting, we first consider the prototypical two-dimensional geometry shown in Figure~\ref{fig:2dstokes_domain}, which is constructed from $32 \times 32$ grayscale voxel data. The ambient domain, which matches the scan window, is taken as a unit square ($L=1$) which is covered by an $8 \times 8$ elements ambient mesh. The viscosity is set equal to $\mu=1$ and the pressure to $\bar{p} = 1$.

\begin{figure}
	\begin{subfigure}[b]{0.45\textwidth}
		\centering
		\includegraphics[width=\textwidth]{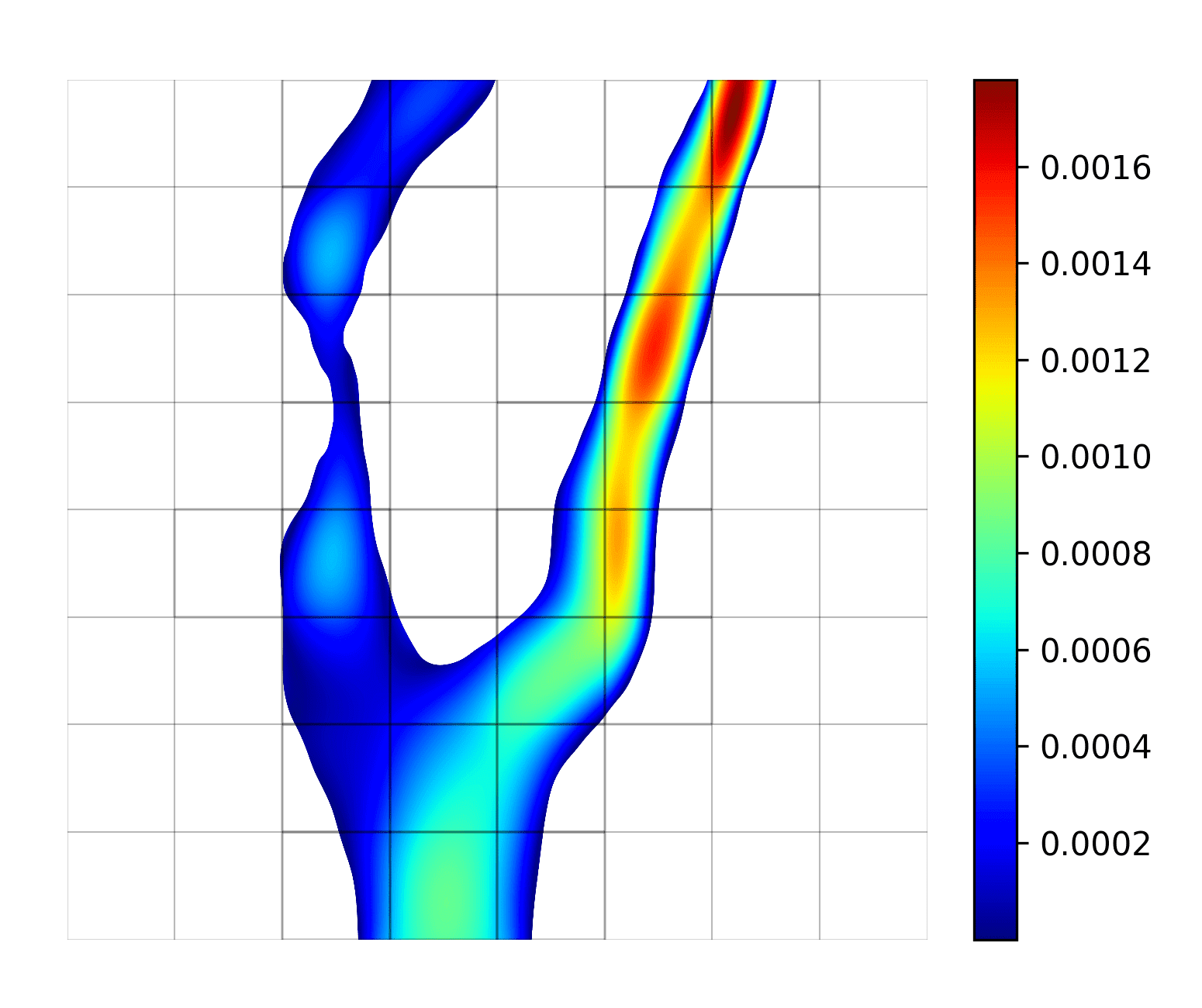}
		\caption{Initial mesh}
		\label{fig:2dstokes_velocity_step1}
	\end{subfigure}\hspace{0.08\textwidth}%
	\begin{subfigure}[b]{0.45\textwidth}
		\centering
			\includegraphics[width=\textwidth]{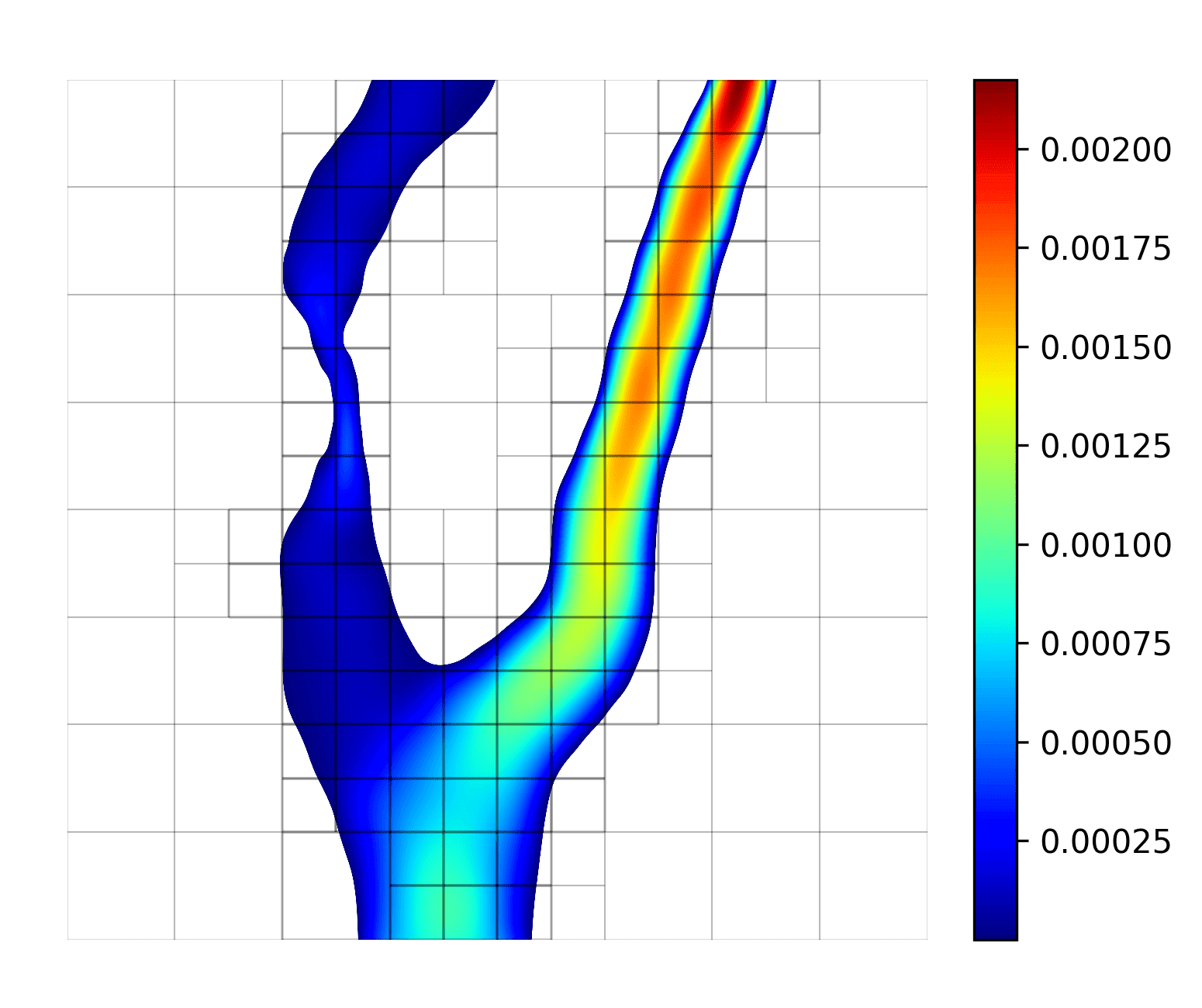}
		\caption{Step $1$}
		\label{fig:2dstokes_velocity_step2}
	\end{subfigure}\\[12pt]
	\begin{subfigure}[b]{0.45\textwidth}
		\centering
		\includegraphics[width=\textwidth]{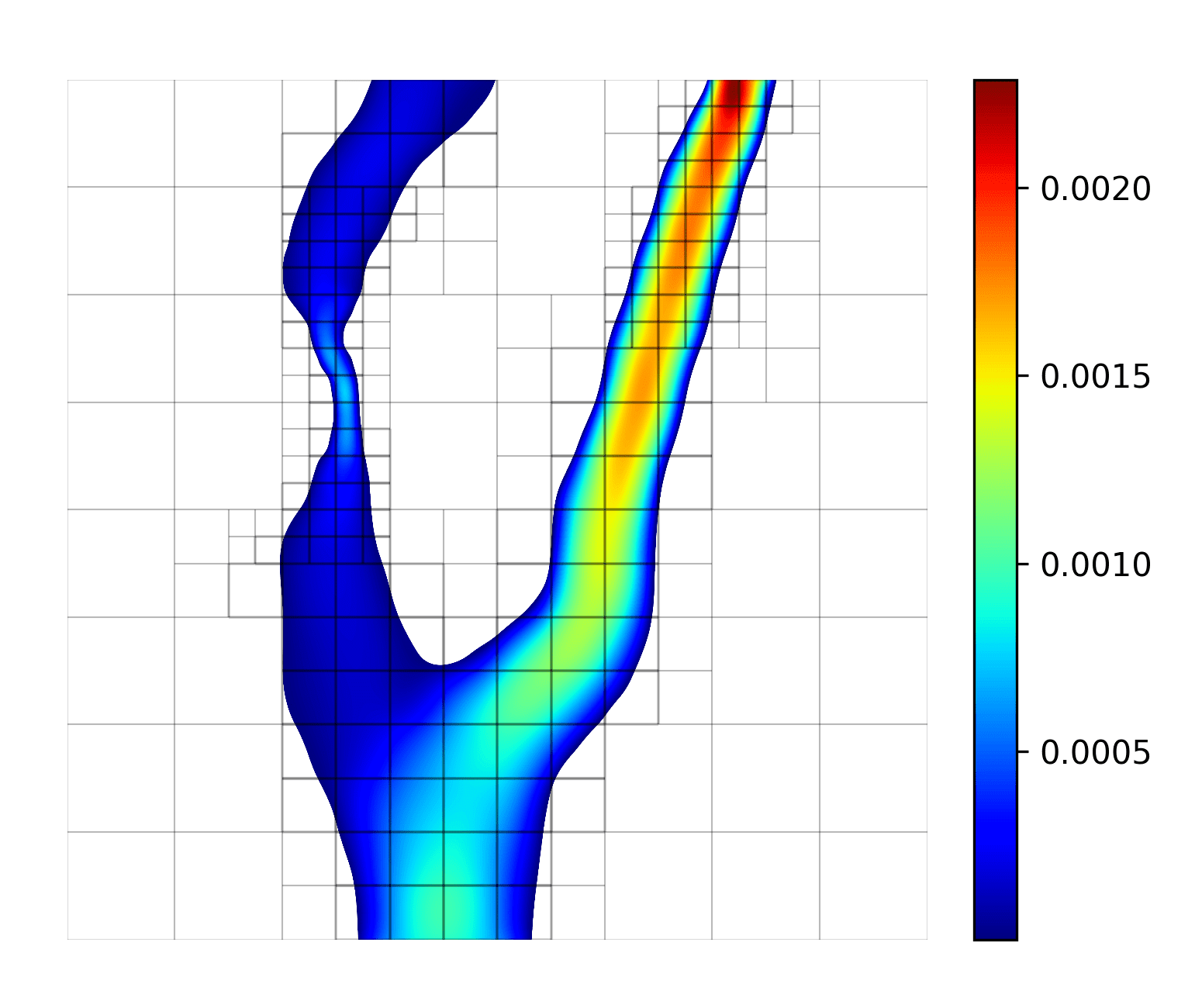}
		\caption{Step $2$}
		\label{fig:2dstokes_velocity_step3}
	\end{subfigure}\hspace{0.08\textwidth}%
	\begin{subfigure}[b]{0.45\textwidth}
		\centering
		\includegraphics[width=\textwidth]{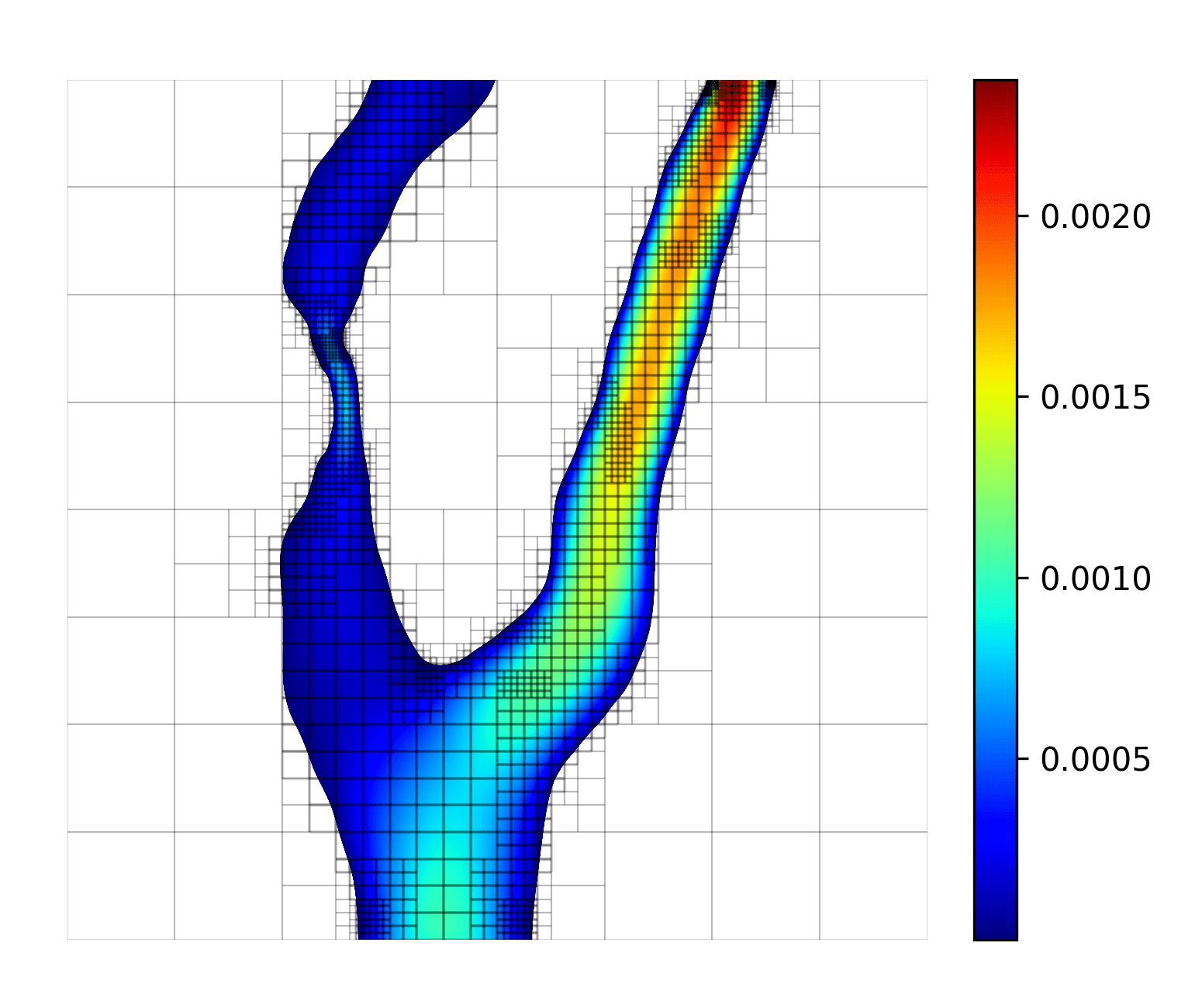}
		\caption{Step $6$}
		\label{fig:2dstokes_velocity_step6}
	\end{subfigure}
	\caption{Evolution of the mesh and (magnitude of the) velocity field during the adaptive refinement process for the viscous flow in two dimensions.}
	\label{fig:2dstokes_velocity}
\end{figure}

Various steps in the adaptive refinement procedure are depicted in Figure~\ref{fig:2dstokes_velocity}. In the first step virtually all elements covering the flow domain are refined, indicating that the initial mesh of only $8 \times 8$ elements is too coarse to resolve the solution globally. After the first refinement step, the refinement strategy starts to focus on the regions where the errors are largest, \emph{i.e.}, near boundaries and narrow sections, as also illustrated in Figure~\ref{fig:2dstokes_mesh}. Under further refinement, the procedure resolves prominent solution details, most importantly the (Poiseuille-like) profile in the carotid part of the artery and the velocity profiles at the inflow and outflow boundaries.

\begin{figure}
	\centering
	\begin{subfigure}[b]{0.45\textwidth}
		\centering
		\includegraphics[width=\textwidth]{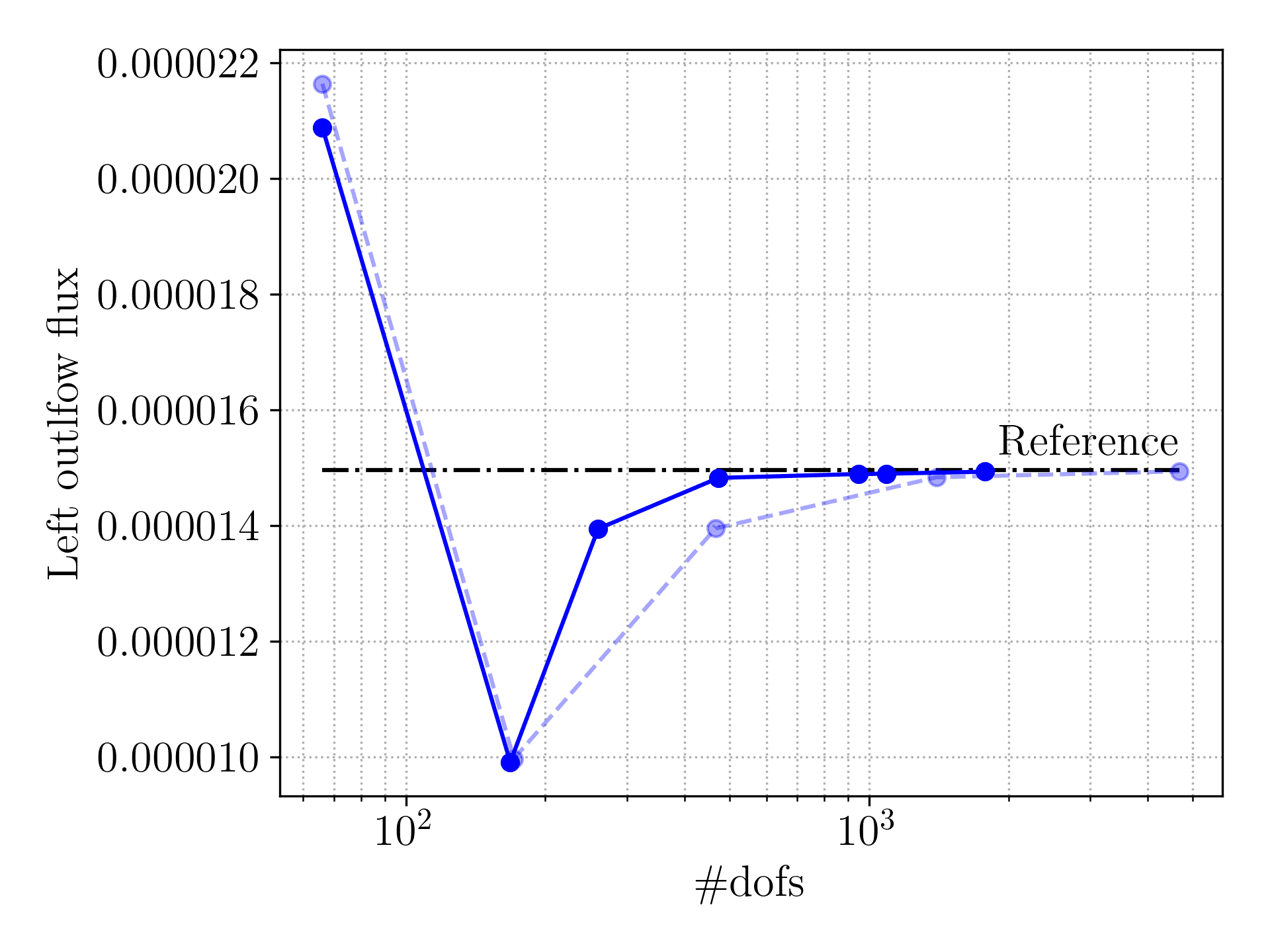}
		\caption{}\label{fig:outflowstokesleft}
	\end{subfigure}\hspace{0.08\textwidth}%
	\begin{subfigure}[b]{0.45\textwidth}
		\centering
		\includegraphics[width=\textwidth]{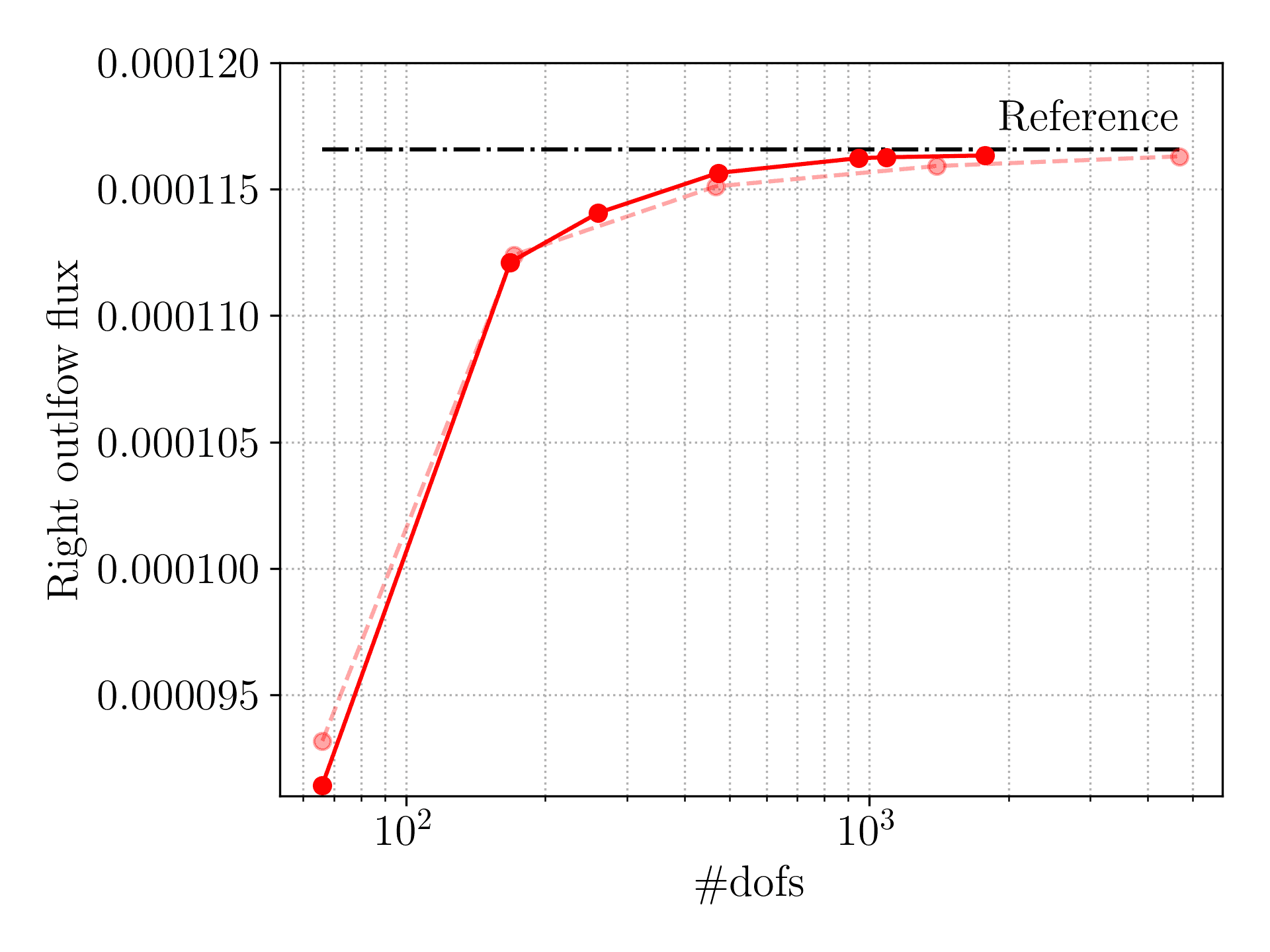}
		\caption{}\label{fig:outflowstokesright}
	\end{subfigure}
	\caption{Mesh convergence of the outflow flux at the (\subref{fig:outflowstokesleft}) left and (\subref{fig:outflowstokesright}) right channel of the domain in Figure~\ref{fig:2dstokes_domain} using adaptive (solid) and uniform (dashed) mesh refinements.}
	\label{fig:2dstokes_conv}
\end{figure}

Further results of the viscous flow problem solved using uniform and adaptive refinements are shown in Figure~\ref{fig:2dstokes_conv} in the form of the flux through the left and right outflow channels.  The minor difference in results on the initial mesh (left-most points) are caused by a different selection of the octree-depth for the uniform and adaptive simulations. Both methods are observed to converge to the same fluxes under refinement, but an excellent approximation of the reference solution (computed on a uniform overkill refinement, consistent with the result reported in Ref.~\cite{divi2021}) is obtained by means of the adaptive mesh refinement procedure using substantially fewer degrees of freedom than for uniform refinements. This is consistent with the observations on the velocity field discussed above, where in particular the ability of the adaptive refinement procedure to resolve the flow in the carotid part is essential.

\subsection{Three-dimensional patient-specific geometry} \label{sec:3dlaplace}
To demonstrate the residual-based adaptivity procedure in a real scan-based setting, we consider the patient-specific carotid artery used in Ref.~\cite{divi2021}. The geometry of the carotid artery is obtainted from CT-scan data containing 80 slices of $85 \times 70$ voxels. The size of each voxel is $300 \times 300\,{\rm \mu m^2}$ and the distance between the slices is $400\,{\rm \mu m}$. The total size of the scan domain is $25.6 \times 21.1 \times 32.0\, {\rm mm}^3$. We set the viscosity to $4$mPa s and pressure to $17.3$kPa ($130$ mm of Hg).

\begin{figure}
	\begin{subfigure}[b]{0.45\textwidth}
		\centering
		\includegraphics[width=\textwidth]{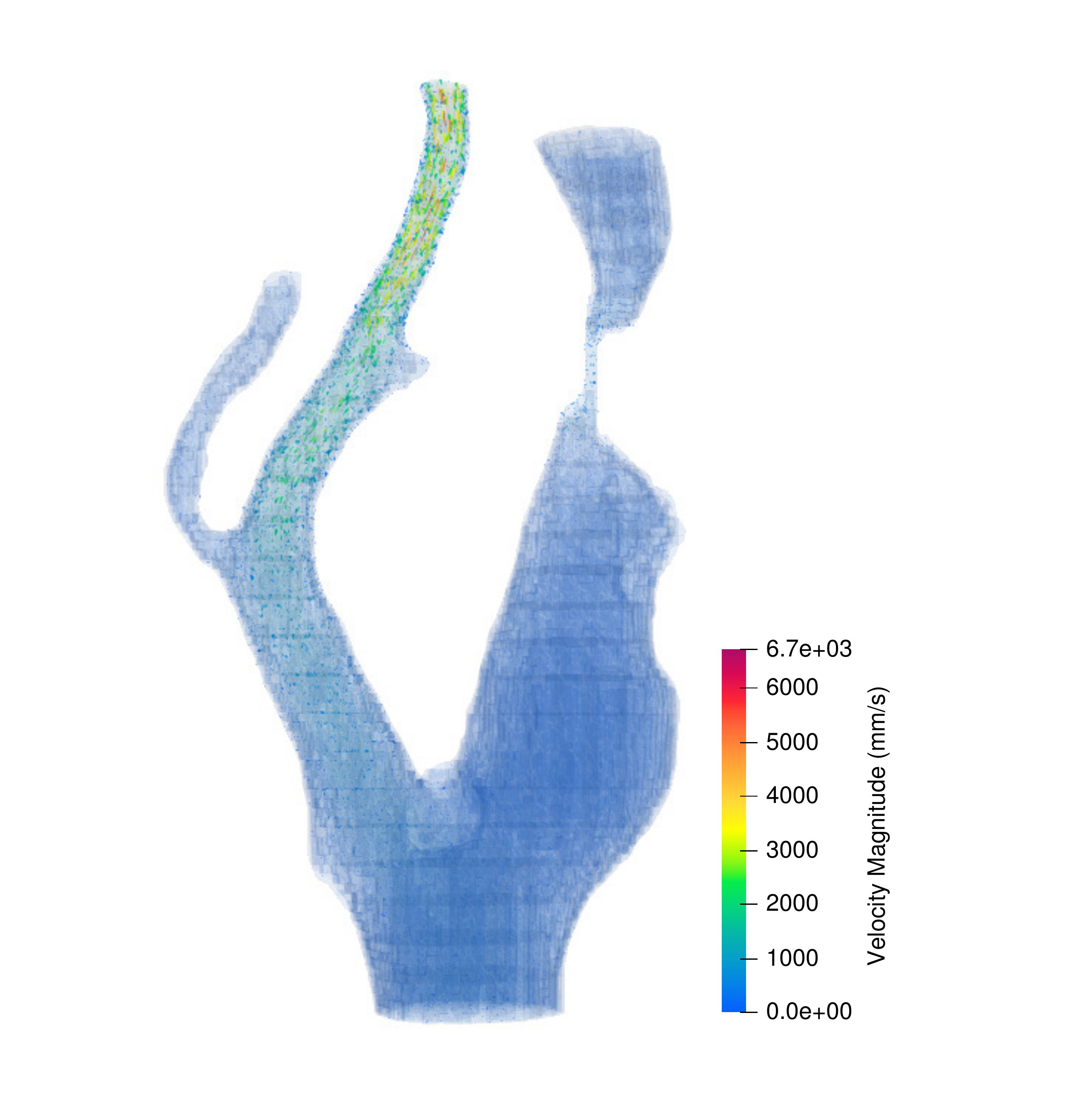}
		\caption{Initial stage}
		\label{fig:3dstokes_velocity_step0}
	\end{subfigure}\hspace{0.08\textwidth}%
	\begin{subfigure}[b]{0.45\textwidth}
		\centering
		\includegraphics[width=\textwidth]{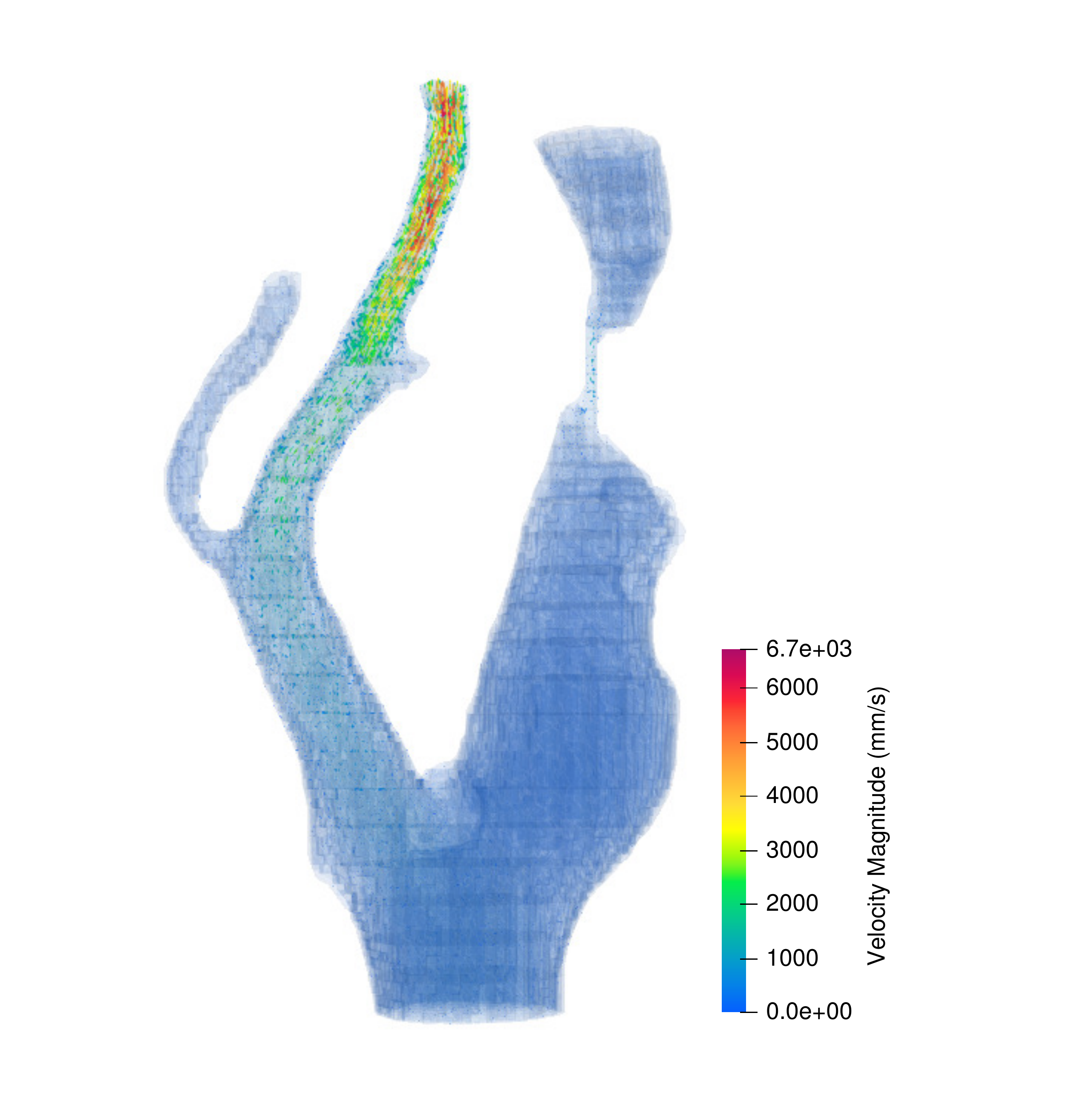}
		\caption{Step $1$}
		\label{fig:3dstokes_velocity_step1}
	\end{subfigure}\\[12pt]
	\begin{subfigure}[b]{0.45\textwidth}
		\centering
		\includegraphics[width=\textwidth]{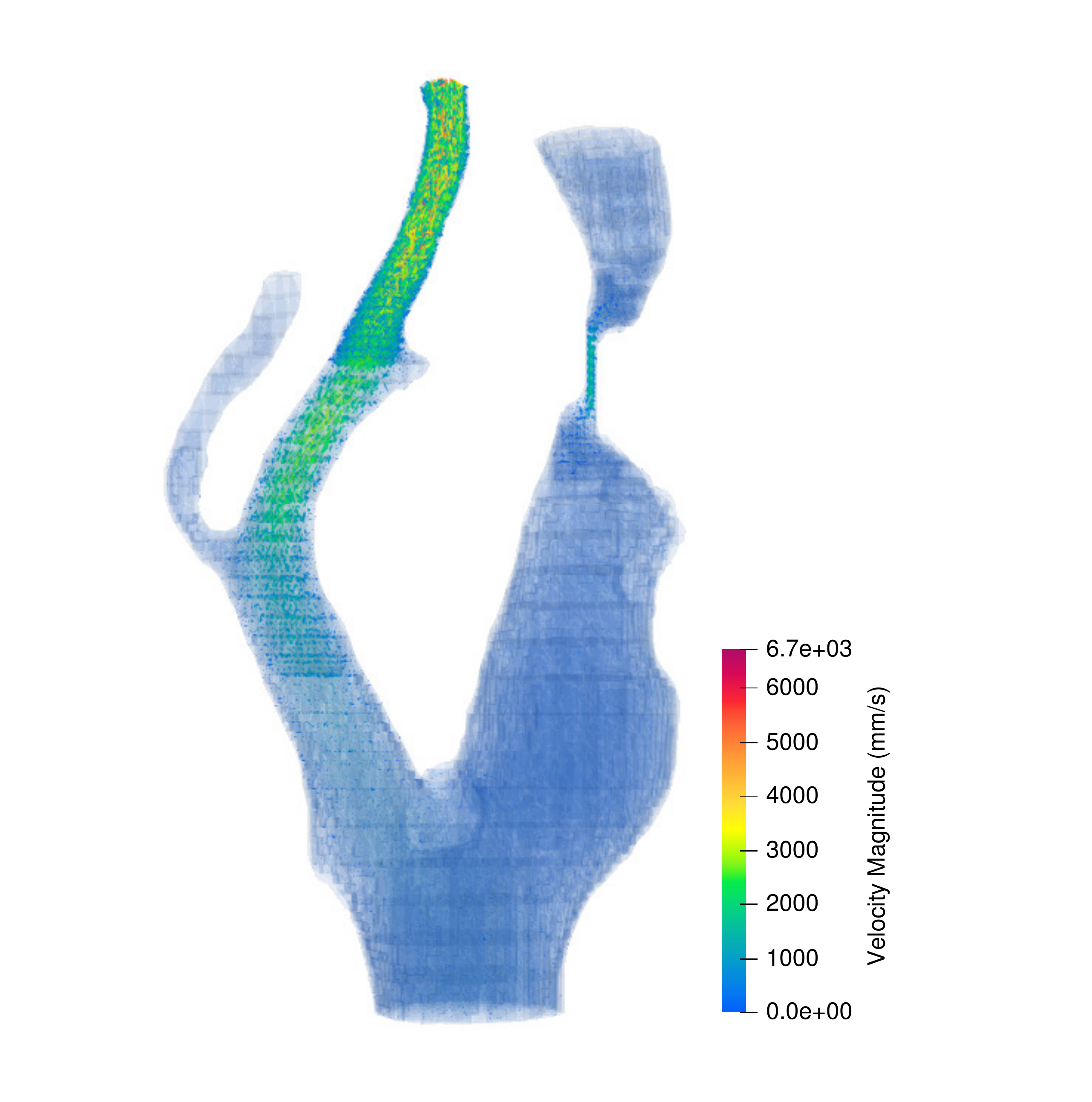}
		\caption{Step $2$}
		\label{fig:3dstokes_velocity_step2}
	\end{subfigure}\hspace{0.08\textwidth}%
	\begin{subfigure}[b]{0.45\textwidth}
		\centering
		\includegraphics[width=\textwidth]{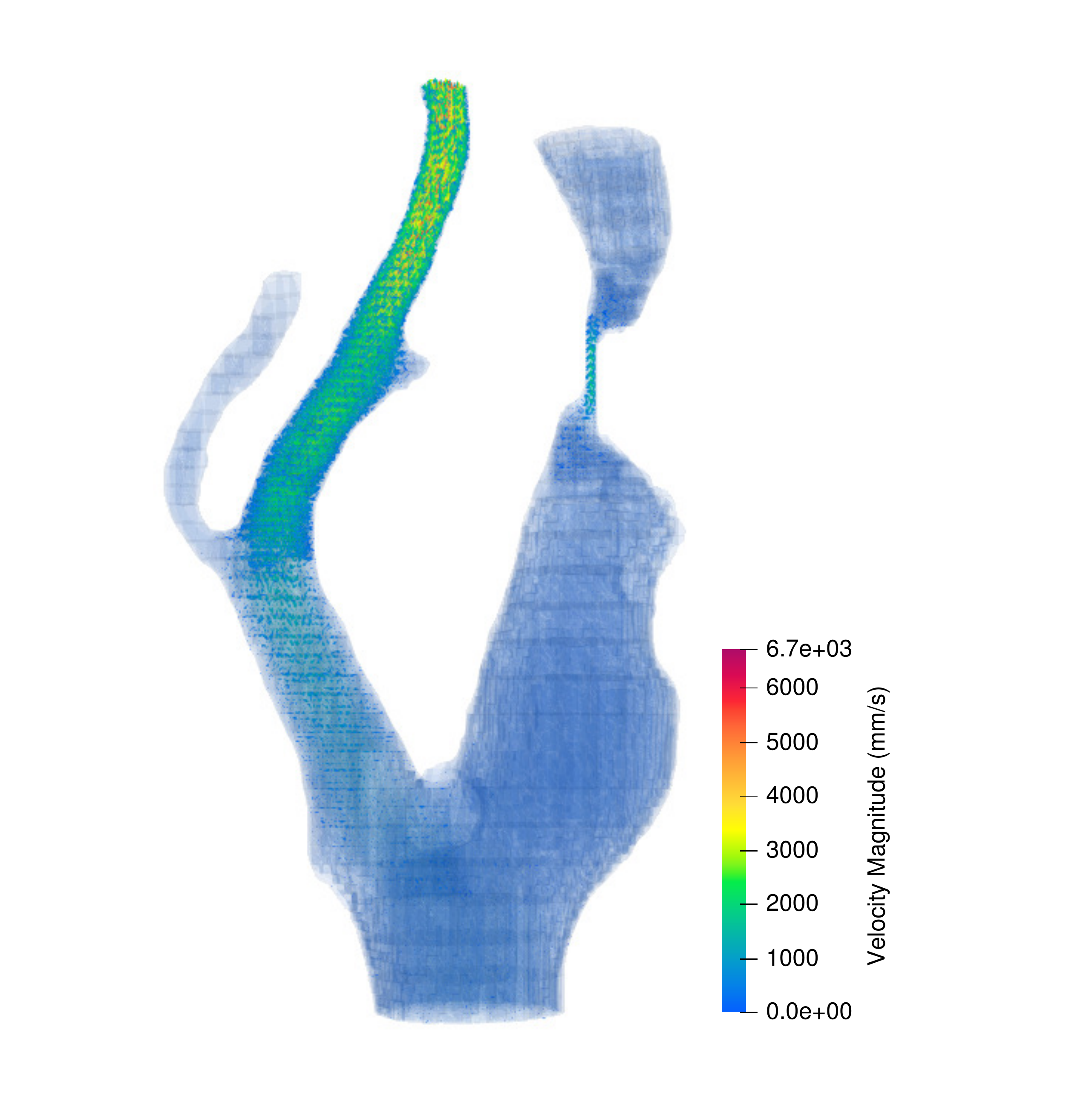}
		\caption{Step $3$}
		\label{fig:3dstokes_velocity_step3}
	\end{subfigure}
	\caption{Velocity magnitude during the adaptive refinement process for the patient-specific viscous flow problem.}
	\label{fig:3dstokes_velocity}
\end{figure}

Simulation results for this problem are shown in Figure~\ref{fig:3dstokes_velocity}. Note that for the considered scan data, the application of the topology-preservation algorithm in Ref.~\cite{divi2021} is essential, as otherwise the narrow channel section in the right artery would disappear. The simulation results are based on a $24 \times 24 \times 24$ ambient domain mesh of $25.6 \times 21.1 \times 32.0\, {\rm mm}^3$ and an octree depth of three. In this setting, after two refinements, an element is of a similar size as the voxels. The need to substantially refine beyond the voxel size is, from a practical perspective, questionable, as the dominant error in the analysis will then be related to the scan resolution and the segmentation procedure. In this sense, the constraint of not being able to refine beyond the octree depth is not a crucial problem in the considered simulations.

\begin{figure}
	\begin{subfigure}[b]{0.45\textwidth}
		\centering
   		\includegraphics[width=\textwidth]{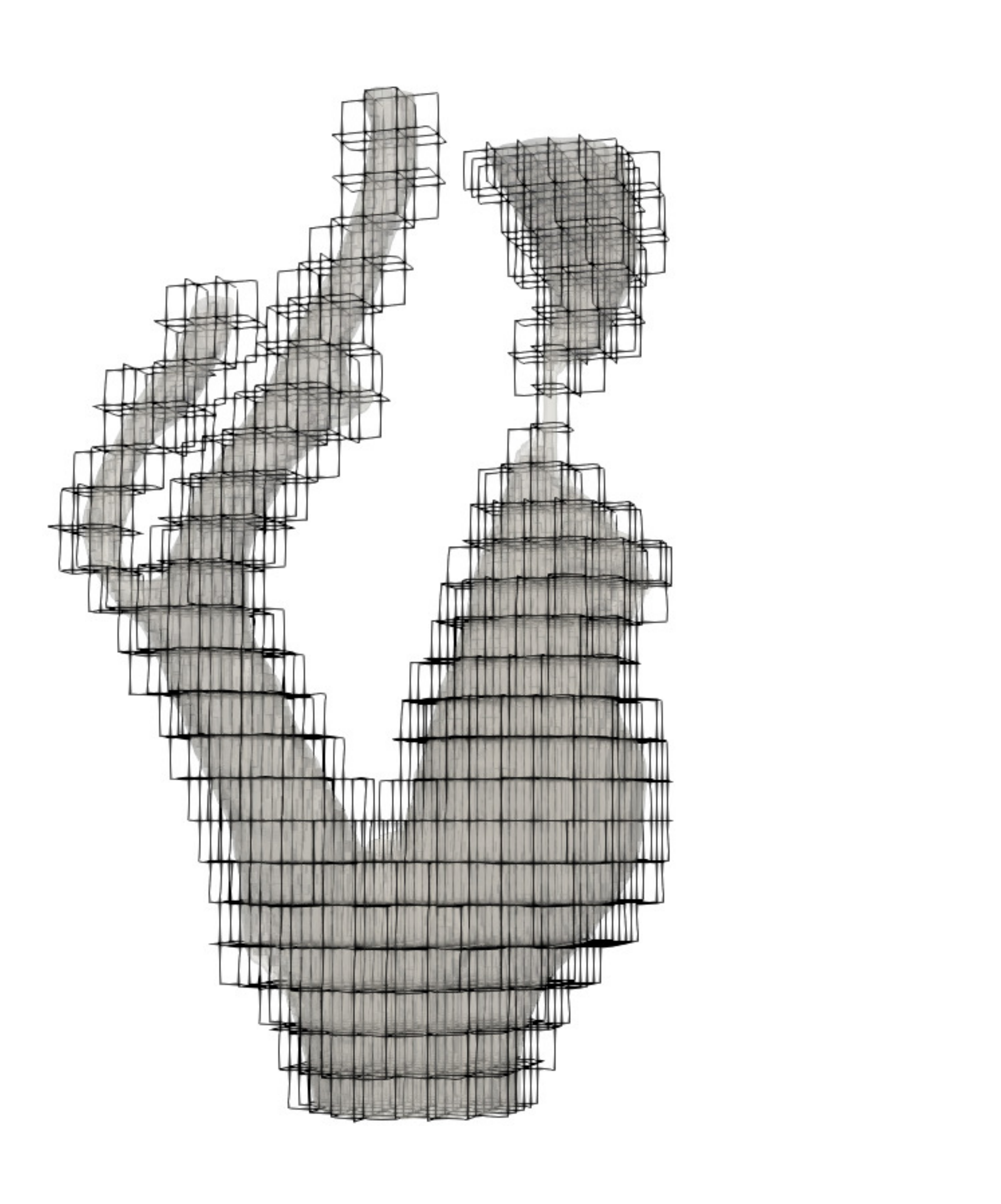}
		\caption{Initial mesh with $3158$ \#DOFs}
		\label{fig:3dstokes_mesh_step0}
	\end{subfigure}\hspace{0.08\textwidth}%
	\begin{subfigure}[b]{0.45\textwidth}
		\centering
   		\includegraphics[width=\textwidth]{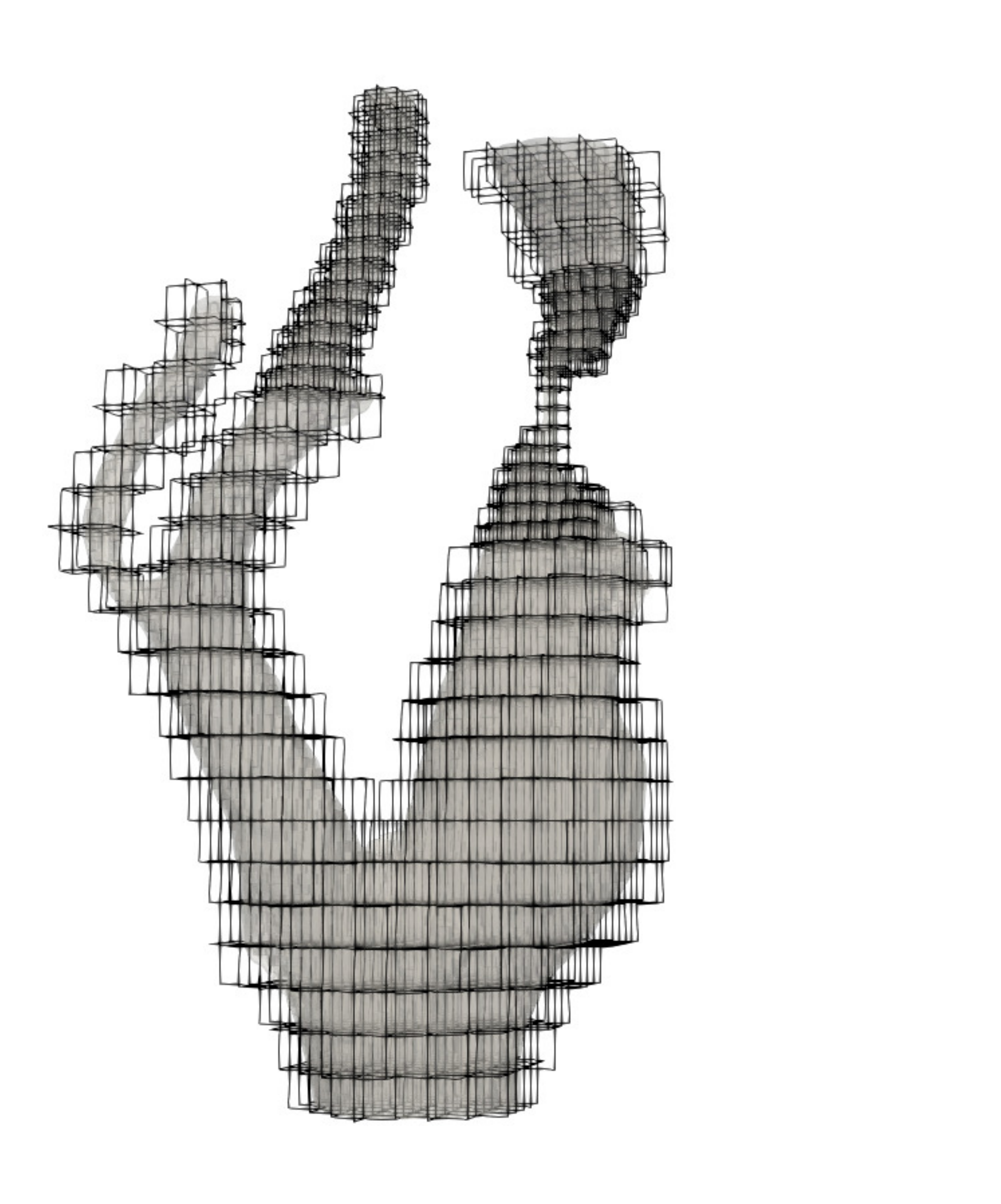}
		\caption{Step $1$ with $4161$ \#DOFs}
		\label{fig:3dstokes_mesh_step1}
	\end{subfigure}\\[12pt]
	\begin{subfigure}[b]{0.45\textwidth}
		\centering
   		\includegraphics[width=\textwidth]{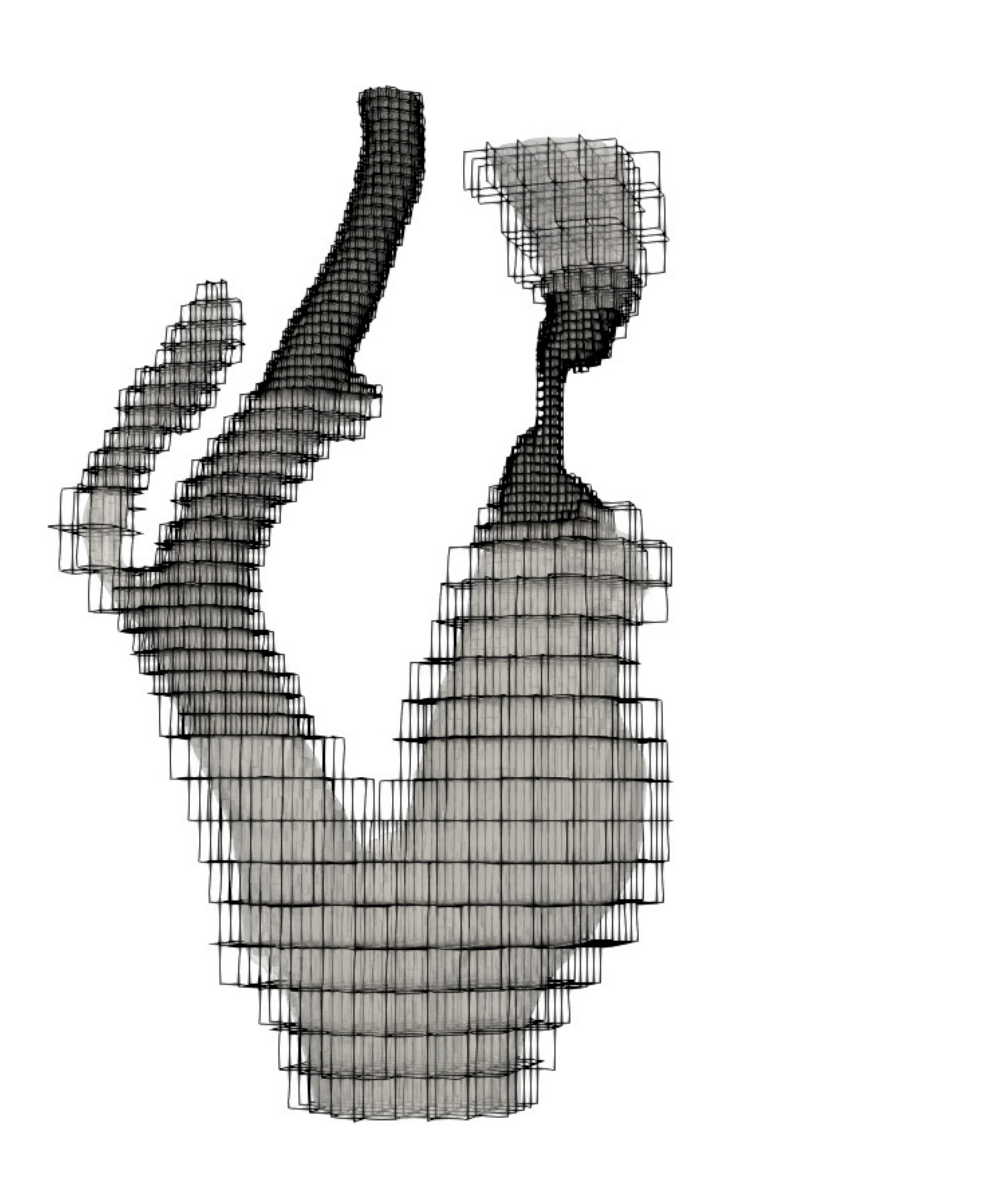}
		\caption{Step $2$ with $8784$ \#DOFs}
		\label{fig:3dstokes_mesh_step2}
	\end{subfigure}\hspace{0.08\textwidth}%
	\begin{subfigure}[b]{0.45\textwidth}
		\centering
   		\includegraphics[width=\textwidth]{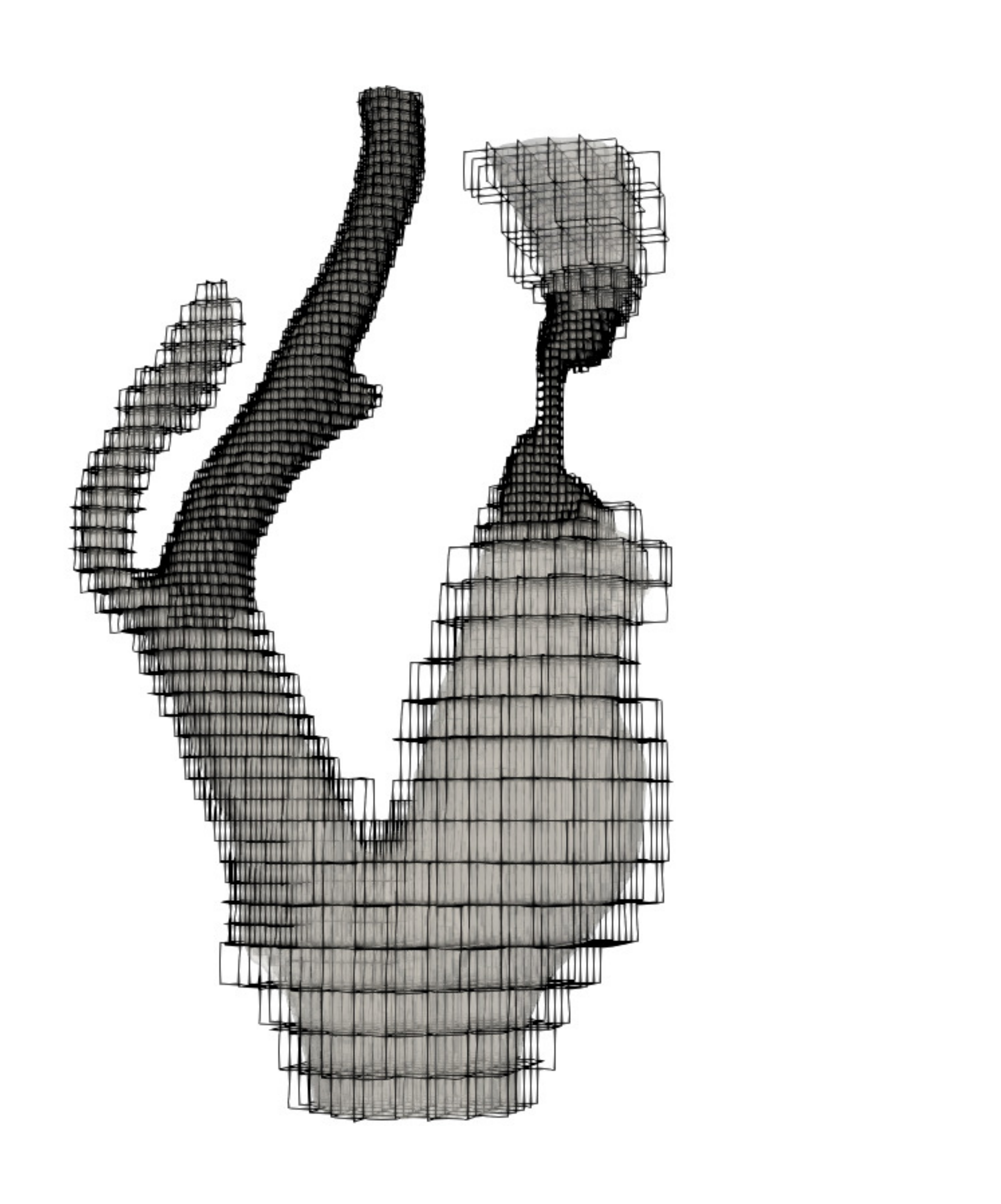}
		\caption{Step $3$ with  $12467$ \#DOFs}
		\label{fig:3dstokes_mesh_step3}
	\end{subfigure}
	\caption{Evolution of the mesh during the adaptive refinement process for the patient-specific viscous flow problem.}
	\label{fig:3dstokes_mesh}
\end{figure}

Different steps in the adaptive refinement procedure are illustrated in Figures~\ref{fig:3dstokes_velocity} and \ref{fig:3dstokes_mesh}. In all the refinement steps, the refinement strategy starts to focus on the regions where the errors are largest, \emph{i.e.}, near the stenosed section (\emph{i.e.}, the narrow region at the right artery) and at the outflow section of the left artery. Under local mesh refinement, the procedure resolves prominent solution details, most importantly the velocity field in the left artery and near the stenotic part of the right artery.

\begin{figure}
	\centering
	\includegraphics[width=0.6\textwidth]{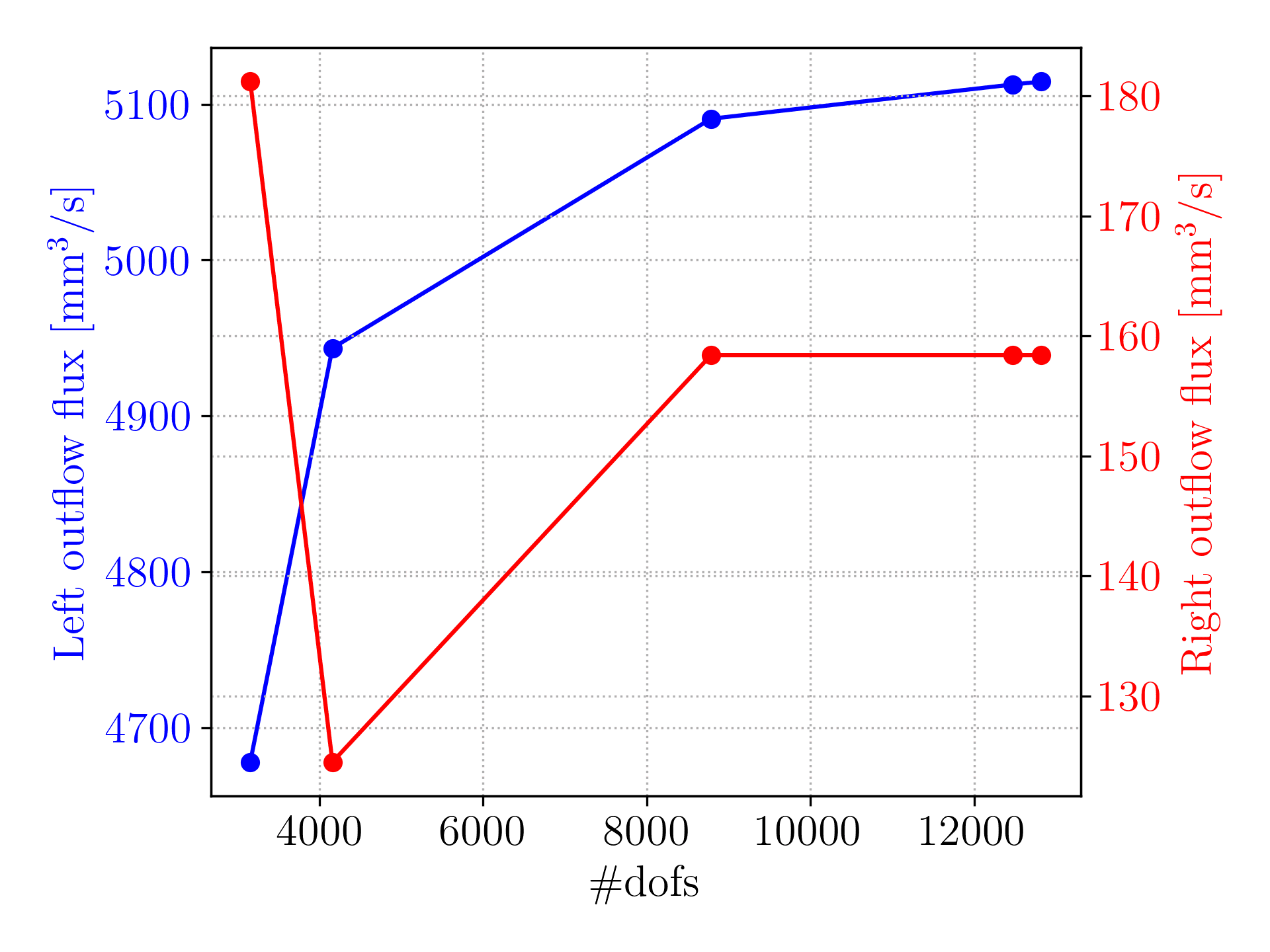}
	\caption{Mesh convergence of the flux at the left and right outflow boundary using adaptive mesh refinements for the patient-specific viscous flow problem.}
	\label{fig:3dstokes_conv}
\end{figure}

The flux at the outlet of the arteries is shown in Figure~\ref{fig:3dstokes_conv}, which is computed with the velocity field obtained by solving the flow problem using adaptive refinements. The solution of the flux in the left artery is observed to gradually converge toward a value of just over $5100$ [mm$^3$/s]. For the right artery, the maximum refinement depth is reached after the second refinement step. As a result, the flux in the right artery does then not substantially change anymore. At this point, the element sizes in the vicinity of the stenotic artery are similar in size to the voxels. The error then becomes dominated by the geometry reconstruction procedure, which also explains why the observed flux in the right artery deviates from the uniform mesh results in Ref.~\cite{divi2021}, \emph{viz.} $\varrho_{\rm max} = 2$ instead of the presently applied $\varrho_{\rm max} = 3$. It is observed that the adaptive procedure terminates after $4$ refinement steps, because of reaching the maximum refinement level in all the elements tagged for refinement. At this point, the adaptive simulation uses $12,816$ DOFs, which is substantially lower than the number of DOFs required using uniform refinements \cite{divi2021}, which amounts to approximately $10^5$.

%% file: sections/conclusions.tex
\section{Concluding remarks}\label{sec:conclusions}
In the immersed (isogeometric) analysis framework, the geometry representation is decoupled from the discretization. This enables the consideration of spline basis functions on complex volumetric domains, for which boundary-fitting discretizations cannot easily be obtained. Moreover, the decoupling of the geometry and the discretization allows one to have a globally accurate representation of the geometry, but only to refine the mesh in places where the errors are large. Such local mesh refinements have the potential to provide a significant efficiency gain compared to uniform meshes. The adaptive simulation strategy proposed in this work automatically refines the elements in places that significantly contribute to the error in the energy norm.

The developed error estimation and adaptivity strategy is based on residual-based error estimation, which is well-established in traditional finite elements and has been successfully applied in boundary-fitting isogeometric analysis. In the considered immersed setting, the residual-based error estimation and adaptivity framework requires the incorporation of the stabilization terms for the weakly imposed Dirichlet boundary conditions, and, in the case of the (mixed) Stokes flow problem, for the treatment of equal-order discretizations of the velocity-pressure pair. Adequate scaling of the stabilization constants with the mesh size is essential for the adaptive procedure to be effective. In particular the order dependence of the stabilization constants and the definition of the local element sizes must be treated adequately.

In contrast to residual-based error estimation for boundary-fitting finite elements and isogeometric analysis, in the stabilized immersed setting it is not evident that the residual-based error estimator bounds the error in the energy norm from above. This is a consequence of the absence of an $h$-independent weak formulation. In this work, it is reasoned, however, that under the assumption of sufficient smoothness, the residual is expected to be useful in the setting of an adaptive refinement strategy. For all numerical simulations considered, including simulations with reduced regularity, it is observed that the error estimator does provide an upper bound to the error in the energy norm. A rigorous study regarding the relation between the residual and the actual error is warranted.

It is demonstrated that the developed adaptive simulation strategy is particularly useful in a scan-based analysis setting, where manual selection of refinement zones is impractical. When used in combination with advanced image segmentation procedures to obtain a smooth geometry representation while preserving small geometric features, the developed adaptive refinement strategy optimally leverages the advantageous approximation properties of splines for geometrically and topologically complex domains. The adaptivity strategy results in a simulation workflow that is capable of obtaining  reliable, error-controlled, results with limited user interaction.

The developed adaptive solution strategy is elaborated for the Laplace problem and the Stokes problem. For other problems, such as, for example,  Navier-Stokes or Cahn-Hilliard problems, the starting point of the derivation of the error-estimator remains the same. The estimators are problem-specific, however, and hence need to be elaborated for such problems. The same holds for the consideration of additional or alternative stabilization techniques, specifically when these alter the Galerkin form of the problem.

%% file: sections/acknowledgement.tex
\section*{Acknowledgement}
We acknowledge the support from the European Commission EACEA Agency, Framework Partnership Agreement 2013-0043 Erasmus Mundus Action 1b, as a part of the EM Joint Doctorate Simulation in Engineering and Entrepreneurship Development (SEED). All the simulations in this work were performed based on the open source software package Nutils (www.nutils.org) \cite{nutils}. We acknowledge the support of the Nutils team. We would like to acknowledge fruitful discussions with Mats G. Larson regarding inf-sup stability and error analysis for immersed isogeometric approximations.